\documentclass[12pt]{article}
\usepackage{amssymb}
\begin{document}
\setcounter{page}{1}
\vspace*{27mm}

\begin{center}
{\LARGE \bf FILLINGS METHOD IN NUMBER \\
\vspace*{3mm}

THEORY}
\vspace*{10mm}

\normalsize

{\large \bf Mikhail V. Antipov}
\vspace*{7mm}

{\small Institute of Computational Mathematics and
Mathematical Geophysics \\
Siberian Branch, Russian Academy of Sciences, Novosibirsk}
\end{center}
\vspace*{7mm}

{\small \bf Abstract.} {\small In offered work the series of problems
of an analytical number theory is surveyed. This problems have direct
and defining reflection in prime numbers distribution. For their
solving the new method called as a method of fillings was developed. In
basis of its the functional imaging lays of the integrated characteristics
and elements of all period of generators - grids on an interval, and
the systems of generators appear as the main object of a research.
A fillings method separates from a traditional sieving process already on
an initial stage because of this fact.

It is necessary to refer possibility of registration and proof of the
important state about distribution of the majorizing characteristics
of numerical grids on an axe to one of main advantages of a method of
fillings. This state is called by the main theorem. Modifications
of a relation of the main theorem which is sequentially improved and
magnified also were by that source, from which one the concrete
outputs about distribution of prime numbers were obtained.

The quite particular versatility of the fillings method has
allowed to receive outcomes in address to such different tasks as
the evaluation of greatest distance between prime numbers, the
Goldbach's conjecture for even, the increase of the Goldbach's
representations, distribution of twins, asymptotics of
distribution of typical configurations primes and some others.
Thus there are all basis to guess, that possibilities of a
method are not exhausted by considered applications, as
the method hold sway over not only primes system.}
\vspace*{5mm}

\vspace*{5mm}
\newpage

\begin{center}
{\Large \bf 1. Introduction}
\end{center}
\normalsize

The concept of an integer with the complete right should be
named one of the first abstract object scientific cognition of
the world, a tool of which has become the mathematics. The
representation about prime number, allocated by distinctive and
well by known features has occurred also in antiquity. The
further researches have delivered problems of prime numbers
distribution, but the progress in the sanction some of them was
scheduled in comparison recently.

The successes in an essential degree have appeared connected
with development of the sieving process, at one time offered by
Eratosthenes. At the same time the series of problems, despite
of elementary mathematical statements, has not yielded all
efforts and progress to new results was essentially slowed. To
the present time is not new directions, leading though
presumably to the for a long time scheduled purposes.

The reason of a created situation can be only one: the sieving
process as any method has restricted capacity, adequate to
limiting level of achievement and conclusions. And this
potential, as far as demonstrate real development, is
practically exhausted. The absence of fresh ideas in so before a
fruitful direction of the number theory also testifies to a
limit of research opportunities of the profound antiquity
method.

The sieving process has a series essential and in essence
ineradicable defects, among which very weak interdependence
between base (prime) numbers, forming a sieve. Some
characteristics of distribution of primes in a sequence natural
is sharp depend from a set of forming values. Therefore there
was a idea to study them on classes of sets, being not limited
only primes. This idea has resulted to creation of a quite new
method of the number theory -- fillings method [1,\,2].

Beginning of the mathematics development as integer arithmetics
is impossible to separate from a problem of divisibility.
Probably, definition of prime number was known long before
Euclid, but only he proved by the simple and elegant theorem the
Cantorian equipotency of countable sets of natural and prime
numbers. The representation of actual infinity, necessary for
realization of this fact, was produced by Phales school
considerably earlier.

The concepts of countability, potential and actual infinity,
series of methods and operations over finite and unlimited
numerical objects, though and were received in half-intuitional
level, have nevertheless found reasonably qualified base for
use. Boiling joy of seizing by infinity axiom and unrestrained
flight of the idea have not allowed to pay attention for a
computability problem, but this permanent sin there were at soul
of the mathematics (and not only mathematics) and now.

Ancient Greek approach to registration and definition of number
- even integer and infinite, contained in core the ineradicable
factor of absolute unrestriction. Naturally, such significant
fact has an effect in speed creation of numerical objects of
actual infinity. At all development the mathematics (and the
science as a whole) have remained at the achieved stage. In
particular, one of examples of a extreme degree unreasonable
idealization of a theory has acted Cantorian theory of sets.

Considerable achievement of a modern science including the
mathematics can not, nevertheless, to hide sorrowful facts of
deep failures, hopeless deadlocks and insoluble contradictions.
Their initial source is the only infinity axiom, transformed a
zone of knowledge to field of an authority of an infinity
quantifier. The consequences of this phenomenon reach much
further problems only sciences.

However all offered development are oriented to the conventional
schemes [3] and will be conducted without references to
non-foundation, indefiniteness and unprovableness of infinity
axiom. Though the series determining and known rules, in
particular connected with concept of set of all primes, does not
make sense outside of mentioned axiom at all its numerous
defects. The ideas, connected with determining, central and
paramount role of infinity axiom, are considered and carefully
justified in work [4--6].

The Euclidean theorem about infinity of primes set can be
considered as the first research, the provable statement of a
modern stage of consciousness - stage of infinity axiom. In this
sense the theorem has acted as the first-born of a present
condition of the science. In a maximum degree is fair, that
those there was the theorem about numbers. The modern stage of
development and level cognition was determined earlier as time
by creation and adaptation of unrestriction idea.

A following major studied stage of primes distribution was the
sieve of Eratosthenes, realizing finite algorithm of reception
all primes, smaller \ $p_n^2$\,, \ if series primes down to
$p_{n-1}$ \ inclusive is known. It sieve will hereinafter appear
great-parent of updating and basis for theoretical research of
the various characteristics and features of primes distribution
accommodation in the natural sequence [3].

The sieving process as well as any method is exhausted. This
cognition law does not know exceptions. The truth, it does not
mean, that moment is (can be precisely established), when the
reception of new results will become impossible. Is not present,
exhaustion has an effect in sharp fall of a level of successes,
their quantity and significance. But any results are possible
nevertheless. Such process, fundamentally appropriate to
reality, is named by a real limit \ $(realimit)$ \ and is
investigated in works [4,\,5].

The offered research essentially uses the theorem about primes.
Namely, there is in such known formulation: Function \ $\pi (x)$ \
(quantity primes, not superior \ $x$\,) is represented by an
expression
$$
\pi (x) \ : \ \frac{x}{\ln x} \ = \ C(x)\,, \qquad p_i \ \le \ x
\qquad and \qquad C(x) \ \ \mathop{\Longrightarrow}\limits^{x \to
\infty} \ \ 1\,, \eqno{(1.1)}
$$
and the borders \ $C(x)$ \ are estimated, limited and
established. It is thus allowable to approximate the value of
prime \ $p_n$ \ depending from its index. In turn, it permits to
evaluate some important integrated characteristics of the finite
sequence \ $\{p_i\}_0^n$\,, \ which pursuant to (1.1) will be
specified at \ $n \to \infty$\,.

All other situations lean on a reasonably known theory or here
provable schemes. The elementary character of constructions in
determined degree is compensated by novelty of a method, on
principle distinguished from the sieving process -- and
algorithmically, and main directive purposes, and initial object
set.

It is especially necessary to note, that study and the decision
of series of problems of primes distribution had by the
beginning not given set primes, as for the majority of
researches such, but set of natural numbers. The way of research
of properties primes through study of natural numbers is not
only proven, but also is structurally used with the help of
fillings method. The expansion of main object set has basic
character and largely predetermines strength of the method.
\vspace*{5mm}

\begin{center} {\Large \bf 2. Some definitions} \end{center}
\vspace*{3mm}

Before to begin a exposition of a method, oriented on research
of some specific properties of integers correlation, it is
necessary give a series of initial definitions. Them obligation
is called by absence of adequate analogues, on which it was
possible to refer. If to take into account the general final
purpose -- distribution of primes, the similar statement appears
and characteristic while a unknown method.

Postulated integrity of initial values and elementary character
of made operations permit to declare, that main relations do not
turn out beyond the borders of fraction-rational numbers.  Even
the concept of numerical infinity at desire can be replaced by
estimable finite values, but the traditions and requirements of
laconic applications forbid to do it. Unique used function,
structurally including representation of infinity is function \
$\pi(x)$ \ from the central primes theorem (1.1).
\vspace*{3mm}

Idea of the method, number of the proof schemes and majority of
the outputs was obtained in 1983 and also the thought arose to use
period's characteristics of some finite numerical forms.
\vspace*{3mm}

A set \ ${\cal L} \ = \ \{\,L_{(s)}\,\}$ \ of
regulated infinite sequences \ $L_{(s)}$ \ of Boolean elements
\ $l_i \ \in \ \{0;\,1\}\,, \ \ i \in {\bf Z}$ \ represents
itself as the primary source object:
$$
{\cal L}\,: \qquad \quad \{\,L_{(s)} \ = \
...,l_{-m},l_{-m+1},...,l_{-1},l_0,l_1,...,l_n,l_{n+1},...
\,\}\,.
$$

Ordering of elements in an assigned sequence \ $L_{(s)}$ \ is
possible to fix by correspondence of each elements
\ $l_i$ \ to an integer \ $i + k_s \in {\bf Z}$\,, \ where \ $| k_s
| < \infty $ \ is some arbitrary constant. Distance \
$\rho_{ij}$ \ is determined for all elements \ $l_i,\,l_j \ \in \
L_{(s)}$\,:
$$
-\infty \ < \ i,\,j \ < \ \infty\,, \quad \qquad \rho_{(i+k)(j+k)} \
\equiv \ \rho_{ij} \ \equiv \
\rho\,(l_i,\,l_j) \ = \ | \ i - j \ |\,.
$$

It does not depend from \ $k_s\,, \ k$ \ or \ $s < \infty$\,. \ Integer
values of distance \ $\rho_{ij}$ \ do not conflict with rules of metric
space. The rule of the triangle has a feature. Greatest distance is
always equal to sum of others: \ $\rho_{ij} = \rho_{ik} + \rho_{kj}$\,.
\vspace*{3mm}

Thus it takes place determined isomorphism of sequence \
$L_{(s)}$ \ and vector \ $\{\rho_{ii'}\}$\,, \ $i \in {\bf Z}$ \
of distances between following one after another next units of
sequence \ $L_{(s)}$\,. \ If \ $l_i = 1$ \ for \ $-\infty < i <
\infty$\,, \ that \ $i' = \min\limits_{k > i} \{l_k = 1\}$\,. \
From here follows, that exact restoration of sequence \
$L_{(s)}$ \ by a vector \ $\{\rho_{ ii'}\}\,$ \ (except for
constant \ $k_s$) is possible. Certainly, with the same success
it is possible to generate and to consider a similar vector of
distances between zeroes, instead units.
\vspace*{3mm}

Now all ready for the direct introduction of definitions of
objects, compulsory and constantly used in further constructions
and researches.
\vspace*{3mm}

{\bf Definition 1.} {\sl Grid} \ $S(a)$ \ {\sl of the module} \ $a \ge 0\,,
\ \ a \in {\bf N}$ \ or \ $a$-{\sl grid} is named as the sequence
\ $L_a \ = \ \{ l_i \}_a$ \ of periodic elements:
$$
L_a: \qquad l_{i+a} = l_{i}\,; \ \ if \ \ l_{j} = 0\,, \ \ that \ \
l_{j+1} = l_{j+2} = ... = l_{j+a-1} = 1\,.
$$
\vspace*{0mm}

The following relations are valid for any grid as periodic object:
$$
\forall\,i \in {\bf Z}: \quad |\,l_{i+a} - l_{i}\,| =
0\,; \ \ if \ \ l_j = 0\,, \ \ that \ \exists \ \ l_k = 0\,: \quad
\rho\,(l_j,\,l_k) = a\,.
$$

The sequence from units only (0-grid) is marked \ $L_0$\,, \
and sequence from zeroes only (1-grid), is marked as \ $L_1$\,. \
The grid \ $L_{\infty}$ \ has equally one zero for
integer axis \ ${\bf Z}$ \ of indexes. The zeroes and units
alternate in first nontrivial 2-grid: \ $l_i = 1 - l_{i-1}$\,.
\vspace*{3mm}

{\bf Definition 2.} The \ $n$-{\sl grating} \ $V_n$ \ is named {\sl
product} of \ $n$ \ grids (optionally different) \
$S(a_i)\,, \ i = 1,2,...,n$\,. \ It is made under the
recurrent scheme
$$
V_0 = L_0\,, \quad V_i = V_{i-1} \& S(a_i)\,, \quad l_j (V_i) =
l_j(V_{i-1}) \& l_j(L_{a_i})\,, \quad j \in {\bf Z}.
$$
\vspace*{0mm}

The itemized logical multiplying of sequences is meant product of
grids. We shall give as an example 4-grating for grids \
$S(4;\,6;\,12;\,12)$\,:
$$
\begin{array}{lc}
V_0 = L_0\,: \qquad & ... \ 1111111111111111111111111 \ ... \\
V_1 = L_4 = S(4)\,: \qquad & ... \ 1011101110111011101110111 \ ... \\
V_2 = V_1 \& S(6)\,: \qquad & ... \ 1001101100111001101100111 \ ... \\
V_3 = V_2 \& S(12)\,: \qquad & ... \ 1000101100111000101100111 \ ... \\
V_4 = V_3 \& S(12)\,: \qquad & ... \ 1000001100111000001100111 \ ...
\end{array}
$$

Commutability, associativity, symmetry and transitivity of grids
products are the direct corollary of these properties for Boolean
elements \ $l_j$\,.
\vspace*{3mm}

{\bf Definition 3.} The \ $n$-{\sl filling} \ $Z_n$ \ is such \
$n$-grating \ $V_n\,, \ Z_n \subseteq V_n$\,, \ in which one at each
index \ $i\,, \quad 1 \leq i \leq n$ \ is executed
$$
V_n \ = \ V_{n-1} (a_1,...,a_{i-1},a_{i+1},...,a_n) \ \& \ S(a_i) \ \neq \
V_{n-1} (a_1,...,a_{i-1},a_{i+1},...,a_n)\,.
$$
\vspace*{-1mm}

The exception of any grid in \ $n$-filling carries on to the
permutation of zeroes and units. It takes place not always
for \ $n$-grating.
\vspace*{3mm}

{\bf Definition 4.} If maximum module \ $\max\limits_{1 \leq i \leq n}
a_i < \infty$\,, \ that anyone \ $n$-filling \ $Z_n$ \ \ (as \ well \ as \
$n$-grating) is periodic, and length of its period \ $PZ_n$ \
does not exceed the least common multiple of modules:
$$
\{\,Z_n\,(a_1,...,a_n)\,; \ V_n\,(a_1,...,a_n)\,\}: \qquad PZ_n
\ \le \ ((a_1,...,a_n))\,.
$$
\vspace*{-1mm}

{ \bf Definition 5.} The {\sl system} of grids \ $SS$ \ is unbounded
sequence of grids \ $S(a_i)$ \ of nondecreasing modules selected
according to some law, if at each \ $n$ \ the grids
\ $\{ S(a_i) \}_1^n\,, \ 1 \le i \le n$ \ form \ $n$-filling \ $Z_n$\,.
$$
SS\,: \qquad n \Rightarrow n + 1\,, \quad a_n \leq a_{n+1}\,,
\quad Z_n \Rightarrow Z_n \& S(a_{n+1}) = Z_{n+1}\,.
$$
\vspace*{-1mm}

The transition from parameter \ $n$ \ to \ $n + 1$ \ gives increase of set
of nonconterminous fillings \ $Z_{n+1}$ \ in product of prior \
$n$-filling \ $Z_n$ \ with the grid \ $S(a_{n+1})$\,. \ A possible
incongruity of fillings \ $Z_{n+1}(k_{n+1})$ \ is determined by
change of the shift parameter \ $k_{n+1}$ \ of new grid.
\vspace*{3mm}

{\bf Definition 6.} The \ $q$-{\sl series} of zeroes is an
interval \ $SR_n(q)$ \
of \ $n$-filling \ $Z_n$\,, \ restricted by units and containing
(including) \ $q$ \ units, where integer \ $q \geq 0$\,. \ {\sl
Length} of \ $q$-series \ $sr_n (q)$ \ is distance between
initial \ $1^{(0)}$ \ and final \ $1^{(q+1)}$ \ units of \ $q$-series:
\ $sr_n (q) = \rho (l_i,\,l_k)\,; \quad l_i = 1^{(0)}\,, \ l_k = 1^{(q+1)}$\,.
\vspace*{3mm}

The irreversible discrepancy between sieving process in all its
modifications and fillings method begins with introduction of
concept of zeroes series. Series of zeroes \ $SR_n(q)$ \
is the object, unknown for sieving process, where the zeroes and
units (that is eliminated and not eliminated numbers), are
strongly connected to the value of concrete natural number.
Then zeroes and units do not differ as supplemental, passing
properties of divisibility of this number. Vice-versa, the role
of zeroes becomes defining in the fillings method.
\vspace*{3mm}

{\bf Definition 7.} The {\sl maximum series} of zeroes \ $MSR_n(q) =
MSR_n(Z_n,\,q)$ \ of the value \ $0 < msr_n(q) < \infty$ \ is such \
$q$-series \ $SR_n(q)$\,, \ length which \ $sr_n (q)$ \ one greatest
in the set of \ $n$-fillings \ $\{Z_n\}$\,:
$$
msr_n (q,\,a_1,\,a_2,\,...,\,a_n) \ = \ \sup_{\{Z_n\}} \ \max_{k,\,j}
\ \{\,sr_n (q)\,[l_k,\,l_{k+j} = 1] \ \in \ Z_n\,\}\,.
$$
\vspace*{-1mm}

The maximum series of zeroes \ $MSR_n(q)$ \ and their value for \ $q =
0$ \ (without units inside series) are most indispensable for
theoretical constructions. These major characteristics of
fillings are marked accordingly \ $MSR_n$ \ and \ $msr_n$\,.
\vspace*{3mm}

{\bf Definition 8.} The \ {\sl regulated} \ $n$-filling \
$ZU_n$ \ = \ $ZU_n (a_1,\,a_2,\,...\,,\,a_n)$ \ distinguishes
algorithm of product of grids: after arbitrary fixing of
beginning of series, the zero of each next grid \ $S(a_i)\,, \ i
=1,...,n$ \ is multiplied with first right unit of the series.
The {\sl unregulated} \ $n$-filling \ $ZN_n $ \ combines all set of
possible fillings.

The {\sl semi-regulated} \ $n$-filling \
$ZP_n(a_1,\,...,\,a_r;\,a_{r+1},\,...,\,a_n)\,,
\ n > 3$\,,

\hspace*{-6mm} $1 < r < n-1 $ \ is an unregulated filling \
$ZN_{n-r} (a_{r+1},\,...,\,a_n)$\,, \ constructed on
regulated filling \ $ZU_r (a_1,\,a_2,\,...,\,a_r)$\,.
\vspace*{3mm}

The values of the greatest series of zeroes obtained by the filling
\ $Z_n \ (ZU_n - ZN_n)$ \ are marked as \ $msr_n(q,\,Z_n) = msr_n (q)$\,. \
It is important characteristics of filling.
\vspace*{3mm}

{\bf Definition 9.} Zero of {\sl multiplicity} \ $k\,, \ (1 \leq k \leq
n)$ \ is an element of \ $n$-filling \ $Z_n$ \,, \ if it is
multiplying of \ $k$ \ zeroes of generating grids \
$S(a_i)$\,.
\vspace*{3mm}

We shall give by the way for example already considered
different fillings \ $Z_4$ \ with vectors of multiplicity of
zeroes
$$
\begin{array}{ccc}
...1011101110111011... & ...1011101110111011... & ...1011101110111011... \\
...1101111101111101... & ...1011111011111011... & ...1011111011111011... \\
...1110111111111110... & ...1110111111111110... & ...1011111111111011... \\
...1111011111111111... & ...1111111111101111... & ...1111111011111111... \\
= & = & = \\
...1000001100111000... & ...1010101010101010... & ...1011101010111011... \\
...v11111vv11vvv111... & ...v2v1v1v1v1v1v2v1... & ...v3vvv1v2v1vvv3vv...
\end{array}
$$

Two 4-fillings and one 4-grating are given here. The grids \
$S(4),\,S(6)$, $S(12)$, \ $S(12)$ \ with the shifts are in
each column, then the outcome of product is given, and the
vectors of multiplicity of zeroes are in the last
line. Unit (zero of the multiplicity of zero) is marked by
sign \ $v$\,.

The examples demonstrate possibility to construct filling with
zeroes only multiplicity 1 (first variant). The second case
testifies change of periodicity (from 12 to 2), but thus there is
zero of the multiplicity two. The filling with zeroes of
multiplicity three and four cannot be constructed for this set
of grids, but such it is possible for 4-grating (third variant
of product). Zero of multiplicity 3 supplies in the last line
for it zeroes of grids \ $S(4),\,S(6)$ \  and \ $S(12)$\,.
\vspace*{3mm}

{\bf Definition 10.} The system \ $SS = \{S(a_i)\}$ \ with modules \
$a_i$ \ as degrees of the same integer \ $d \geq
2\,, \ (a_i = d^{k_i}, \ k_i \geq 1)$\,, \ is the {\sl
degree}-system. Integer value \ $d$ \ is named as the {\sl
basis} of the degree-system \ $SS = SS_d$ \ and appropriate
fillings \ $Z_n$\,.
\vspace*{3mm}

{\bf Definition 11.} The system \ $SS$ \ is named as the system {\sl
without multiple zeroes} \ $SS = SS'$\,, \ if for anyone \ $n$ \ there are
fillings \ $Z_n$\,, \ in which one there will be no zeroes of the
multiplicity above than unit.
\vspace*{3mm}

{\bf Definition 12.} The system of fillings \ $SS=VP$ \ is named {\sl
coprime}, if for anyone \ $1 \leq i \neq j < \infty$ \ the
modules of these grids satisfy to the relation \ $(a_i\,,\,a_j) = 1$\,.
\vspace*{3mm}

{\bf Definition 13.} The nonsingular system of grids and fillings \ $SS$
\ is named as {\sl mixed} \ $(SS = SM)$\,, \ if it does
not belong to any circumscribed types.
\vspace*{3mm}

{\bf Definition 14.} The system \ $VP = \{2,3,5,...,p_i,...\}$ \ of modules
as primes is named {\sl 0-prime system} \ $SP_0$\,. \ The system
of the pair primes \ $\{S(2),\,2 \cdot S(p_i)\}$: \ $SS = SM =
\{2,3,3,5,5,...,p_i,p_i,...\}$ \ is {\sl 0-double} system of primes \
$SW_0$\,.
\vspace*{3mm}

{\bf Definition 15.} Some non-singular system of grids \ $SS =
\{a_1,a_2,...,a_i,...\}$ \ belongs to {\sl first \ $SS_{(I)}$\,, \ second}
\ $SS_{(II)}$ \ or {\sl third type} \ $SS_{(III)}$\,, \ if accordingly
$$
\left\{ C_{ss} \ = \ \lim\limits_{n \to \infty} \sum_{i=1}^n\,
\frac{1}{a_i} \right\}: \quad \qquad C_{ss} \ \leq \ 1\,, \qquad
1 < C_{ss} < \infty\,, \qquad C_{ss} \ = \ \infty.
$$
\vspace*{1mm}

For example, the sieve of Eratosthenes on each step
is the regulated filling of the system \ $SP_0$ \ for an interval
\ $I \ \ll \ PZ_n$\,, \  if the first zeroes of each grid after
termination of algorithm are exchanged by units.
\vspace*{3mm}

But before to address to the main representatives of systems of
the third type \ $SS_{(III)}$\,, \ it is necessary to consider systems
without multiple zeroes \ $SS'$ \ (all zeroes have multiplicity 1).
And already then to formulate and to present
the main statement, in which as evidently as numerically the
idea of a fillings method is concentrated. The fact is that
without exact formulas of model fillings the transition to more
complex and important systems is unreasonable.
\vspace*{5mm}

\begin{center} {\Large \bf 3. Systems without multiple zeroes}
\end{center}
\vspace*{3mm}

Systems without multiple zeroes \ $SS'$ \ are systems, for which at each \
$n$ \ is present \ $n$-filling \ $Z_n$\,, \ and in it on a
period, hence, on axis there is no zero of multiplicity higher
unit. Such systems as model for other more important and
necessary systems, nevertheless, have beside features. Just they
permit to find their characteristics. For example, the
construction of maximum series for fillings of such systems is
determined only by ordering of grid modules.

For intervals \ $J_n$ \ about a period and more density for them
perfectly is approximated by value \ $(\gamma_n)$\,. \ Decrease
of the interval leads to growth of zeroes density because of
fall them multiplicity down to \ $(\gamma_n^*)$\,, \ which is
limiting. But for systems without multiple zeroes \ $SS$ \ there
is the only zeroes frequency, it is density of zeroes \
$\gamma_n = (\gamma_n^*) $ \ for the period and consequently
axis for given filling \ $Z_n \subset SS$\,. \ To the first and
the most indicative class of such systems concern degree-systems
\ $SS_d$\,.

The zeroes frequency (density) for the period and whole numerical axis
of any \ $n$-filling \ $Z_n \subset SS_d$ \ is equal
$$
SS_d: \qquad \gamma_n \ = \ 1 \ - \ \frac{E_n}{PZ_n} \ = \
\frac{H_n}{PZ_n} \ = \ \sum_{i=1}^n \ \frac{1}{a_i} \ =
\ \sum_{i=1}^n \ \frac{1}{d^{k_i}}\,. \eqno{(3.1)}
$$

As far as research is supposed only non-singular ($msr_n < \infty$ \ for
each \ $n$) systems, that
at all \ $n$ \ value \ $\gamma_n < 1$\,. \ Besides the condition
of sequentially non-decreasing modules of grids for a system
\ $SS_d$ \ testifies, that according to (3.1) we have obvious
relation \ $PZ_n \le a_n$\,. \ In such case expediently to
consider a number system with the base-radix \ $d$ \ and
representation of frequency \ $\gamma_n$ \ in this system.

As is known, number system with basis \ $d$
\ has all digits less \ $d: \ 0 \le m_j \le d - 1$\,. \
However there is one obstacle on a way of representation of
zeroes density \ $\gamma_n$ \ kind \ $d$-th fraction, which is
consisted in arbitrary quantity of grids of one module. In this
connection we shall allocate subsets from the class of
degree-systems.
\vspace*{3mm}

{\bf Definition 16.} The {\sl correct} system \ $SS = SS_d$
\ of degree grids \ $S(a_i)$ \ with the basis \ $d \ge 2$ \ is
named system, for which the quantity of grids of one module \ $a_i$
\ does not surpass \ $d - 1\,: \quad S(a_k) = S(a_{k+1}) = ... =
S(a_{k+s-1})\,, \quad s \le d-1$\,.
\vspace*{3mm}

{\bf Definition 17.} The filling \ $Z_n(a_1,..., a_n)
\subset SS$ \ is named as {\sl saturated}, if filling \ $Z_{n +
1}$ \ is singular, that is \ $msr_n \ = \ msr_n (q) \ = \ \infty$\,,
$$
Z_n (a_1,a_2,...,a_n) \ \& \ S(a_{n+1}') \ = Z_{n+1}
(a_1,a_2,...,a_n,a_{n+1}') \ = \ L_1  \eqno{(3.2)}
$$
for any grid \ $S(a_{n + 1}')$\,, \ let and not included into
the system \ $SS$\,, \ but satisfying relation for period \ $PZ_n$ \ of
\ $n$-filling
$$
( \ PZ_n\,, \ a_{n+1}' \ ) \ = \ a_{n+1}'\,, \qquad PZ_n \ \le \
((a_1,a_2,...,a_n))\,.  \eqno{(3.3)}
$$
\vspace*{1mm}

Here, as above as further, the expression in brackets of the
first equality (3.3) means the greatest common divisor, and in
the second equality -- least common multiple. From the
definition clearly, that the period of saturated filling for \
$n = k_s$ \ grids contains equally one unit.
\vspace*{3mm}

{\bf Definition 18.} {\sl Sparse} frequency \ $\gamma$ \ is
value, received by exception from initial frequency \ $\alpha$ \
frequency and \ $d$-th digits of saturated filling:
$$
\alpha = 0.d^*...d^*m_{r+1}m_{r+2}...\,, \quad d^* = d-1; \qquad
\gamma \ = \ 0.\,m_{r+1}\,m_{r+2}\,...\,m_v\,...~, \eqno{(3.4)}
$$
where because of correctness of a system value \ $m_{r + 1} < d - 1$\,.
\vspace*{3mm}

{\bf Definition 19.} The grid \ $S(a_i)\,, \ 1 \le i \le n$ \ of
filling \ $Z_n$ \ is named as {\sl essential}, if value of the
maximum series is \ $msr_n (0) \ge a_i$\,.
\vspace*{5mm}

{\bf Theorem 1.} Length \ $msr_{\alpha} \{n\}$ \ of the maximum
series \ $MSR_{\alpha} \{ n\}$ \ for degree \ $n$-filling \ $Z_n$ \
of degree-system \ $SS_d$ \ with the basis \ $d \ge 2$ \ at zeroes
frequency \ $\alpha$ \ is expressed by the formulas (here parameter
\ $q$ \ is equal \ $0$):
\vspace*{2mm}

{\bf A.} Quantity of grids \ $n \ \leq \ k_s \ = \ (d - 1)\,[-\log_d (1 -
\alpha)] \ = \ (d-1)\,r$\,. \ Then
$$
MSR_{\alpha} \{n\}: \qquad msr_{\alpha} \{n\} \ = \
d^{[\frac{n}{d - 1}]} \left\{ 1 + n
- (d - 1)\left[\frac{n}{d - 1}\right] \right\}\,. \eqno{(3.5)}
$$

{\bf B.} Quantity of grids \ $n \ > \ k_s \ = \ (d - 1\,[-\log_d (1 -
\alpha)] \ = \ (d - 1)\,r$\,. \ \ Then \ $msr_{\alpha} \{k_s\}
\ = \ d^r\,, \ k \ = \ n - k_s\,, \ \gamma \ = \ 1 - d^r\,(1 - \alpha)$\,,
$$
MSR_{\alpha} \{n\}: \qquad \qquad msr_{\alpha} \{n\} \ = \ d^r \cdot
msr_{\gamma} \{k\}\,. \qquad \eqno{(3.6)}
$$

{\bf C.} Value of the maximum series \ $MSR_{\gamma} \{k\}$ \ for
frequency \

\hspace*{-6mm} $\gamma \ = \ 1 - d^r\,(1 - \alpha) \ = \ 0.\,m_1\,m_2\,...
\,m_v\,m_{v+1}\,...$ \ is equal
$$
msr_{\gamma} \{k\} \ = \ \frac{k - t_1}{1 - \gamma_v} \ + \
\frac{t_1 - t_2}{1 - \gamma_{v-1}} \ + \ ... \ + \ \frac{t_{v-1} - t_v}{1 -
\gamma_1} \ + \ t_v \ + \ 1\,, \eqno{(3.7)}
$$
where
$$
t_j \ = \ \left( t_{j-1} - \sum_{i=1}^{s-j+1} \ m_i\right)
{\rm mod}\,(d^{v-j+1} - \gamma_{v-j+1}\ d^{v-j+1}) +
\sum_{i=1}^{v-j+1} \ m_i\,,
$$
$j \ = \ 1,\,2\,,...\,,\,v\,; \ \ t_0 \ = \ k\,; \ \ \gamma_j \ = \
\sum\limits_{i=1}^j \ m_i d^{-i}$ \ under condition
$$
d^v \ - \ \sum_{i=1}^v m_i (d^{v-i} - 1) - 1 \ < \ k \ \le \
d^{v+1} - \sum_{i=1}^{v+1} m_i (d^{v-i+1} - 1) - 1\,.
$$
\vspace*{3mm}

{\sl Proof} of the formulas we shall give consistently. Let
\ $n \leq k_s$ \ and frequency \ $\alpha$ \ is expressed as (3.4).
Then
$$
msr_{\alpha} \{n\} = \left\{
\begin{array}{ll}
d^r \ & for \ n = k_s = r\,(d - 1)\,, \  \\
d^{[\frac{n}{d - 1}]} \left\{ 1 + n
- (d - 1)\left[\frac{n}{d - 1}\right] \right\} & for \ n < r(d
- 1)\,, \end{array} \right.\eqno{(3.8)}
$$
that it is enough clearly from concrete appendices and
definition of regulated filling. For example, the second variant
(3.8) for \ $n = (r-1 )\,(d - 1)$ \ is reduced to first. From
here expression (3.5) and statement {\bf A} follows.
\vspace*{3mm}

Expression {\bf B} and relation (3.6) follow from a obvious
conclusion, that the availability of saturated grids increases a
maximum series that filling, but without saturated grids,
equally in \ $d^r$ \ time, that is increasing of series occurs by value
of the period of filling from saturated grids. The particular case of equality
(3.6) can be noticed in the second relation (3.8).
\vspace*{3mm}

The most compound, but also the major variant is submitted by
expression {\bf C}. It finally permits find exact value \
$msr_{\gamma} \{k\}$ \ of maximum series \ $MSR_n (q)$ \ at all
parameters of degree-systems \ $SS_d$\,.
\vspace*{3mm}

Let the basis of correct \ degree-system \ \ $SS_d\,, \ d \ge 2$ \ with
zeroes density \ $\gamma = 0.\,m_1\,m_2\, ... \,
m_v\,m_{v + 1}\,...$ \ is given. Then reception of the maximum
series \ $MSR_n$ \ of filling \ $Z_n \subset SS_d$ \ of length \
$msr_n = d^v$ \ will be required equally the quantity \ $n$
$$
\{MSR_n (0) \subset Z_n\}: \quad \qquad n \ = \ d^v \ - \ \sum_{i = 1}^v \
m_i\,(d^{v-i} \ - \ 1) \ - \ 1 \eqno{(3.9)}
$$
$d$-th grids. This number is minimum, but from them only
$K_1 = \sum_{i = 1}^v \ m_i$ \ grids are essential.
Really, quantity of units in the period \ $PZ = d^v$\,, \ formed by grids
\ $K_1 = \sum_{i = 1}^v \ m_i$ \ equally \ $d^v -
\sum_{i = 1}^v \ m_i\,d^{v-i}$\,. \ Product executes
following grids for the scheme of regulated filling and
since these grids will be inessential, that they
can be replaced by infinite, we shall receive value \
$msr (q) = d^v$\,, \ where \ $q$ \ by one unit (boundary)
less than units in the period of filling \ $Z (K_1)$\,. \
Thus, the common quantity of grids is \
$$
n \ = \ d^v \ - \ \sum_{i = 1}^v \ m_i\,d^{v-i} \ + \ \sum_{i =
1}^v \ m_i \ - \ 1\,, \
$$
that coincides with the formula (3.9). Such quantity of grids is
minimum, that it follows from definitions for correct systems \
$SS_d$ \ for fixed \ $d$\,.

Then we have, after designation of the maximum series length for given
frequency and grids quantity \ $msr_{\gamma} \{n\}$:
$$
msr_{\gamma} \left\{d^v - \sum\limits_{i = 1}^v
m_i\,(d^{v-i} - 1 )\right\} \ = \ d^v + msr_{\gamma}
\left\{\sum\limits_{i = 1}^v m_i \right\}. \eqno{(3.10)}
$$
It follows from the formula (3.9), \ whence we receive immediately the
equality \ $msr_{\gamma} \left\{d^v - \sum_{i = 1}^v m_i\,(d^{v-i} - 1) -
1\right\} = d^v$\,. \ Value and the maximum series, formed by \
$K_1 = \sum_{i = 1}^v \ m_i$ \ grids, are repeated for an
initial interval of the period \ $d^v$\,. \ Moreover it is
necessary to note, the grid \ $S(d^{v + 1})$ \ is inessential
for the considered interval, as maximum series \ $msr_{\gamma}
\left\{\sum_{i = 1}^v \ m_i \right\} < d^v$\,.

In conditions of relation (3.10) we have equality
$$
msr_{\gamma} \left\{d^v - \sum\limits_{i = 1}^v
m_i\,(d^{v-i} - 1 ) + k - 1\right\} \ = \ d^v + msr_{\gamma}
\left\{\sum\limits_{i = 1}^v m_i + k - 1\right\}, \eqno{(3.11)}
$$
where \ $1 \ \le \ k \ \le \ \left\{d^{v+1} - \sum\limits_{i = 1}^v
m_i\,(d^{v+1-i} - 1 ) - d^v + \sum\limits_{i = 1}^v
m_i\,(d^{v-i} - 1 )\right\} \ = $

$= \ (d - 1)\,\left\{\,d^v \ - \
\sum\limits_{i = 1}^v \ m_i\,d^{v-i}\,\right\}$\,.

Really, there are \ $d^{v+1} - d\,\sum_{i = 1}^v \
m_i\,d^{v-i}$ \ units in the interval by length \ $d^{v + 1} -
1$ \ as a result of product of \ $\sum_{i = 1}^v \ m_i$ \ grids,
as for this interval already grids of kind \ $S(d^{v + 1})$ \
are inessential (3.9). Thus, the formula (3.11) will be
valid so long as
$$
k \ + \ d^v \ - \ \sum\limits_{i = 1}^v m_i\,(d^{v-i} - 1 ) \
\le \ d^{v+1} \ - \ d\sum\limits_{i = 1}^v m_i\,d^{s-i} \ + \
\sum\limits_{i = 1}^v m_i\,,
$$
whence we receive the border for values \ $k$ \ from above. For
large \ $k$ \ it should in the expression (3.11) replace
parameter \ $v$ \ to \ $v + 1$\,. \ At \ $k  = 1$ \ we are
return to the formula (3.10), and at \ $k = 0$ \ (outside of
conditions) in a left-hand part (3.11) we receive grids quantity
(3.9) and length of series \ $msr_{\gamma} = d^v$\,.
\vspace*{1mm}

Let are given again the basis of correct degree filling \ $d$ \ and
density of zeroes \ $\gamma \ = \ 0.\,m_1\,m_2\,...\,m_v\,... \ $\,, \ but
already at \ $0 \le m_1 < d - 1$\,, \ that is for sparse. Then
$$
msr_{\gamma} \{k\} \ = \ \frac{k \ - \ t_1}{1 \ - \ \gamma_v} \ + \
msr_{\gamma} \{t_1\}\,, \eqno{(3.12)}
$$
where \ $t_1 \ = \ (k - \sum_{i=1}^v \ m_i)\,{\rm mod}\,(d^v -
\gamma_v d^v) + \sum_{i=1}^v \ m_i$\,, \ frequency \

$\gamma_v \ = \
0.\,m_1\,m_2\,...\,m_v \ = \ \sum_{i=1}^v \ m_i d^{-i}$ \ under
condition
$$
d^v \ - \ \sum_{i=1}^v m_i (d^{v-i} - 1) - 1 \ < \ k \ \le \
d^{v+1} - \sum_{i=1}^{v+1} m_i (d^{v-i+1} - 1) - 1\,, \eqno{(3.13)}
$$
and if \ $k \ \le \ d - 1$\,, \ then \ $msr_{\gamma} \{k\} \ = \
k + 1$\,.

From expression (3.11) it is possible to conclude, that when the
value \ $k$ \ lies in borders, specified by relations (3.10, \
3.11), the maximum series is equal
$$
msr_{\gamma} \{k\} = d^v + msr_{\gamma} \left\{k - d^v + \sum_{i=1}^v
m_i d^{v-i} \right\} = d^v + msr_{\gamma} \{k - d^v (1 - \gamma_v) \}.
\eqno{(3.14)}
$$

We apply consistently the formula (3.14) \ $j$ \ time so that
the value \ $k - j\,d^v\,(1 - \gamma_v)$ \ has not become less \
$\sum_{i=1}^v m_i$\,. \ At the same time \ $j$ \ should be
greatest of possible. The limits of change \ $k$ \ are
established in view of product of the first \ $\sum_{i=1}^v m_i$
\ grids. Under these conditions the value \ $t_1 = k -
j\,d^v\,(1 - \gamma_v)$ \ can be found only as
$$
t_1 \ = \ \left( k - \sum_{i=1}^v \ m_i \right)\,{\rm mod}\,(d^v -
\gamma_v d^v) \ + \ \sum_{i=1}^v \ m_i \ =
$$
$$
= \ \left\{\frac{k - \sum_{i=1}^v
\ m_i}{d^v(1 - \gamma_v)}\right\}d^v(1 - \gamma_v) \ +
\ \sum_{i=1}^v \ m_i \,, \eqno{(3.15)}
$$
where \ $\{\cdot\}$ \ means fractional part of function. From here
follows, that \ $j = \frac{k - t_1}{d^v ( 1 - \gamma_v)}$\,, \ and
then we shall receive
$$
msr_{\gamma} \{k\} \ = \ jd^v \ + \ msr_{\gamma} \{t_1\} \ = \
d^v\,\frac{k - t_1}{d^v (1 - \gamma_v)} \ + \ msr_{\gamma} \{t_1\}\,,
\eqno{(3.16)}
$$
that coincides with the statement (3.12). The necessity of
conditions (3.13) at search of value \ $t_1$ \ (3.15) is obvious, as
differently becomes impossible filling by \ $\sum_{i=1}^v m_i$ \ given
grids. Value \ $j$ \ in (3.16) is common quantity of equality
applications (3.14). The last condition (3.13) with the maximum
series are also obvious, as far as variant is submitted here,
when all grids are inessential.
\vspace*{1mm}

Now we apply expressions (3.12,\,3.13 ) recursively for
reception of maximum series \ $MSR_{\gamma} \{k\}$ \ value, that
is we are addressed to the formulas (3.12,\,3.16 ) at
first at greatest allowable \ $v$\,, \ determined by condition,
then at \ $v - 1,\,v - 2\,,...,\,2$\,, \ consistently finding
values \ $t_j$ \ from expression:
$$
t_j \ = \ \left( t_{j-1} - \sum_{i=1}^{v-j+1} \ m_i\right)
{\rm mod}\,(d^{v-j+1} - \gamma_{v-j+1}\ d^{v-j+1}) +
\sum_{i=1}^{v-j+1} \ m_i\,, \eqno{(3.17)}
$$
repeating condition of the theorem. In equality (3.17) the
parameters are in the borders \ $j \ = \ 1,\,2\,,...\,,\, v\,; \
\ t_0 \ = \ k\,; \ \ \gamma_j \ = \ \sum\limits_{i=1}^j \ m_i
d^{-i}$\,, \ and initial value \ $v$ \ is founded from condition
(3.13).

As value \ $t_v = (t_{v-1} - m_1)\,{\rm mod}\,(d - m_1) +
m_1 \ \le \ d - 1$\,, \ hence maximum series is \ \
$msr_{\gamma} \{t_v\} \ = \ t_v + 1$\,. \ The theorem is proven.
\hfill $\Box$
\vspace*{3mm}

As an example we shall consider the degree-system \ $SS_3$ \ with given
density \ \ $\gamma = 5/8 = 0.121212...$\,. \ We shall find \
$msr_{\gamma} \{16\}$\,. \ Thus \ $k$ \ value \ $v = 3$\,, \ as
with condition (3.13) we receive: \ $14 < 16 = k < 37$\,. \ If
to take into account, that from relation (3.17) and for \
$\gamma_j$ \ we have: \ $t_1 = 5\,, \ 1 - \gamma_3 = 11/27\,; \
t_2 = 5\,,$ \ $1 - \gamma_2 = 4/9\,; \ t_3 = 1\,, \ 1 - \gamma_1 =
2/3$\,, \ with the help of expression (3.7) we shall receive value of
maximum series \ $msr_{\gamma}\{16\} = 27 + 6 + 1 + 1 = 35$\,.

It should note, that in given statement the task of parameter \
$q$ \ inexpedient, as units automatically enter in value \ $k$ \
because from inessential grids of next filling. However it do
not without this parameter, equivalent to quantity of infinite
grids, at the task of zeroes density of filling in a kind of
finite fraction.

Thus, if density of zero \ $\gamma < 1$ \ and basis \ $d \ge 2$
\ of degree-system \ $SS_d$ \ are given, intervals of maximum
series \ $MSR_{\gamma} \{k\}$ \ always are determined precisely
with the help of all three relations of the theorem 1. Their
application does not assume any restrictions relatively included
saturated or essential grids in fillings.
\vspace*{3mm}

The theorem 1 permits to generate the important conclusions.
\vspace*{3mm}

{\bf Theorem 2.} The value of maximum series \ $msr_{\alpha} \{n
+ q\}$ \ with \ $q$ \ units allows unimprovable valuation in any
degree-system \ $SS_d$ \ with zeroes density \ $0 < \alpha < 1$
$$
MSR_{\alpha} \{n+q\} \subset SS_d: \quad \qquad msr_{\alpha} \{n+q\} \ < \
\frac{n + q}{1 - \alpha} \ + \ 1\,.  \qquad \eqno{(3.18)}
$$
\vspace*{0mm}

{\sl Proof} reasonably transparent follows from the formulas of the
theorem 1. Not too complex to show, that availability \ $k_s \ge d
- 1$ \ grids of saturation only eases the formulation of the theorem.
Therefore we shall consider case \ {\bf C.} \ For it
$$
msr_{\gamma} \{k\} \ = \ \frac{k - t_1}{1 - \gamma_v} \ + \
\frac{t_1 - t_2}{1 - \gamma_{v-1}} \ + \ ... \ + \ \frac{t_{v-1} - t_v}{1 -
\gamma_1} \ + \ t_v \ + \ 1\,,
$$
the sum of non-negative numerators of fractions is equal \
$k$\,, \ and denominators \ $1 - \gamma_{v-i} \ge 1 - \gamma$ \
for all \ $i$\,. \ Thus, the series value (3.7) is the closer to
valuation (3.18), the closer frequency \ $\gamma$ \ to zero.
From here follows unimprovable valuation. A final kind
expression (3.18) acquires after replacement \ $k$
\ in (3.7) to \ $n + q$\,. \hfill $\Box$
\vspace*{3mm}

The main merit of the theorem 1 consists in
important generalization for a class of systems.

Earlier all states and conclusions of this chapter
were formulated for class of degree-systems \ $SS_d$\,. \
At the same time rather easily to look after, that reasoning at
designing and algorithmization of constructions of maximum
series \ $MSR_n (q) $ \ for system without multiple zeroes
are analogous considered by the theorem 1. Similar though
naturally little more complex and difficult will be and
formula relations of the type (3.4 -- 3.7).
\vspace*{3mm}

{\bf Theorem 3.} The value of maximum series allows absolute
unimprovable majorant in any system \ $SS$ \ without multiple
zeroes with the density \ $0 < \alpha < 1$ \
$$
MSR_{\alpha} \{n+q\} \subset SS: \qquad \qquad msr_{\alpha}
\{n+q\} \ < \ \frac{n + q}{1 - \alpha} \ + \ 1\,.  \qquad
\eqno{(3.19)}
$$
\vspace*{3mm}

{\sl Proof.} A class of systems without multiple zeroes is
essentially wider of degree-class: \ $SS \supset SS_d$\,. \ For
example, to such class systems of grids concern:
$$
S(a_i) \ \subset \ SS: \qquad \quad (\,a_i,\,a_{i+1}\,) \ = \ a_i\,,
\qquad \forall \ i \ \ge \ 1\,. \eqno{(3.20)}
$$

At the same time not only the systems, satisfying to relation
(3.20), enter in such class. Besides such systems \ $SS$ \
always the first type, but not all systems of the first type are
systems without multiple zeroes.

Unimprovable valuation from below, indicated in such transparent
form (3.19), is reasonably clear, as far as \ $msr_{\alpha} \{n
+ q\} \le n + q + 1$\,, \ and the value \ $\alpha$ \ can be
near from zero. It means, that at all inessential grids
easily find a border of density \ $\alpha < \alpha_0$\,, \ for which
$$
msr_{\alpha} \{n + q\} \ = \ \left[\frac{n+q}{1 - \alpha}\right]
+ 1 \ = \ n + q + 1\,, \qquad [\,\cdot \,] \ - \ integer \ part.
$$

The relations {\bf A} and {\bf B} of theorem 1 for saturated
grids of with evidence are transferred for systems without
multiple zeroes. Therefore the special attention is deserved
case {\bf C} and formula (3.7). However easily to see, that
algorithmical features of constructions of maximum series in
systems without multiple zeroes and degree-systems coincide. It
means, that recurrent formulas for calculation of the maximum series
in systems \ $SS$ \ should be the type (3.7) and to differ only
reception of values \ $t_j$\,, \ which \ we \ shall \ designate \
$tt_j$\,.

In a result we shall receive transformed from a
relation (3.7) formula, in which the given zeroes density
\ $\alpha$ \ is consistently submitted approximations \
$\alpha_i$\,. \ They are similarly connected by inequalities
$$
\frac{1}{1 - \alpha} \ > \ \frac{1}{1 - \alpha_v} \ \ge \
\frac{1}{1 - \alpha_{v-1}} \ \ge \ ... \ \ge \ \frac{1}{1 -
\alpha_1}\,, \
$$
while corresponding non-negative values \ $tt_j$ \ give \
$$
n + q - tt_1 + tt_1 - tt_2 +... + tt_v \ = \ n + q\,.
$$

If given filling maintains \ $k_s$ \ saturated grids,
we shall act by analogy with the degree-filling.
After substitution of valuation for sparse density \
$\gamma$\,, \ taking into account a period of filling by
saturated grids, equal \ $D \ = \ \frac{1 - \gamma}{1 -
\alpha}$\,, \ we shall receive:
$$
msr_{\alpha} \{n\} \ = \ D \cdot msr_{\gamma} \{n - k_s\} \ < \
\frac{1 - \gamma}{1 - \alpha} \left(\frac{n - k_s}{1 - \gamma} \
+ \ 1\right)\,,
$$
whence inequality (3.19) follows with evidence. \hfill $\Box$
\vspace*{3mm}

The estimation of theorems 2 and 3 at availability of saturated
grids has some redundancy, which essentially less, if initial
zeroes density is sparse. The valuation will be also redundant
and in case of incorrect degree-system, that is in variant given
non-canonical decomposition of zeroes frequency (density). As
follows from appendices [1], the best approximation of
valuations (3.18) in the majority of cases is maximum series of
correct binary system \ $SS_2$\,.

Quite similar state is observed for valuations (3.19) of
maximum series and for systems without multiple zeroes \ $SS$\,,
\ though the formulas for them are not given here. However
role of zeroes density for such systems is same, hence, and the
high valuations should coincide, and obstacle can not act
concrete formula realization for values of maximum
series \ $MSR_{\alpha} \{n\} [SS]$\,.
\vspace*{3mm}

The reception of algorithms for calculation of exact values of
the maximum series (theorem 1) and upper generalized estimation
(theorem 3) permit to make conclusions for a specific class of
systems without multiple zeroes, which can serve by necessary
spring-board at reception essentially important generalization.
\vspace*{5mm}

\begin{center} {\Large \bf 4. Imaging principle and main theorem}
\end{center}
\vspace*{3mm}

The expressions and high valuations of maximum series of the
theorems 1--3 are found for systems without multiple zeroes.
However systems with the multiple zeroes present greatest
interest just. All systems of the second and third type without
fail have zeroes of multiplicity higher unit.  In particular, if
in filling \ $Z_n$ \ there will be though one pair of grids with
modules \ $(a_i,\,a_j) = 1$\,, \ in the period multiple zeroes
will meet. At the same time fulfilment of relation
$$
S(a_i),\,S(a_j) \ \subset \ Z_n: \qquad (\,a_i,\,a_j\,) \ \neq \ 1\,,
\qquad \forall \ 1 \le i \neq j \le \ n \eqno{(4.1)}
$$
not yet guarantees absence of multiple zeroes. For the most important
in appendices systems of the third type the share of zeroes
multiplicity higher unit increases to unit at \ $n \to \infty$\,.

Results of the previous chapter prove, that the problem of
search of maximum series of fillings \ $Z_n$\,, \ their values \
$msr_n (q)$ \ and particularly valuations is decided
successfully for systems without multiple zeroes. From
here there is idea of approximation of fillings for any systems
by fillings without multiple zeroes, let even for level of the
majorizing characteristics. This purpose some mental construction,
named imaging principle corresponds.

We are addressed to fillings with the extreme characteristics
for a given set of \ $n$ \ grids. Naturally, only non-singular systems
are implied.
\vspace*{3mm}

The unregulated fillings are not random, but found during
exhaustive search or different way for an evaluation of
interesting numerical characteristics. For example, maximum
series generally can be guaranteed are found, identified and
are appreciated only on the class of unregulated fillings. The
modifications of the sieving process do not removal from the
numerical nature of the worked up sequence, as against it is in
the class \ ${\cal L}$\,. \ The call to Boolean elements allows
to use completely other constructions.

The main difference of one method from other consists in an
evaluation of the majorizing characteristics of the \ $n$-filling.
It can and even conveniently be passed from the integer analysis
of outcomes of each grid effect to learning cooperative
influence of zeroes frequencies of \ $n$-filling. It happens for
all period, and not just on an initial interval of length about
\ $p_{n+1}^2$\,, \ which one restricts itself the sieving process.
Therefore fillings method can be interpreted as method of the
analysis of frequencies of Boolean zeroes -- results of grids
products.
\vspace*{2mm}

{\bf Definition 20.} If \ $E_n$ \ there is quantity of units on period
of length \ $PZ_n $ \ for nonsingular \ $n$-filling \
$Z_n(a_1,...,a_n)$ \ of the system \ $SS$\,, \ and \ $H_n =
PZ_n - E_n$ \ is number of zeroes, the value of {\sl total zeroes
(sum of multiplicity)} \ $H_n^*$ \ is
$$
H_n^* (Z_n) \ = \ PZ_n\, \ \sum_{i=1}^n \ \frac{1}{a_i}\,,
\qquad PZ_n \ = \  ((a_1,a_2,...,a_n))\,.
$$
\vspace*{1mm}

{\bf Definition 21.} The main object of the research is two-dimensional
strip region of binary values of volume \ ($\,n\,; \ \infty\,$) \
of grids

$$
\begin{array}{ccc}
S(a_1): \qquad & ...l_{j-1}^{(1)} \ l_j^{(1)} \ l_{j+1}^{(1)}... &  \\
S(a_2): \qquad & ...l_{j-1}^{(2)} \ l_j^{(2)} \ l_{j+1}^{(2)}... & \forall
i \ \{ l_j^{(i)} = 1\}: \ \ t_j = 1; \\
.......... \qquad & ......................... &  \\
S(a_n): \qquad & ...l_{j-1}^{(n)} \ l_j^{(n)} \ l_{j+1}^{(n)}... & \exists
i \ \{ l_j^{(i)} = 0\}:\\
\Longrightarrow \qquad & - \to - \to - \to - \to - & 1 \le k
\le m \le n \\
W_n: \qquad & ...t_{j-1} 0_j^{(1)} 0_j^{(2)}...\,0_j^{(k)} t_{j+1}...\ , &
\end{array} \eqno{(4.2)}
$$
where \ $m$ is quantity of zeroes in column \
$l_j^{(1)},\,l_j^{(2)},\,...\,,\,l_j^{(n)}$\,. \ The imaging \
$W_n$ \ of grids \ $S(a_i)$ \ on the sequence of the same values (4.2)
$$
S(a_1)\,\&\,S(a_2)\,\&\,...\,\&\,S(a_n) \
\stackrel{f(k)}{\Longrightarrow} \ W_n\,; \qquad k \ \le \ m
$$
is created on rules: the unity element \ $t_j = 1$ \
corresponds to unity column; the imaging \
$\stackrel{f(k)}{\Longrightarrow} \ W_n$ \ transforms zeroes to
the sequential series from \ $k \ (1 \le k \le m)$ \ zeroes, if
column has \ $m$ \ zeroes \ $ (1 \le m \le n)$\,.
\vspace*{3mm}

Variant of four grids of prime modules we shall give
as an example of the imaging:
$$
\begin{array}{cl}
S(3): \quad & \qquad ...101101101101101101101101... \\
S(4): \quad & \qquad ...101110111011101110111011... \\
S(5): \quad & \qquad ...101111011110111101111011... \\
S(7): \quad & \qquad ...101111110111111011111101... \\
\Longrightarrow \quad & \qquad - \to - \to - \to - \to - \to - \to - \to - \\
W_4: \quad & \qquad ...100011000000001001000010100001...\ .
\end{array}
$$

Here value \ $k = 3$ \ is the imaging of column of four zeroes
in first case, where zero of product is substituted by three
zeroes. Further product of two zeroes represented by pair of
zeroes in imaging \ $W_4 $ \ in all four cases.

The imaging of multiple zeroes should be realized at
the expense of increase of period length, as the quantity of
units \ $E_n$ \ on period is constant. It in sufficient
measure the conditional increase of period is followed to perceive only
as the tool of obtaining of necessary numerical
characteristics. Length of period of imaging \ $W_4$ \ in an
example is conditional also. But it exceeds the value \
$PZ_4$\,.

Sense of such conditional increase of period at the expense of
multiple zeroes consists in an evaluation of the local interval of
product of grids. Then the product of grids with multiple
zeroes can by suitable shifts be resulted with diminished quantity
of such multiple zeroes or even with their complete liquidation.
Series and maximum series of zeroes place just on intervals of the
similar type, that is highly small length, it is much less value
of  period \ $PZ_n$\,. \ But the redistribution is possible not
always.
\vspace*{3mm}

{\bf Definition 22.} If the imaging \ $W_n $ \ takes into account all \
$m$ \ of zeroes, that is in all cases \ $k = m$\,, \ that
this complete imaging of zeroes \ $(W_n^*)$\,. \ If the
value \ $k > 1$ \ even in one case, but is not always executed
\ $k = m$\,, \ we have incomplete imaging \ $(W_n^{**})$\,. \
Direct imaging of zeroes \ $(W_n)$ \ is obtained at \ $k = 1$\,,
\ given in all period. The imagings are \ $W_n^{*} \ = \ W_n^{**} \ =
\ W_n$ \ for systems \ $SS'$\,.
\vspace*{3mm}

Definitions 20 -- 22 represent itself as the scheme of the
imaging principle.
\vspace*{3mm}

{\bf Definition 23.} The frequency of zeroes \ $\gamma_n$ \ on period
of imaging -- from direct \ $(W_n)$ \ up to complete \
$(W_n^*)$\,, \ serves for the basis adequate ratings of the maximum
series \ $MSR_n(q)$ \ for all systems with multiple zeroes:
$$
1 - \frac{E_n}{PZ_n} = \frac{H_n}{PZ_n} = (\gamma_n) \ < \
(\gamma_n^{**}) \ < \ (\gamma_n^*) = 1 - \frac{E_n}{PZ_n^*} =
\frac{H_n^*}{PZ_n^*}\,.
$$
\vspace*{1mm}

Clearly, the evaluations of the maximum series \ $MSR_n(q)$ \
should be constructed because of frequencies of zeroes \
$(\gamma_n)$ \ of direct imaging for values \ $q$\,, \ near from
quantity of units \ $E_n$ \ on period. Quite other position
develops with evaluations of series \ $MSR_n (q)$ \ for small
\ $q$ \ or even for \ $q = 0$\,. \ The intermediate frequency \
$(\gamma_n^{**})$ \ certainly is frequency of incomplete imaging
\ $(W_n^{**})$\,. \ It is necessary to mark, if the parameters
\ $\gamma_n(q)\,, \ n$ \ and \ $q$ \ are given, it is
possible to forget about concrete set of grids \ $\{a_i\}$\,.

Really, these parameters are sufficient for obtaining unknown
quantities, but approximate ratings. However precise definition
of maximum series \ $MSR_n (q)$ \ requires of greater.
\vspace*{3mm}

Let's formulate main definition touching the means of learning of
introduced systems and fillings.
\vspace*{5mm}

{\sl Definition} {\bf 24. The fillings method is research
of properties and characteristics of \ $n$-fillings and systems
\ $SS$ \ with multiple zeroes, and also obtaining of the series
of fundamental numerical estimations with the help of imaging of
zeroes \ $W_n$ \ of all types.}
\vspace*{5mm}

Generally filling is reduced to imaging of two-dimensional strip
region of elements on one-dimensional with partial or complete
conversion of multiple zeroes in single. Such extended filling
on changed period is base of learning of properties of source
distribution of zeroes and units as products of grids for
different classes of systems.

The explained principles of the fillings method have allowed to
reveal central relation, all rests are  consequences from which. The
formula reflects statement, is foolproof enough expressed
mathematically and claiming to be main for the rather vast class
of the tasks of number theory. The versatility of this relation
do not know exceptions on set of nonsingular systems and
fillings.

The offered thesis does not imply dependence from principle of imaging.
The principle only explains paths and sources of the approach,
sense and parents of appearance. At the same time at all riches
of applications and importance of the obtained outputs, the
thesis can be surveyed in the different forms with direct,
incomplete or complete imaging. First of all call to this or
that form of the main theorem is determined by the degree of
correspondence to content and fundamental essence of the
fillings method.
\vspace*{6mm}

{\sl The main theorem.} {\bf Estimation of maximum series \
$MSR_n(q)$ \ is valid for anyone \ $n$-filling $Z_n$
of the arbitrary nontrivial system $SS$:}
$$
msr_n(q) \ < \ \frac{n + q}{1 - \gamma_n^*} \ + \ 1 \ = \ {\cal
M}_n(q)\,, \qquad \gamma_n^* \ = \ \frac{H_n^*}{E_n + H_n^*}\,,
\eqno{(4.3)}
$$
{\bf where \ $E_n$ \ is quantity of units of filling's period \
$PZ_n$ \ and \ $H_n^*$ \ is total (sum of multiplicity)
zeroes. Value of density \ $\varrho_n$ \ for each \ $1 \ \leq
\ n \ < \ \infty$ \ will be discover always}:
$$
msr_n(q) \ \leq \ \frac{n + q}{1 - \varrho_n} + 1\,, \quad
\sup_{SS} \ \sup_{\{Z_n\} } \ \max_{0 \leq q < \infty} \ \{ \varrho_n
(Z_n \subset SS) \} < \gamma_n^*\,, \eqno{(4.4)}
$$
$$
\inf_{SS} \ \sup_{\{Z_n\} } \ \max_{0 \leq q < \infty} \ \{ \
\varrho_n (Z_n \subset SS) \ \} \ = \ \gamma_n\,, \qquad
\gamma_n \ = \ \frac{H_n}{PZ_n}\,,  \eqno{(4.5)}
$$
{\bf where \ $\sup$ \ and \ $\inf$ \ are in the class of
nontrivial systems. But concrete kind of \ $n$-filling can define
incomplete imaging and appropriate frequency of zeroes \
$\gamma_n^{**}$\,, \ and consequently unimprovable estimation
for some system \ $SS$\,}
$$
msr_n(q) \ \leq \ \frac{n + q}{1 - \varrho_n} \ + \ 1\,, \quad
\sup_{\{Z_n\} } \ \max_{0 \leq q < \infty} \ \{ \varrho_n
(Z_n \subset SS) \} < \gamma_n^{**}\,. \eqno{(4.6)}
$$

{\bf The following inequalities take place for all classes of
nonsingular \ $n$-fillings \ $SS$\,, \ where} \ $0.5 < C =
C(SS) \le 1$\,:
$$
\gamma_n \in Z_n \subset SS\,: \qquad C\,\frac{n +
q}{1 - \gamma_n} \ + \ 1 \ < \ msr_n(q) \ < \ 2\,\frac{n + q}{1
- \gamma_n} \ + \ 1\,. \eqno{(4.7)}
$$

{\bf The dependence of the zeroes density \ $\varrho_n$ \
from \ $q $ \ leads to the form}
$$
\exists\,q;\,\varrho_n (q): \quad msr_n(q) \ = \ \frac{1 + q}{1
- \gamma_n} \ + \ 1\,; \qquad \lim_{q \to \infty} \ \varrho_n
(q) \ = \ \gamma_n\,. \eqno{(4.8)}
$$
\vspace*{6mm}

The first part (4.3) of theorem states about existence absolute
majorizing frequency (density) of zeroes for arbitrary filling
of any nonsingular system \ $SS $ \ of grids. The density of
zeroes of complete imaging \ $\gamma_n^*\,, \ (\gamma_n^* \ge
\gamma_n^{**} \ge \gamma_n )$ \ represents itself as such
frequency. It determines unconditional and even an inaccessible
upper-bound estimate \ ${\cal M}_n(q)$ \ of an appropriate
maximum series \ $MSR_n (q) [SS]$\,.

The logic and constructibility of such evaluation form \ of \
\ maximum \ series \ \ $MSR_n(q)$ \ is justified by complete
coincidence with an evaluation for fillings without multiple
zeroes. The transition to the relation (4.3) for arbitrary systems,
including with multiple zeroes, becomes well-grounded after
operation of complete imaging of multiple \ zeroes \ of \
$n$-dimensional strip region of binary elements (zeroes and
units of grids). In an outcome all zeroes on complete
(extended) period \ $PZ_n^*$ \ have multiplicity of unit.

Each system \ $SS$\,, \ any more not speaking about \
$n$-filling, has the majorizing density of zeroes, which one
here is marked \ $\varrho_n$\,. \ It provides an evaluation of
the inequality (4.4), replicated main form (4.3) at all values \ $n
\ge 1\,; \ q \ge 0$\,. \ Nevertheless, top and bottom boundary
of densities of zeroes on the class of all acceptable systems
of grids are, accordingly, the frequencies \ $\gamma_n^*$ \ and
\ $\gamma_n$ \ of relations (4.4,\,4.5).

However redundancy of an evaluation \ $MSR_n(q)$ \ for systems
without multiple zeroes, especially has an effect in variant of
arbitrary fillings. It is explained to the redundancy of an
evaluation (4.3) for small intervals (that is \ $n $ \ and \
$q$\,), in the total reduces in such interval, for which one
there is no filling without multiple zeroes. Thereby some
multiple zeroes appear superfluous in data conditions and
consequently is acceptable to be restricted to incomplete
imaging of zeroes. So frequency of zeroes \ $\gamma_n^{**}$ \
lesser what \ $\gamma_n^*$ \ but exceeding \ $\gamma_n$ \
occurs.

All these reasons reduce to appropriate densities of zeroes and
evaluations (4.6). The transition to more precise modification of
the method of fillings gives detection of multiple zeroes which do not
lead to increase of majorizing density. The rather
outstanding part of such inefficient multiple zeroes of all
period can appear for number of systems. Naturally, it
essentially will decrease value \ $\gamma_n^{**}$ \ about \
$\gamma_n^*$\,. \ At the same time it is impossible to guess,
that the value \ $\gamma_n^{**}$ \ will reach value \ $\gamma_n$
\ for great many of multiple zeroes.

The unimproving evaluations (4.7) of the maximum series \ $MSR_n
(q)$ \ was obtained from the study of axis configurations of the
prime system \ $SP_1$\,.

The realizability of equality (4.8) for some values \ $q$ \ (for
example, for \ $q = E_n - 1$\,) \ at \ $\varrho_n =
\varrho_n'(q) = \gamma_n$ \ is the quite definite
characteristic of majorizing density of zeroes \ $ \varrho_n $
\ in the expression (4.6). In this case maximum series coincides an
evaluation for \ $n = 1$ \ and it is equal to length of period \
$msr_n (E_n - 1) = PZ_n$\,. \ From here it is clear, only value
\ $\gamma_n$ \ can be by limit (4.8) for constant \ $n \ (n > 1)$ \
and \ $q \to \infty$ \ for density of zeroes \ $\varrho_n$ \ in an
estimation (4.6).

Thus, majorizing estimation of the maximum series \ $MSR_n(q)$ \
is connected inversely proportional dependence with the density
of zeroes in the period. It appears by the adequate characteristic
of imaging (complete or incomplete) multiple zeroes of the strip
region of elements. Zeroes frequency of complete imaging \
$\gamma_n^* \ = \ \frac{H_n^*}{E_n + H_n^*}$ \ in elongated
period \ $PZ_n^* = E_n + H_n^*$ \ (sometimes \ $PZ_n^* \gg PZ_n$) \
is thus natural absolute majorant for any systems and fillings.

So on the first view indisputable on logic and validity
the thesis nevertheless, requires the proof of impossibility
of sieve substitution of multiple zeroes of units on an
interval of the maximum series \ $MSR_n $\,. \ It is really
impossible as well as in variant of fillings without multiple
zeroes. As the distribution of zeroes of grid is uniform, the
transition of units in zeroes is obliged to lead in restoring
units on adjacent places. As the frequency \ $\gamma_n^*$ \
registers and takes into consideration zeroes of all multiplicity
without exception.

At the same time one of proofs of the main theorem is received
from the analysis of fixed distributions so named axial series
for the period of \ $n$-fillings. It does not lean on the imaging
principle, but confirms legitimacy of its introduction and
consideration. Thus the majorizing constant two is found for the third
form of the main theorem. This constant is unimproved, as it is given
below.
\newpage
\vspace*{1mm}

\begin{center} {\Large \bf 5. Premises of an evidence of the main
theorem} \end{center}
\vspace*{3mm}

The offered approach is not uniquely possible.
\vspace*{3mm}

{\bf Definition 25.} \ \ Algorithm \ of \ construction of maximum series
of zeroes \ $MSR_n (q)$ \ shall be name one-sided \ $Z_n^{(1)}$\,, \ if the
filling is conducted by half-grids from chosen beginnings in one party
(for example, right).
\vspace*{3mm}

Regulatedness, semi-regulatedness and unregulatedness of a mechanism
of fillings here remain in complete force at formation of any series of
zeroes \ $SR_n (q)\,, \ q \ge 0$\,, \ including and maximum \ $MSR_n (q)$\,.
\ One-sidedness does not depend from other characteristics of
filling algorithm. The half-grids is one-sided infinite grids.

It is possible to remind, at fixing of beginning strict
regulatedness means product the first (from the right) unit of current
filling \ $Z_{n-1}$ \ and boundary unit of a series \ $SR_{n-1}$ \
with the zero of a grid \ $S(a_n)$\,. \ The received product then will
be filling, but maximality of a formed series is observed not always
and requires separate consideration. We shall suggest another scheme of
fillings mechanism.
\vspace*{3mm}

{\bf Definition 26.} Algorithm of construction of series and maximum
\ $MSR_n (q)$ \ series of zeroes shall be name two-sided \ $Z_n^{(2)}$\,,
\ if the filling is conducted by complete grids (two-sided infinite)
till both party from a point of beginnings.
\vspace*{3mm}

The sense of the introduction of one-sidedness of fillings mechanism
clears up for regulated and semi-regulated algorithms of search and it
is in direct dependence on expressions of next modules of allocated
filling.
\vspace*{3mm}

{\bf Consequence 1.} The values of maximum series of zeroes \ $msr_n(q)$ \
of unregulated fillings coincide at one-sided and two-sided filling.
$$
MSR_n (q,\,Z_n^{(1;2)},\,ZN_n,\,SS)\,: \qquad msr_n (q,\,Z_n^{(1)},\,SS)
\ = \ msr_n (q,\,Z_n^{(2)},\,SS)\,. \eqno{(5.1)}
$$
\vspace*{1mm}

{\sl Proof.} At unregulated filling search of maximum series of zeroes
equivalent complete selection of all possible products of grids. Then
variants of a configuration of zeroes of such series
\ $MSR_n (q,\,Z_n^{ (k)})$ \ can not coincide. Obviously their
arrangement on a numerical axis (or half-axis) differs, but the lengths of
such series coincide on a sense of construction. Thus the expression (5.1)
is executed for any parameters \ $(n,\,q)$\,, \ fillings \ $Z_n \subset SS$
\ and systems \ $SS$\,. \hfill$\Box$
\vspace*{4mm}

{\bf Theorem 4.} The construction of maximum series of zeroes \ $MSR_n (q)$
\ is oriented to one-sided filling \ $Z_n^{ (1)}$ \ for systems of grids
$SS_{[ 2]} (a_i)$, the unequal next modules of which are connected by the
expression \ $a_{i+1} \ \ge \ 2 \cdot a_i$ \ (for \ $a_{i+1} \neq a_i$)\,.
$$
MSR_n (q,\,SS_{[2]})\,: \qquad  MSR_n \left(q,\,Z_n^{(1)},\,SS_{[2]}\right)
\ \Rightarrow \ msr_n (q)\,, \quad a_i \in SS_{[2]}\,. \eqno{(5.2)}
$$
\vspace*{1mm}

{\sl Proof.} It is necessary at once to note, that for grids of equal
modules the factor of one-sidedness or two-sideness is away by virtue of them
indistinction. Therefore is allowable to be limited to consideration of a
system with growth of modules of all grids. Thus there is inessential mutual
prime of modules of this system of the first type -- that is with a sum of
values of reverse modules less unit, if to exclude a trivial binary
system.

Let down to some stage, that is up to parameter \ $n - 1$ \ inclusive, the
search of maximum series \ $MSR_k (q) $ \ at \ $k \le n-1$ \ was maintained
pursuant to the scheme of one-sided fillings. It means, that all grids
down to \ $S(a_{n-1})$ \ participated in one-sided filling, and
concentration of zero from the right from a index point appreciably
higher, than at the left, if instead of half-grids of the statement
temporarily to consider complete grids.

In such case the second zero of the half-grid \ $S(a_{n-1})$ \ on a
constructed interval of a series will meet earlier, than the first zero
of a next half-grid \ $S (a_{n})$\,. \ It means, that the one-sided
algorithm of filling contains higher potential of growth of a series,
including maximum series \ $MSR_n (q) $\,. \ Unique difficulty is made
in possible earlier occurrence of a zero of multiplicity two, that
is in crossing of zero, that can partly deform a reasoning.
Just is here allowable some variability in filling, not changing the to
common scheme of one-sidedness.

Really, the primary occurrence of zeroes of grids of smaller modules
forces to address to the one-sided scheme, as the attempts of two-sided
filling are obliged to result in downturn of concentration of zeroes in a
interval because of a inequality \ $a_n / 2 \ \ge \ a_{n-1}$\,. \ Zeroes
of a grid \ $S(a_n)$ \ are just on distance \ $a_n / 2$ \ from a central
point of a prospective two-sided interval.

But then, if to take into account specific character of such fillings
for the first type system, value \ $msr_{n-1} < a_{n-1}$\,, \ and the
second zero of a grid \ $S(a_{n-1})$ \ by no way can not enter in
a series \ $SR_{n-1}$ \ or \ $MSR_{n-1}$\,. \ For this reason the zero
of a grid \ $S(a_n)$ \ as grid of two-sided filling can not already
participate in formation of a series of the heaviest length. And on a next
step \ ($n + 1$) \ significance of the second zero of a grid
\ $S (a_{n-1})$ \ becomes determining. It liquidates the unjustified
claims of grids \ $S(a_{n})$ \ and \ $S(a_{n + 1})$ \ for
two-sided algorithm of filling.

This conclusion is reflected by expression (5.2). The orientation to
one-sided filling is wholly determined by growth of a next module. In
such case attempts of two-sided filling \ $Z_n^{(2)}$ \ are inexpedient,
as the one-sided filling \ $Z_n^{(1)}$ \ essentially more effective
conducts to formation of maximum series \ $MSR_n (q)$ \ for systems of
grids of a kind \ $SS_{[2]} (a_i)\,, \ a_{i + 1} \ge 2a_i$\,, \
if grids are unequal \ $(a_i \neq a_{i+1})$\,. \hfill$\Box$
\vspace*{3mm}

The narrowing of a class of systems \ $SS \subset SS_{[2]} $ \ admits
concreteness of the statement. This statement connects important
properties of fillings.
\vspace*{1mm}

{\bf Theorem 5.} Construction of maximum series of zeroes \ $MSR_n (q)$
\ at any \ $q \ge 0$ \ is provided with regulated one-sided algorithm of
filling \ $ZU_n^{(1)}$ \ for all non-singular degree-systems \
$SS_d(a_i)$ \ of any basis \ $d \ge 2$\,.

$$
MSR_n (q,\,SS_d): \quad MSR_n \left\{q,\,ZU_n^{(1)},\,SS_d(a_i)\right\}
\Rightarrow msr_n (q), \quad ZU_n^{(1)} \subset Z_n^{(1)}.
\eqno{(5.3)}
$$
\vspace*{1mm}

{\sl Proof.} Obviously, the degree-system \ $SS_d$ \ at any basis \
$d \ge 2$ \ enters in a class of systems of primary growth of modules \
$SS_d \subset SS_{[2]}$\,. \ Such system has not multiple zeroes at any
\ $n$ \ from definition.

But absence of multiple zeroes at such correlation of modules of grids
transform one-sidedness and regulatedness of filling into the compulsory
rules. Any infringement of such rules conducts to distortion of maximum.
The expression (5.3) establishes inevitability of algorithm of regulated
one-sided filling \ $ZU_n^{(1)}$ \ at formation of all maximum series
\ $MSR_n (q)$\,. \hfill$\Box$
\vspace*{3mm}

{\bf Theorem 6.} The one-sided filling remains effective for systems
without multiple zeroes \ $SS'$\,. \ Strict regulatedness can be sometimes
infringed at realizations of grids of modules
\ $a_{i+1} < 2 \cdot a_i$\,, \ if \ $a_{i+1} \neq a_{i}$\,.
$$
MSR_n (q,\,SS'): \quad MSR_n \left\{q,\,SS'(a_i)\right\}
\ \Rightarrow \ MSR_n \left\{q,\,Z_n^{(1)} = ZU_n^{(1)}
[a_{i+1} \ge 2 \cdot a_i]\right\}. \eqno{(5.4)}
$$
\vspace*{1mm}

{\sl Proof.} A system of grids \ $SS'$\,, \ all fillings of which do not
contain zeroes of multiplicity above unit, is in the correlation
\ $SS_d \subset SS'$\,. \ But inclusion of a system \ $SS'$ \ in a class of
systems \ $SS_{[2]}$ \ optionally, that is probably \ $SS' \not\subset
SS_{[2]}$\,. \ For example, grids \ $S(4),\,S(6)$ \ will form such
2-filling, which can enter in a system \ $SS'$\,. \ At the same
time fitting to a system of the first type \ $\sum\,1 / a_i \ \le \ 1$ \
precisely indicates on basic character of degree-filling at creation of
systems \ $SS'$ \ without multiple zeroes.

Really, mutual non-primeness of modules of any grids $(a_i,\,a_j) \neq 1$
of the systems demonstrates or very the large modules $a_i$ of all
grids as products prime multipliers for initial fillings, or reasonably
foreseeable growing simularity with a degree-system. But the first
variant can ensure smaller density of zeroes for such fillings, at the same
time compelling pass to the degree-schemes at growth of parameter \ $n$\,. \
Second, more natural and in a limit the same variant is stipulated by
conditions of the theorem.

In this case any filling of a system without multiple zeroes \ $SS'$ \ can be
submitted on the basis of consecutive representations of some filling of
a degree-system. Actually, if the determined filling is realized by grids \
$\{S(a_i)\}$\,, \ at first degree, and consequently not possessing
multiple zeroes, formation of other filling, in which modules \ $a_i$ \
change by product with any, including identical multipliers \

\hspace*{-6mm}$ b_i\,: \ a_i \Rightarrow b_i \cdot a_i\,, \ b_i \in
{\bf N}$\,, \ leaves filling in class without multiple zeroes.

Such way can essentially increase a period of filling and in parallel to
decrease density of zeroes. Therefore for growth of initial and received
density the introduction of grids is allowable, but already not any
degree, but agreed with current modules, not to admit multiple zeroes.

Or else, in the basis of any filling, included in system without multiple
zeroes, lies some degree filling (can be, even not unique). The number of
modules of its grids (may be, even all) are transformed by product with
multipliers \ $b_i$\,, \ in common case not connected with the basis \
$d$ \ of degree filling. Clearly, that such filling nor will give multiple
zeroes. It though because are used \ $1 / b_i$ \ zeroes of each grid \
$S(a_i)$\,. \ But the first zero of these two grids of formed series can
coincide, and distinction multipliers \ $b_i$ \ in a condition to infringe
strict regulatedness of degree filling. It should note, that \
$b_i$ \ grids \ $\{S(b_i \cdot a_i)\}$ \ restore density of zeroes,
given by a grid \ $\{S(a_i)\}$\,, \ as well as itself grid \ $\{S(a_i)\}$\,.

The described way of construction of fillings without multiple zeroes fixes
result, lying closer to results for degree-systems, than less various
multipliers \ $b_i$\,, \ different from unit, in a transformed
set of modules \ $\{b_i \cdot a_i\}$\,. \ Thus there can arise an
expression for modules
$$
b_i \cdot a_i \ < \ 2 \cdot b_j \cdot a_j\,, \quad b_i \cdot a_i \ > \
b_j \cdot a_j\,; \qquad \{S(4),\,S(6): \ \ 6 < 8\,; \ 6 > 4\}\,,
$$
which partial infringement of regulatedness admits at formation of
maximum series. In indicated examples with two grids, in particular,
one infringement of regulated algorithm takes place in case \
$S (4),\,S(6)$\,: \ $8 = msr_2 [3,\, ZN^{(1)}] \ > \ sr_2 [3,\,ZU^{(1)}]
= 7$\,.

If to take into account limited opportunities of growth of a set of various
multipliers \ $b_i$ \ because of a essential increase of transformed
modules of grids (that decreases opportunities of occurrence of
close modules), from here and statement of the theorem and expression
(5.4) follows. In conditions of a system of the first type infringements
of regulated algorithm are reasonably rare on a common background.
At primary growth of next modules obviously use of regulated fillings
for all variants of series with units \ $q \ge 0\,: \ \
ZU_n^{(1)}[a_{i+1} \ge 2 \cdot a_i]$\,.
\hfill$\Box$
\vspace*{3mm}

{\bf Theorem 7.} Just the algorithm of one-sided filling \ $Z_n^{(1)},
\ n > n_0$ \ is obliged to demonstrate efficiency for systems of grids \
$SS_{(I)} = \{S(a_i)\}$ \ of the first type, in any case, since some
parameter \ $n_0$ \ for \ $n_0 > 1$\,.
$$
MSR_n [q,SS_{(I)}]: \ \ MSR_n \left\{q,SS_{(I)}\left(\sum
\frac{1}{a_i} \le 1\right)\right\} \Rightarrow
MSR_n \left\{q,Z_n^{(1)}\right\}, \ n > n_0. \eqno{(5.5)}
$$
\vspace*{1mm}

{\sl Proof.} Degree-systems \ $SS_d$ \ and the systems without multiple
zeroes \ $SS'$\,, \ naturally, satisfy to conditions of the statement,
that confirm the proved theorems 5 and 6. About it speaks and theorem 4,
considering grids with growing twice and more modules. Systems of
the first type \ $SS_{(I)}$ \ include all listed, but  wider of them.
In particular, this system \ $SS_{(I)}$ \ grids can quite include grids with
modules, lying in the borders \ $a_i \ < \ a_{i+1} \ < \ 2 \cdot a_i$\,.

At the same time relative quantity of such abnormal inequalities for modules
of next grids can not be appreciable. Otherwise the sum of reverse modules
appears more unit. But in such case there are no essential handicapes
for a establishment of one-sided filling \ $Z_n^{(1)}$\,, \ which can
appear effective in general at all \ $n$\,, \ and if similar is not observed,
will eventually set in for all \ $n > n_0$ \ at some \ $n_0$\,. \
Integrated factor of one-sidedness of algorithm of filling not to overcome
by separate infringements of primary increase of modules.

Steady two-sidedness of algorithm of filling requires firm advantage of a
set of slowly growing modules of grids. Otherwise the advancing growth
of set of units on a interval of a prospective maximum series will
transform two-sided algorithm in inefficient. It reflects expression (5.5).
The one-sided character of filling at formation of the heaviest series
\ $MSR_n \left\{q,\,Z_n^{(1)} \right\} $ \ is established for all \
$n > n_0$\,. \hfill$\Box$
\vspace*{3mm}

The question should arise about a principle of formation of maximum
series for systems with in comparison slow growth of modules of next grids.
\vspace*{3mm}

{\bf Theorem 8.} The two-sided filling \ $Z_n^{(2)}$ \ is obliged to
appear effective for systems of grids \ $\{S(a_i)\}$\,, \ the
overwhelming part of modules of which are connected by a correlation \
$a_{i+1} \ < \ 2 \cdot a_i$ \ (not including equal), since some \ $n_0$\,.
$$
MSR_n (q,\,SS): \quad MSR_n \left\{q,\,SS(a_{i+1} < 2 \cdot a_i)\right\}
\ \Rightarrow \ MSR_n \left\{q,\,Z_n^{(2)}\right\}, \ n > n_0. \eqno{(5.6)}
$$
\vspace*{1mm}

{\sl Proof.} The primary affinity of grids \ $a_{i+1} \ < \ 2 \cdot a_i$
\ permits to conclude, that such systems concern to second and, mainly,
to the third type. The possible grids of a identical module even
are not taken into account. Zeroes of multiplicity above unit not only
become ineradicable, but also can form a overwhelming set, as far as the
talk goes only about non-singular fillings and systems. As a example it is
possible to result a system primes \ $SP_0$\,. \ $n$-filling of this system
in a condition to have on a period a sum of multiplicities of zeroes,
in a many times superior period.

In such case the algorithm of two-sided filling becomes expedient.
The high density of zeroes in a period \ $n$-fillings and, accordingly, in
comparison small quantity of units in a vicinity of a formed maximum series \
$MSR_n (q)\,, \ q \ge 0$\,, \ admit a opportunity some other algorithm
of filling \ $Z_n^{(2)}$\,. \ The algorithm extreme effectively uses
redundant frequency of zeroes and redistribution of zeroes of high
multiplicity.

At the same time particular search in variants of initial values of
parameter \ $n$ \ as quantities of grids quite we allow one-sided algorithm
of filling. However the conditions of the statement prove, that in result
moment of saturation of an interval of a series by zeroes will come.
Then it is possible to recollect a reserve, about the forgot the second
side, not taken into account at one-sided algorithm \ $Z_n^{(1)}$\,.

The expression (5.6) emphasizes inevitability of the address
to two-sided algorithm of filling \ $Z_n^{(2)}$\,. \ It will be realized
for all \ $n > n_0$ \ with some \ $n_0$ \ by search of maximum series \
$MSR_n (q)$ \ in systems of predominary slow growth of modules of grids
and, hence, high density of zeroes. \hfill$\Box$
\vspace*{1mm}

Sometimes the researched system of grids \ $SS$ \ extremely slightly differs
from considered systems. The characteristics them are also reasonably close.
\vspace*{3mm}

{\bf Theorem 9.} Inclusion or exception of any finite set of grids, not
infringing non-singular system, results in that will be always found
such \ $n_0$\,, since which is restored one-sideness or two-sideness of
fillings for any system of grids \ $SS = \{S(a_i)\}$\,.
$$
MSR_n (q): \quad MSR_n \left\{q,\,SS(a_i),\,Z_n^{(k)}\right\}
\ \Rightarrow \ MSR_n \left\{q,\,SS(a_i \cup b_j),\,Z_n^{(k)}\right\},
\ n > n_0. \eqno{(5.7)}
$$
\vspace*{1mm}

{\sl Proof.} The condition of non-singular system and fillings obviously
in a result of inclusion of some fixed set of grids.
One-sidedness or two-sidedness of fillings is determined by wittingly
advantage of inequalities \ $a_{i+1} \ge 2 \cdot a_i$ \ or \ $a_{i+1}
< 2 \cdot a_i$ \ in a chain of comparisons of modules. It remains and
after inclusion of any finite set of grids in a system. That is if for
one system the kind \ $Z_n^{(k)}\,, \ k = 1;2$ \ of fillings algorithm
is established at \ $n > n_0$\,, \ for other considered system the
same kind of algorithm is observed at other parameters \ $m > m_0$\,.

Really, another set of grids \ $\mathop{\cup}\limits^N\,S(b_j)$\,, \ let
even it and finite, up to some moment, determined by value of parameter
\ $n_0$\,, \ in a condition seriously to deform the scheme of algorithm
of filling \ $Z_n^{(k)}$\,. \ But the scheme \ $Z_n^{(k)}$ \ is restored
(5.7) without fail, as far as the system on a condition initially
has unlimited total stabilizing effect. \hfill$\Box$
\vspace*{3mm}

From here immediately the next conclusion follows for a class of systems
of the second type. This class \ $SS_{(II)}$ \ is real and spacious in
the fillings method.
\vspace*{3mm}

{\bf Theorem 10.} The one-sideness of fills \ $Z_n^{(1)}$ \ is established
always for systems of grids \ $SS_{(II)} = \{S(a_i)\}$ \ of the second
type, since some \ $n_0$\,.
$$
MSR_n (q,\,SS): \quad MSR_n \left\{q,\,SS_{(II)}\right\}
\, \Rightarrow \, MSR_n \left\{q,\,Z_n^{(1)},\,SS_{(II)}\right\}, \
n > n_0. \eqno{(5.8)}
$$
\vspace*{1mm}

{\sl Proof.} From definition the system of grids of the second type \
$SS_{(II)}$ \ has a finite sum of reverse modules. Such sum is more unit.
It obviously means, that multiple zeroes are inevitable for all fillings
\ $Z_m$\,, \ since some \ $m_0$\,, \ that is when \ $m > m_0$\,. \ At
the same time finiteness of a sum of reverse modules for a considered
system shows, that there will be value \ $r$\,, \ for which
certainly correctly this expression
$$
SS \subset SS_{(II)}\,: \qquad \sum_{i=1}^{\infty}
\ \frac{1}{a_i} \ < \ \infty \quad \Rightarrow \quad \exists\,(r)\left\{
\sum_{i=r + 1}^{\infty} \ \frac{1}{a_i} \ < \ 1\,, \qquad n_0 \ > \
1 \right\}\,.
$$

Or else, the system of grids \ $\{S(a_i)\}$ \ at \ $i > r$ \ appears
by a system of the first type \ $SS_{(I)}$\,. \ But then pursuant to the
previous theorem exception the first \ $r$ \ grids from the initial system
or the inclusion the same \ $r$ \ grids in a system of the first type does
not change final one-sided algorithms of fillings. That is there will be
such value \ $n_0$\,, \ that at all parameters \ $n > n_0$ \ effective will
be, as well as in the expression (5.8), one-sided algorithm of filling
\ $Z_n^{(1)}$\,. \hfill$\Box$
\vspace*{3mm}

Sense of the introduction of one-sided concept of algorithm of filling
clears up the following statement, important for the analysis of the
characteristics of all systems.
\vspace*{3mm}

{\bf Theorem 11.} Unit acts by a majorizing coefficient \ $\tau_n$ \
valuations of a series length \ $msr_n (q) $ \ in the main theorem for
maximum series \ $MSR_n (q) $\,, \ generated by algorithm of one-sided
filling \ $Z_n^{(1)} \subset SS$\,.
$$
MSR_n (q,\,Z_n^{(1)},\,SS)\,: \qquad msr_n (q) \ \le \ \tau\,\frac{n + q}{1
- \gamma_n} + 1\,; \quad \gamma_n \ = \ \frac{H_n}{PZ_n}\,,
\quad \tau_n \le \tau = 1\,. \eqno{(5.9)}
$$
\vspace*{1mm}

{\sl Proof.} Here \ $\gamma_n$ \ there is the density (average frequency)
of zeroes \ $H_n$ \ for a period \ $PZ_n$ \ of filling \ $Z_n$\,. \ A
one-sided algorithm of filling \ $Z_n^{(1)}$ \ at formation of maximum
series \ $MSR_n (q) $ \ imposes on a system \ $SS$ \ and concrete modules
of grids rather severe restrictions. They are described by conditions
of the theorems 4--7. If not dependence between modules of grids, under all
other conditions the two-sided filling, would seem, is obliged to ensure
just the heaviest series. But the modules of next grids do not permit
to run in a double interval of a prospective series of zeroes.

The particular appendix of this situation degree-systems and systems
without multiple zeroes follows from the theorems 1--3. The system with
primary increase of modules requires the special consideration. We shall
offer majorizing variant of such system, providing minimum of mutual
prime modules of grids without grids of saturation. Then the modules of
similar grids are consistently and unequivocally: \
$3,\,7,\,16,\,37,\,79,\,...$\,. \ Naturally, the equal modules are also
excluded, as they do not change the characteristic and properties of a
system.

Redundancy of valuation of a maximum series \ $MSR_n(q)$ \ for initial
parameters \ $n$ \ at all \ $q \ge 0$ \ is rather simple directly to
establish. \ Expediently \ to \ consider \ $2$-filling and heaviest interval,
where will not meet a multiple zero at maximum concentration of zeroes of
multiplicity unit. As far as \ $\gamma_2 = \frac{3}{7}$\,, \ the heaviest
growth of a series and valuation will appear at \ $msr_2(4) = 11\,; \
\frac{6}{1-\gamma_2} + 1 = 11.5$\,. \ And then on intervals \ $I_2 \ge 21$
\ zeroes of multiplicity two are inevitable, and the difference between
valuation and series grows.

And let such growth is slowed down, but it is inevitable, and transition to
next filling \ $Z_3$ \ with a grid \ $S(16)$ \ at all desire to achieve
heaviest concentration of zeroes can grant only variant \ $msr_3 (8) = 20$\,,
\ when \ $\frac{11}{1-\gamma_3} + 1 \sim 21.5$\,, \ as far as \ $\gamma_3
= \frac{13}{28}$\,. \ The so appreciable difference is connected that
on this interval a multiple zero is inevitable. And further the quantity
of multiple zeroes accrues, even more increasing gap between estimation
and series.

If to take into account, that the relation between modules is saved
(about two), occurrence of new multiple zeroes is inevitable also. It
promotes a further divergence between valuation and maximum series of
zeroes. Thus completely it is necessary to allocate complete equality
of multiple zeroes. They in a identical degree cause fall of value of
the heaviest series.

It is impossible not to note, the consideration of maximum series \
$MSR_n(0)$ \ without units $(q=0)$ sharply simplifies the proof of the
statement. The essential part of grids forms density of zeroes \
$\gamma_n$\,, \ but participates in creation of a series of zeroes \
$MSR_n(0)$ \ by only  one zero, being compared with a infinite grid.
The statement does not deform and occurrence in a system of grids of
equal modules, already not speaking about liquidation of the rule of
their mutual prime.

Thus, one-sided algorithm of fillings \ $Z_n^{(1)}$ \ causes a
establishment of a majorizing coefficient \ $\tau = 1$ \ in the
formula (5.9) of main theorem. It precisely determines thus rather
extensive class of systems. \hfill$\Box$
\vspace*{3mm}

Now all is prepared for distribution of the statement of the main theorem
for all classes of non-singular systems of grids and fillings.
\vspace*{5mm}

\begin{center} {\Large \bf 6. Proof of the main theorem} \end{center}
\vspace*{3mm}

The imaging principle of explains logic and validity of origin
of the main theorem in the offered form, but does not remove
necessity of its proof. The series of such proofs is found in
monograph [1], using various feature of fillings method,
including absence of imaging principle's support. Finishing
variant giving the strongest valuation is indicated here.

We are addressed to three coprime \ $VP$ \ systems of the third
type and one mixed \ $SM$ \ system, playing central role in
some problems of primes distribution:
$$
\begin{array}{ll} SP_0 \ = \ \{2,\,3,\,5,\,7\,,...,\,p_i\,,...\}\,, & SP_1
\ = \ \{3,\,5,\,7\,,...,\,p_i\,,...\}\,, \\
SP'_1 \ = \ \{3,\,4,\,5,\,7\,,...,\,p_i\,,...\}\,, &
SM'_{04} = \,\{2,\,3,\,4,\,5,...,p_i,...\}.  \end{array}
\eqno{(6.1)}
$$

All these related systems at any \ $n > 2$ \ and united \ $q \ge 0$ \ form
fillings \ $Z_n [SS]$\,, \ which have the dependent characteristics of
series.
\vspace*{3mm}

{\bf Lemma 1.} The values of all series \ $SR_n [SS] = SR_n (q)$ \ in
periods of fillings \ $Z_n [SS]$ \ of systems of expression (6.1) are
connected by next equalities:
$$
sr_{n+1} [SP_0] \ = \ 2\cdot sr_n [SP_1]\,, \quad sr_{n+2}
[SM'_{04}] \ = \ 2\cdot sr_{n+1} [SP_0] \ = \ 4\cdot sr_n
[SP_1]. \eqno{(6.2)}
$$
\vspace*{0mm}

{\sl Proof.} The grid \ $S(2)$ \ for system \ $SP_0$ \ and grids \
$S(2)\,, \ S(4)$ \ for system \ $SM'_{04}$ \ are saturated
grids. It means, that all series of their fillings (including
period with corresponding number of units \ $q \gg 0$) \ are increased,
accordingly, in two and four times in
relation to series of filling of system \ $SP_1$\,.

Certainly, this conclusion is valid and for maximum series \
$MSR_n$\,. \ From here borders possible find for values of
the maximum series for important coprime system \ $SP'_1$\,. \ This
system have grid \ $S(4)$ \ and \ $(4,\,2) \neq 1$ \ for grid \
$S(2) \subset SM'_{04}$\,.
$$
\frac{msr_{n+1} [SP_0]}{msr_n [SP_1]} \ = \ 2\,, \quad
\frac{msr_{n+2} [SM'_{04}]}{msr_n [SP_1]} \ = \ 4\,, \quad
\frac{msr_n [SP'_1]}{msr_n [SP_1]} \ < \ C\,, \eqno{(6.3)}
$$
where \ $1 < C < 2$\,. \ First two equalities (6.3) follow direct from
equalities (6.2), and third inequality follows from obvious reasons, as far
as the addition of the same grid \ $S(2)$ \ to fillings \ $Z_n [SP'_1]$ \ and \
$Z_n [SP_1]$ \ leads them into fillings \ $Z_{n + 1} [SM'_{04}]$ \ and \
$Z_{n + 1} [SP_0]$ \ accordingly. It means limitation of value \
$C < 2$\,. \hfill $\Box$
\vspace*{3mm}

{\bf Theorem 12.} If non-singular systems \ $SS_1\,, \ SS_2$ \ differ by
finite set of grids, the correlation of values of their maximum series
has constant boundaries:
$$
n_0 \ge 1\,, \qquad C_1 \ < \ \frac{msr_{n+n_0}
[SS_1]}{msr_{n+n_0} [SS_2]} \ < \ C_2\,, \qquad \forall \ n >
0\,. \eqno{(6.4)}
$$
\vspace*{1mm}

{\sl Proof.} The values of constants \ $C_1 < C_2$ \ are
determined by concrete set of non-coincide grids. Their quantity
gives value \ $n_0$\,, \ that permits some to decrease a
difference \ $C_2 - C_1$\,. \ Inequalities (6.4) become clear
after addition of grids (optionally coincide) in each of systems
\ $SS_1\,, \ SS_2$\,, \ that sets of discrepancy have formed
saturated fillings. Quantity of such additions is finite and in
a result systems will be formed, distinguished only by sets of
saturated grids with periods \ $D_1$ \ and \ $D_2$\,. \ In such
case we receive values of maximum series \ $D_1 \cdot msr_n
[\overline{SS}]$ \ and \ $D_2 \cdot msr_n [\overline{SS}]$\,, \
where \ $\overline{SS} = SS_1 \cap SS_2$ \ is common grids part of
initial systems. \hfill $\Box$
\vspace*{3mm}

{\bf Definition 27.} {\sl Criterion of mixing} \ $K_n(SS)$ \ for filling
\ $Z_n$ \ of any non-singular system \ $SS \ = \ \{S(a_i)\}_1^{\infty}$\,,
\ where \ $(Z_n \ \subset \ \{S(a_i)\}_1^n \ \subset \ SS)$ \ is expression
$$
\frac{\gamma_n^*}{\gamma_n} \ = \ K_n(SS) \ \geq \ 1\,, \qquad
\gamma_n^* \ = \ \frac{H_n^*}{PZ_n^*}\,, \quad \gamma_n \ = \
\frac{H_n}{PZ_n}\,, \quad Z_n \ \subset \ \{S(a_i)\}_1^n\,. \eqno{(6.5)}
$$
\vspace*{1mm}

One would think, such criterion can act ratio \ $\frac{1 -
\gamma_n}{1 - \gamma_n^*} \ = \ \frac{PZ_n^*}{PZ_n}$ \ or \
$\frac{H_n^*}{H_n}$\,, \
however they are non-informative, as far as they characterize
not mixing, but availability of multiple zeroes, their plenty.
These values will be maximum for saturated fillings for system \
$\{(p_i - 1) \cdot S(p_i)\}$\,, \ where \ $p_i$ \ are primes.
The criterion \ $K_n(SS)$ \ of expression (6.5) is intended to
allocate that fact, that mixing is simplified step to systems
\ $SS'$ \ without multiple zeroes.
\vspace*{3mm}

{\bf Lemma 2.} The system \ $SP'_1 = \{ 3,4,5,7,..., p_i,...\}$\,,
\ where \ $p_i$ \ are primes, corresponds to majorizing sequence of
mixing criteria, that is relation
$$
\sup\limits_{SS} \ \max\limits_{\gamma_n^1 \sim \gamma_m(SS)} \
K_m(SS) \ = \ K_n(SP_1')\,, \quad \gamma_n^1 \ \in \ Z_n^1 \
\subset \ SP_1'\,, \quad \gamma_m \subset SS. \eqno{(6.6)}
$$
\vspace*{2mm}

The formulation of lemma means, that the upper estimates of
maximum series values, found with the help of the analysis of
zeroes frequencies (in particular, main theorem) for \
$n$-fillings of the system \ $SP_1'$\,, \ will be valid for any
other system.

Really, any inclusion in \ $n$-filling of grids with modules \
$(a_i,\,a_j) \neq 1$ \ means decrease of zeroes of multiplicity
higher unit, as far as such grids can be considered in this
filling as one grid with frequency \ $(\frac{1}{a_i} +
\frac{1}{a_j})$\,.

The address to the system \ $SP_1'$ \ is predetermined by
step-by-step consideration of \ $n$-fillings. Let \ $n = 2$\,, \
and the modules \ $a_1,\,a_2$ \ at coprime should be (as well as
above) close for maximum of \ $K_2(SS)$\,. \ In such case
setting \ $a_1 = m\,, \ a_2 = m + 1$\,, \ we shall receive \
$K_2 = 1 + \frac{m - 1}{2(m^2 + m + 1}$\,, \ whence follows,
that the maximum of this ratio is reached for module \ $m =
3$\,. \ The similar maximum for fillings at \ $n = 3$ \ takes
place for grids with initial modules \ $3,\,4,\,5$\,.

Further largely the coprime of modules enters, minimum which
is provided by primes \ $p_i$\,. \ Thus there is generated
system \ $SP_1'$\,, \ for which relation (6.6) is executed. It
is possible to note, that for \ $\gamma_n^1 \to 1$ \ the ratio
(5.5) aims to unit. Therefore value \ $K_n$ \ has maximum, which
is reached in a system \ $SP_1'$ \ for \ $n = 5$ \ and is equal
\ $K_5 = 1 + \frac{ 36456}{325367} \sim 1.11205$\,. \ Thus it
appears \ $K_4 \sim 1.1106\,, \ K_6 \sim 1.1116$\,.

The existence of maximum \ $\max\limits_n \ K_n$ \ does not
contradict that \ $\gamma_n^* \ge \gamma_n$ \ and \
$\gamma_n^* - \gamma_{n-1}^* \ge \gamma_n - \gamma_{n-1}$ \
for any system, and if at such transitions fresh multiple zero
will be fixed, the inequalities should be replaced to strict.

So, according to construction, as well as from definition of
mixing criterion \ $K_n(SS)$\,, \ if \ $\gamma_n^1(SP_1') \sim
\gamma_m(SS)$\,, \ when are close direct (average) density of
zeroes for two systems, the mixing criterion \ $K_n(SS) =
K_n(SP_1')$ \ will be more for \ $n$-filling of
system \ $SP_1'$\,. \ Clearly, approximation of density nearness
dictates and some approximation of criterion advantage, as far
as the compared fillings can differ rather slightly for many, if
not to all parameters.

Unfortunately, natural requirement of density equality for two
fillings incorrectly because from various sets of modules.
However for close density and at a essential divergence in sets
of grids, advantage of value \ $K_n(SP_1')$ \ will be without
fail displayed. Such advantage will become obvious, if the
appreciable part of modules of filling grids for system \ $SS$ \
will appear not coprime.

Therefore further system \ $SP_1'$ \ will be considered as
determining system of limiting concentration of multiple zeroes
concerning zeroes \ $H_n$ \ of direct imaging.
It is received by withdrawal of saturated grid \ $S(2)$ \ and
then inclusion of grid \ $S(4)$ \ with a minimum
even module, large two. It should add, that practically
such majorizing system (rather close), is standard prime system
\ $SP_1$\,. \hfill $\Box$
\vspace*{1mm}

Systems of grids with modules -- primes \ $SP_0$ \ and \ $SP_1$
\ are determining on a way of the proof of the main theorem. Or
else, if the theorem is valid for these systems, it is valid and
for a class of all non-singular systems. It is called by that
coprime system \ $SP_1'$ \ appears by a limiting system
according to mixing criterion (lemma 2), and it only by one grid
differs from mentioned prime systems.
\vspace*{3mm}

{\bf Definition 28.} Next primes, connected by equality \ $(r \geq 1):$ \
$p_n = p_{n-1} + 2 r$\,, \ are
named as {\sl kinsfolk of rank} \ $r \in {\bf N}$\,. \ We
designate thus \ $p_n = p_n^{(r)}$\,.
\vspace*{3mm}

Thus, twins are kinsfolk of first rank. Clearly, that the search
kinsfolk very large rank produces to significant difficulties.
As the kinsfolk is determined large: \ $p_n^{(r)} = BR_s^r = p_n
= p_{n-1} + 2r$\,, \ then all primes (except 2 and 3) are
kinsfolk of one from ranks. Next kinsfolk of various ranks can be
also incorporated and are considered as independent object.
\vspace*{3mm}

{\bf Definition 29.} {\sl Configuration } of \ $m\,, \ (m \geq
2)$ \ primes is named group of next primes as vector kinsfolk of
ranks \ $r_i$ \ of dimension \ $m - 1\,: \
(r_1,\,r_2,\,...,\,r_{m-1})$\,.
$$
p_n^{(r_1)} = p_{n-1} + 2r_1\,, \ p_{n+1}^{(r_2)} = p_n +
2r_2\,,...\,,\,p_{n+m-2}^{(r_{m-1})} = p_{n+m-3} + 2r_{m-1}\,.
\eqno{(6.7)}
$$
\vspace*{1mm}

At the same time quite clearly, any vector \ $\{r_i\}^m$ \ of
dimension \ $m - 1$ \ at \ $m \ge 3$ \ does not guarantee, that
there will be an appropriate configuration. For example, for a
vector \ $(1,\,1,\,1)$ \ the configuration of primes of kind
(6.7) does not exist.
\vspace*{3mm}

The problem about the upper estimate of the maximum series \
$MSR_n$ \ or \ $MSR_n(q)$ \ acts central in the fillings
method. Thus majorant of systems \ $SP_1$ \ and \ $SP'_1$ \ in
class of all systems by mixing criterion acquires decisive
character. Therefore main investigated system will become just
the system \ $SP_1$\,, \ though for the researchers, not aware
about fillings method, always unique was the system \ $SP_0$\,,
\ in which so it is conveniently to build sieve of Eratosthenes.

Object of the fillings method is whole period of grids product,
and sometimes the study is not limited even by period. It is
testified already repeatedly, that the period of \ $n$-filling
in \ $SP_1$ \ has a length, equal to product of all \ $n$ \ odd
prime.  Period disintegrates by series of zeroes \ $SR_n$\,, \
the lengths of which vary from values \ $sr_n = 1$\,, \ (there
are no zero between units), up to value \ $msr_n$\,.

In illustrations series \ $MSR_n (q)$ \ and sequences of series \
$SR_n$ \ are submitted and the periodicity has allowed them
to close in a ring (Fig. 1). According to lemma 1 transition
from a system \ $SP_0$ \ to \ $SP_1$ \ means transformation
kinsfolk of rank \ $r$ \ in series of this length, and
configuration of primes in configuration of ranks.

Fig. 1 from "Graphic Illustrations" clearly demonstrates symmetry
of filling series \ $Z_3$\,, \ which is present for any \
$Z_n \subset SP_1$\,. \ Thus each series of the period of filling
has double, except two series \ $SR_n^{(I)}\,, \ SR_n^{(II)}$\,, \
submitted in single specimen -- series of length of unit and two.
We shall designate symmetric axial configurations of series with
these series in center as \ $Kf_n^{(I)}$ \ and \ $Kf_n^{(II)}$\,.

Each grid \ $S(p_k) \in Z_n$ \ has central series, consisting
from two next units, that is length of unit. It clearly, as far
as any period of grid \ $S(p_k)$ \ includes \ $p_k~-~1 \ge 2$ \
units. Thus the series \ $SR_n^{(I)}$ \ unit length of the first
axial configuration \ $Kf_n^{(I)}$ \ will be saved for each
step. It can be concluded from the obvious relation
$$
\frac{1}{2}\left\{p_n\,\prod_{i=1}^{n-1} p_i \ \pm \ 1\right\} \
\not\equiv \ 0 \ ({\rm mod} \ p_n) \ \Longrightarrow \
\frac{1}{2}\,(p_n \ \pm \ 1) \ \not\equiv \ 0 \ ({\rm mod} \ p_n)\,,
$$
and serial unit \ $SR_n^{(I)} $\,, \ lying equally in distance of a
half-period from central series \ $SR_n^{(II)}$ \ of the second
configuration, remains in constancy.

In the period \ $PZ_n$ \ length \ $\prod_{i=1}^n \ p_i$ \
there is unique zero of multiplicity \ $n$\,, \
received by product of zeroes of all grids. It enters in
central axial series \ $SR_n^{(II) }$ \ of length \ $sr_n^{(II)}
= 2$\,. \ Clearly, that all other zeroes of grids lie symmetric
concerning mentioned \ $n$-th zero. Hence, and all remaining
series of the configuration \ $Kf_n^{(II)}$ \ place symmetric,
except one central series \ $SR_n^{(I)}$\,.

The clear sense has consideration of those sequences of next
series, which will meet too in subsequent fillings at increase
of parameter \ $n$\,. \ Such configurations are named typical.
The first typical axial configuration \ $Kf^{(I)}(Z_n) \ = \ Kf^{(I)}_n$
\ is submitted uniquely by series in this scheme:
$$
Kf^{(I)}_n \ = \ Kf\left\{...,\,\frac{p_{n+2} -
p_{n+1}}{2},\,\frac{p_{n+1} - 1}{2},\,1,\,\frac{p_{n+1} -
1}{2},\,\frac{p_{n+2} - p_{n+1}}{2}\,,...\right\}, \eqno{(6.8)}
$$
where central series \ $SR_n^{(I)}$ \ of unit length \ $(sr_n^{(I)} = 1)$
\ surround by two series of the greatest length for this configuration
\ $Kf^{(I)}_n$\,.

We shall put in conformity zero of axis and average point of two
central units of each grid. \ Then representation of configuration (6.8)
is proved by that in points of projection \ $Kf^{(I)}_n \in {\bf Z}$\,: \
$\pm \left(\frac{p_{n -j} + 1}{2} \ + \
k\,p_{n-j}\right)$ \ zeroes stand and at \ $p_{n + i} \ < \ p_{n + 1}^2$\,:
$$
\frac{p_{n -j} + 1}{2} \ + \ k\,p_{n-j} \ \neq \ \frac{p_{n +i} +
1}{2}, \quad \forall\,j,i,k: \ 0 \le j \le n-1, \ i,k \ge 1. \eqno{(6.9)}
$$

The second axial configuration \ $Kf^{(II)}_n$ \ is submitted by next series:
$$
Kf^{(II)}(Z_n) \ = \ Kf\{\,
...,\,p_{n+3} - p_{n+2},\,p_{n+2} - p_{n+1},\,p_{n+1} - 2^m,\,2^{m-1},...
$$
$$
...,\,4,\,2,\,1,\,2,\,1,\,2,\,4,...,\,2^{m-1},\,p_{n+1} - 2^m,\,p_{n+2} -
p_{n+1},\,p_{n+3} - p_{n+2},...\,\}\,, \eqno{(6.10)}
$$
where \ $m = [\,\log_2\,p_{n+1}\,]$\,, \ that is \ $m =
\max\,\{i: \ 2^i < p_{n+1}\}$\,, \ and typicalness of each
configuration, dependent from \ $n$\,, \ should be established
especially. Central series \ $SR_n^{(II)}$ \ of configuration \
$Kf_n^{(II)}$ \ has length two \ $sr_n^{(II)} = 2$\,. \ Besides units,
correspond to numbers \ $2^{m + k}\,, \ k \ge 1$\,, \ will meet
in configuration.

We shall put in conformity zero of axis \ ${\bf Z}$ \ and zero of each
grid. Then representation (6.10) is proved by that in points \
$\pm 2^r\,, \ r \ge 0$ \ of projection \ $Kf^{(I)}_n \in {\bf Z}$ \
units stand and by analogy with (6.9) we have: \ $k\,p_{n - j} \
\neq \ p_{n + i}$\,.

The maximum series \ $MSR_n \equiv MSR_n (0)$ \ can
enter in a core of typical axial configuration, and
it is not (remaining by a typical series), but in virtue
of indicated expressions it is possible to make the
conclusions about its value.
\vspace*{3mm}

{\bf Theorem 13.} Values \ $msr_n (0)$ \ of maximum series \ $MSR_n (0)$
\ in prime systems \ $SP_0 \supset \{S(p_i)\}_0^n$ \ and \ $SP_1 \supset
\{S(p_i)\}_1^n$ \ have following lower estimates
$$
\{\,msr_n(SP_1) \ \ge \ p_{n-1}\,\} \ \ \Rightarrow \ \
\{\,msr_{n+1}(SP_0) \ \ge \ 2\,p_{n-1}\,\}, \qquad  n \ge 1.
\eqno{(6.11)}
$$
\vspace*{1mm}

{\sl Proof.} The second inequality (6.11) is consequence of the
first according to lemma 1. It easily be convinced in validity
of given equality for parameter \ $n \le 7$\,: \ \ $msr_n(SP_1) \ = \
p_{n-1}$\,. \ For parameter \ $n > 7$ \ the equality begins to be
infringed. According to
representation (6.8) value of series \ $SR_n(2)$ \ for \ filling
\ $Z_n$ \ with two units in system \ $SP_1$ \ is equal \
$sr_n (2) \ = \ \frac{p_{n+1} - 1}{2}\,+\,1\,+\,\frac{p_{n+1} - 1}{2}
= p_{n+1}$\,. \ It is received from three central series
of configuration. As far as for the same series \ $sr_n(2) \le
sr_{n + 2} (0)$ \ and the maximum series majorize of any, the
statement of the theorem follows from replacement of parameter \
$n$ \ to value \ $n-2$\,. \hfill $\Box$
\vspace*{3mm}

Appeal to second \ $Kf^{(II)}_n = Kf^{(II)}(Z_n)$ \ axial
configuration demonstrates, that lower estimate of maximum
series is unjustifiable rough at sufficiently large \ $n$\,.
\vspace*{3mm}

{\bf Theorem 14.} Lower estimate of maximum series \
$MSR_n(0)$ \ in prime system \ $SP_1$ \ for \ $n \ge 26$ \ surpasses
estimation \ $p_{n-1}$\,, \ as it is expressed by the formula
$$
msr_n(SP_1) \ \ge \ 2 \cdot p_{n - 2m -1}\,, \qquad
m \ = \ [\,\log_2\,p_{n-2m-1}\,]\,, \qquad n \ge 26\,. \eqno{(6.12)}
$$

{\sl Proof.} Here the expression for search of intermediate
parameter \ $m$ \ provides the decision of a small integer
equation, which however can not call difficulties and does not
in essence change common kind of estimation. We shall consider
expression \ $Kf^{(II)}(Z_n)$ \ of second axial configuration
(6.10). Values \ $sr_n (q)$ \ of central series \ $SR_n(q)$ \ with
\ $q$ \ units hence it follow immediately:
$$
sr_n (2m) \ = \ 2^{m+1}\,, \qquad m =
[\,\log_2\,p_{n+1}\,]\,, \quad p_n = \max\,\{p_i: \ p_i < 2^m\}\,;
$$
$$
sr_n (2m + 2) \ = \ 2\cdot p_{n+1}\,; \qquad \ sr_{n+k} (2m + 2)
\ = \ 2\cdot p_{n+k+1}\,, \quad k \ge 1\,. \eqno{(6.13)}
$$

From these expressions follows, that as far as any grid can
zerofill not less units, than infinite, for \
$(n+k+2m+2)$-filling with grids \ $S(p_i)$ \ is always executed
$$
\{\,sr_{n+k+2m+2} (0) \ \ge \ 2\cdot p_{n+k+1}\,\} \quad
\Rightarrow \quad \{\,sr_n (0) \ \ge \ 2\cdot p_{n-2m-1}\,\}\,,
\eqno{(6.14)}
$$
and the last inequality is received after replacement of \
$(n+k+2m+2)$ \ to \ $n$\,. \ It is thus necessary to take into
account, that value \ $m$ \ is found for parameter \
$(n+k+2m+2)$ \ instead of \ $n$\,. \ It predetermines
necessity of the equation decision for search \ $m$ \ in the
formulation of theorem.

The concrete check [1] of axial configurations \ $Kf^{(I)} \subset SP_1$
\ and \ $Kf^{(II)} \subset SP_1$ \ finds out the first and minimum
value \ $n = 26$ \ for which is executed
$$
sr_{n - 2q} (2q)\,[Kf^{(II)}] \ \ > \ \ sr_{n - 2}
(2)\,[Kf^{(I)}]\,, \qquad q \ = \ q(n)
$$
for the greatest series with units of corresponding axial
configurations. Then for \ $n > 26$ \ the sign of inequality
does not already change. Twin quantity of units \ $q > 1$ \ is
determined from minimum and symmetry conditions according to
which the greatest initial series of the configurations should
enter in investigated series.

It is simple find, that limit of the lower estimates of theorems 13 and 14
is equal two at increase of filling parameter \ $n \to \infty$ \ in
prime system \ $SP_1$\,:
$$
\lim_{n \to \infty} \ \frac{2 \cdot p_{n - 2m -1}}{p_{n - 1}} \
= \ 2 \ \lim_{n \to \infty} \ \left( 1 - \frac{c_1 \ln
n}{n} \right) \left( 1 - \frac{c_2}{n} \right) \ = \ 2\,.
\eqno{(6.15)}
$$

So, the lower bound ( 6.12 - 6.15) of maximum \ $MSR_n (0)$ \
series for large \ $n$ \ qualitatively surpasses similar values
for small filling parameters \ $n$\,. \hfill $\Box$
\vspace*{3mm}

However essentially greater significance for the subsequent
research had upper estimates of maximum series \ $MSR_n$ \ and \
$MSR_n(q)$\,. \ It is better to have exact values \ $msr_n(q)$ \
for all \ $q \ge 0$\,. \ There is, at all complexity and
importance of this problem, it is solvable within the framework
of the fillings method for system \ $SP_1$ \ and the following
statements serve necessary step for it.
\vspace*{3mm}

{\bf Theorem 15.} Greatest series \ $SR_n (q)$ \ with \ $q \ge 0$
\ units belonging to core of the first \ $Kf_n^{(I)} =
Kf_n^{(I)}(Z_n)$ \ axial configuration in a system \ $SP_1$ \
have the length
$$
SR_n (q) \subset Kf_n^{(I)}: \qquad sr_n (0) = \frac{p_{n+1} -
1}{2}\,; \ \qquad sr_n (1) = \frac{p_{n+2} - 1}{2}\,;
$$
$$
sr_n (2) = p_{n+1}\,; \quad sr_n (q) = \max_{1 \le i \le [q/2]}
\left\{p_{n+i} + \frac{p_{n+q-i} - p_{n+i}}{2}\right\}, \quad q
\ge 3\,. \eqno{(6.16)}
$$
\vspace*{1mm}

{\sl Proof.} The representation of expression (6.8) of the first axial
configuration \ \ $Kf_n^{(I)} = Kf_n^{(I)}(Z_n)$ \ quite
determines and fixes relations of theorem. It should take into
account, that values of the first four series are given by a
kind of the configuration directly, and their receipt does not
require in search, as for \ $q = 0;\,2$ \ it obviously, for $q
= 3$ \ follows from symmetry, and \ $sr_n (1) = \frac{p_{n+2} -
1}{2} \ge \frac{p_{n+1} + 1}{2}$\,. \ But also linear search
during simple selection of maximum does not result to large
retrieval of values, as far as it is connected with local
non-uniformity of primes distribution for a sequence of indexes from \
$n + 1$ \ up to \ $n + q - 1$\,.

The length of series with \ $q$ \ units of axial configuration
is equal to sum of values of making series, separated by \ $q$ \
commas. So, the length of series \ $SR_n (3)$ \ with three
units, which contains both greatest series, is equal \ $p_{n+1}
+ \frac{p_{n+2} - p_{n+1}}{2}$\,. \ According to condition of
the theorem, such series is greatest, as far as its length is
unique. Already for \ $q = 4$ \ the situation changes.

Really, in this case according to condition has to choose in
expression (6.16) from two variants of series \ $SR_n(4)$ \ of lengths \
$sr_n (4)$
$$
SR_n (4)\,: \qquad \left\{\,sr_n (4) \ = \ p_{n+2}\,; \quad
sr_n' (4) \ = \ p_{n+1} \ + \ \frac{p_{n+3} - p_{n+1}}{2}\,\right\}.
$$

Depending on \ $n$ \ the advantage can have this or that
variant. For example, for \ $n=2$ \ maximum \ $SR_2(4)$
\ determines by first variant (as \ $11 > 10) $ \ and at \
$n=3$ \ for \ $SR_3(4)$ -- by second: \ $(13 < 14)$\,. \ From
here we have \ $msr_2(4) = 11, \ msr_3(4) = 14$\,.

From unique construction of axial configurations for any \ $n$
\ and condition of choice of the greatest series follows, that
quantity \ $[q/2]$ \ exists in common case of various series \
$SR_n(q)$ \ as claimants for a role maximum. In common case them
lengths of these series can not coincide too. It is thus necessary to
take into account typical configurations and value \ $q$ \ is limited
naturally. From here search of variants of maximum series is small,
estimated and obvious.

The given formula (6.16) acts for values \ $q$\,, \ satisfying \
inequality \ $p_{n+q-1} < p_{n+1}^2$\,. \ It is obligatory condition for
system \ $SP_1$ \ and fillings method. For large \ $q$ \ there is reminder,
that is considered regulated \ $n$-filling, instead of primes
distribution. That is in the configuration \ $Kf_n^{(I)}(Z_n)$ \
there will be units, not corresponding primes, and last search
it is necessary to change. The theorem is proven. \hfill $\Box$
\vspace*{3mm}

{\bf Theorem 16.} Greatest series \ $SR_n(q)$ \ with \ $q \ge 0$ \
units belonging to core of the second \ $Kf_n^{(II)}$ \ axial
configuration in a system \ $SP_1$ \ have the length
$$
SR_n (q) \ \subset \ Kf_n^{(II)}\,: \qquad sr_n (q) \ = \
\max_k\,\left\{\,\sum_{i=1}^q sr_n (0)\,[k+i]\,\right\}\,,
\eqno{(6.17)}
$$
where series \ $SR_n (0)\,[j] \subset Kf_n^{(II)}$ \ of corresponding
lengths for consecutive parameters \ $j = k+i$ \ are next series of the
configuration \ $Kf_n^{(II)}$\,.
\vspace*{3mm}

{\sl Proof.} Some indeterminacy of task of parameter \ $k$ \ at
search of the greatest series \ $SR_n(q)$ \ is removed by that
at small \ $q$ \ in series one of two groups of greatest initial
series \ $SR_n(0)\,[k+i]$ \ is obliged enter and at \ $q \ge 2m$
-- both groups. It is enough obviously from representation of
the configuration. Besides as well as in the theorem 7 search is
conducted so long as \ $p_{n + q-1} < p_{n+1}^2$\,.

From the formula of theorem (6.17) directly follows, that
$$
sr_n (0) \ = \ \max\{2^{m-1},\,p_{n+1} - 2^m\}\,; \
sr_n (1) \ = \ \max\{3\cdot 2^{m-2},\,p_{n+1} - 2^{m-1}\}\,;
$$
$$
sr_n (2) \ \ = \ \ \max\,\{\,7\cdot 2^{m-3}, \ p_{n+1} -
2^{m-2},\,p_{n+2} - 2^{m-1}\,\}\,;\,...
$$
$$
...,\,sr_n (2m) = \max\{2^{m+1},\,2^{m-1} + p_{n+1}\}\,; \quad
sr_n (2m+1) \ = \ 2^m + p_{n+1}\,;
$$
$$
sr_n (2m+2) \ = \ \max\,\{\,2\cdot p_{n+1}, \ 2^m + p_{n+2}\,\}\,;\,...
$$
$$
...,\,sr_n (2m+w) \ = \ \max_{1 \le i \le [w/2]} \ \{\,2^m +
p_{n+w}, \ p_{n+w-i} + p_{n+i}\,\}\,, \quad w \ge 2\,.
$$

The decisive significance for further conclusions in research
has series \ $SR_n$ \ with quantity of units, reached and
exceeded border \ $2m$\,, \ when in object both greatest groups
of initial series of the configuration are involved. Value \
$sr_n(2m + 1)$ \ is deprived  of the search factor and is
calculated directly after the task of central parameter \ $n$ \ of
filling. Examples of increasing values of such greatest series for \
$q = 2m$ \ and \ $q = 2m + 2$ \ are in monograph [1].

In the second axial configuration as regulated filling \ $Z_n$ \
all elements corresponding to values \ $2^s$ \ for \ $0 \le s
\in {\bf Z}$ \ are units. It is not taken into account by last
relation for series \ $sr_n(2m+w)$ \ in which such unit of the
greatest number stands in the point \ $2^m$ \ near from series bound.
Therefore we shall continue representation of the second axial
configuration \ $Kf_n^{(II)}$\,.

In such case after the found parameter of degree \ $m =
[\log_2\,p_{n + 1}]$ \ of filling \ $Z_n$ \ for prime system \
$SP_1$\,, \ whence by minimum
value \ $n' \le n$\,, \ at which unit in point \ $2^m$ \ becomes
the representative of central and already constant core of
configuration, will be \ $n' = \min \{i: \ p_ {i+1} > 2^m\}$\,.
\ From here under the given scheme we shall find following
such unit of the configuration in point \ $2^{m+1}$\,, \ lying
directly after unit corresponding to prime \ $p_u\,, \ u =
\min \{ j: \ p_{j+1} > 2^{m+1}\}$\,. \ Now we shall present a
half-configuration (for shortening of notation):
$$
2,1,2,4,...,2^{m-1},p_{n+1} - 2^m,p_{n+2} - p_{n+1},...,2^{m+1}
- p_{u}, p_{u+1} - 2^{m+1},...
$$
in which and the subsequent such units are similarly. Received
specified kind of configuration should take into account during
search of its series.

However this kind of configuration representation \
$Kf_n^{(II)}$ \ will be infringed yet earlier than unit
corresponding to \ $p_{n+1}^2$ \ will meet. Regulated filling
of configuration except permanent units corresponding to values
\ $2^k$ \ has also consistently zerofilled (at $n \to n + s$)
units in place \ $2^k p_{n+s}\,, \ k,s \ge 1$\,. \ Naturally,
first such unit the place \ $2p_{n+1}$ \ determine. Alongside
with taken into account units of representation, these units of
configuration play essential role, appreciably complicating
formula \ $Kf_n^{(II)}$ \ at aspiration to expand observed set \
$q$\,.

But additional units of series interval can not be the factor,
promoting to increase of the greatest (maximum) series of filling
\ $Z_n$ \ for any system. \hfill $\Box$
\vspace*{3mm}

The theorems 15 and 16 permit to notice important distinction
between axial configurations \ $Kf_n^{(I)}$ \ and \
$Kf_n^{(II)}$\,. \ In first \ $(I)$ \ there is only one
permanent series \ $SR_0$ \ of length unit, but all units in
interval up to \ $p_{n + 1}^2$ \ correspond exclusively primes.
For second \ $(II)$ \ quantity of permanent units and
hence series increases, but in the same interval additional
units corresponding to composite numbers \ $2^k p_{n+s}$ \ meet.
Thus role as that, as other configuration in formation of the
most main objects of fillings is impossible overrate.
\vspace*{3mm}

{\bf Theorem 17.} In system \ $SP_1$ \ for parameter \ $n \ge 3$ \
of the filling \ $Z_n$ \ the greatest series from considered
axial configurations are series \ $SR^{(I)}_n(q)$ \ for parameter \
$0 \le q < m + 3$\,, \ but at exception of series parameter \
$q = 1\,: \ SR_n (1)$\,. \ Then, for \ $m + 3 < q < 2 m - 3$ \
indeterminacy of advantage is accompanied by proximity
of series values. For parameters \ $q > 2 m - 3$ \ the advantage
passes to series \ $SR^{(II)}_n(q)$ \ of the second
configuration \ $Kf_n^{(II)}$\,, \ and it quickly increases with
growth \ $n$\,:
$$
sr^{(I)}_n (0) \ \ge \ sr^{(II)}_n (0)\,; \quad  sr^{(I)}_n (1) \
< \ sr^{(II)}_n (1)\,; \quad sr^{(I)}_n (2) \ > \ sr^{(II)}_n (2)\,;
$$
$$
sr^{(I)}_n (3) \ \ge \ sr^{(II)}_n (3)\,; \ ... \quad ... \ sr^{(I)}_n
(q) \ > \ sr^{(II)}_n (q) - \Delta_n\,, \quad q < m + 3\,;
$$
$$
n > 15: \quad ...\,,\,sr^{(II)}_n (q) \ > \ sr^{(I)}_n (q)\,,
\qquad 2m - 3 < q < 3n\,. \eqno{(6.18)}
$$
\vspace*{1mm}

{\sl Proof.} The last condition \ $q < 3 n$ \ is given as the
plenty of units \ $q$ \ loses informative sense at complication
of representation of the second configuration. In this case
maximum series value begins promptly to approach with average
value of such series in period. Intermediate parameters \ $m + 3
< q < 2 m - 3$ \ are omitted from consideration, as for such
variants the greatest series of both configurations are far from
relative maximum. Certainly, the advantage of series of one of
configurations can be established for each concrete \ $q$ \ at
increase \ $n$\,, \ but the special necessity is not present, as
far as for \ $q = 2m,\,2 m + 2$ \ the series value \
$SR^{(II)}_n(q)$\,, \ that is value \ $sr^{(II)}_n(q)$\,, \ becomes
determining.

Unsteady advantage of series \ $SR^{(I)}_n(q)$ \ for parameter \
$q < m + 3$ \ is reflected by the introduction in corresponding
relation (6.18) of essentially small \ $\Delta_n > 0$ \
concerning series value. For parameter \ $q$ \ approaching to
bound \ $2m$\,, \ when in evaluated series all central units of
kind \ $2^k$ \ enter for \ $n > 15$\,, \ advantage of series of
the second configuration appears obvious.

For proof of the first inequality it is enough to compare
theorems 15 and 16 concerning the greatest series of
configurations, designated as \ $sr^{(I)}_n(0)$\,, \
$sr^{(II)}_n(0)$\,, \ and then to see, that
$$
2^{m-1} \ \le \ \frac{p_{n+1} - 1}{2}\,, \qquad p_{n+1} - 2^m \ \le \
\frac{p_{n+1} - 1}{2}\,.
$$
If to take into account relation \ $m = [\log_2\,p_{n + 1}]$\,,
\ we shall receive given conditions \ $2^m + 1 \le p_{n+1} \le
2^{m+1} - 1$\,. \ Both these variant take place for various
primes: \ $5,\,17\,; \ \ 3,\,7,\,31,\,127$\,, \ that
predetermines possible and attainable equality.

Following value \ $q = 1$ \ results to explicable inequality \
$3 \cdot 2^{m-2} \ > \ \frac{p_{n+2} - 1}{2}$\,, \ when value \
$p_{n + 1}$ \ is reasonably close to lower bound \ $2^m + 1$\,. \
Let it not so and \ $3 \cdot 2^{m-2} \ < \ \frac{p_{n+2} -
1}{2}$\,. \ Then for \ $p_{n + 1} = 2^m + R$ \ and
\ $\Delta = p_{n + 2} - p_{n + 1}$ \ we receive: \ $R > 2^{m-1}
- \Delta + 1$\,. \ However in this case according to theorem 8
is executed \ $p_{n+1} - 2^{m-1} \ > \ \frac{p_{n+2} - 1}{2}$\,,
\ for that it is enough \ $R > \Delta - 1$\,. \ This condition
follows from the earlier received assumption, as far as \
$2^{m-1} > 2 \Delta - 2$\,.

The parameter \ $q = 2$ \ is especially important for the first
configuration and whole further as uniting both its greatest
series. Thus is reasonably obviously executed
$$
p_{n+1} > \max\,\{\,7\cdot 2^{m-3}, \ p_{n+1} -
2^{m-2},\,p_{n+2} - 2^{m-1}\,\}\,, \quad p_{n+1} \ge 2^m + 1
$$
because of obligatory last inequality, which and proves
advantage of the first configuration for such \ $q = 2$\,. \ It
will be saved and for subsequent parameters \ $q$\,, \ the proof
of the statement for which similarly.

Such reasonably stable situation begins to change at approach
to the dependent value \ $q = 2 m$\,, \ when advantage passes to the second
configuration \ $Kf_n^{(II)}$\,. \ However this exclusively important fact
takes place for rather large values \ $n$\,, \ thus according to
large \ $m$ \ also. From theorem 6 follows, that the complete
definiteness arises for \ $n + q = 26$ \ and then in accordance
with growth of parameters \ $n$ \ and \ $m$ \ advantage of
series of the second configuration becomes decisive.

The given by theorem condition \ $n > 15$ \ is determined by that value \
$n - 2m$\,, \ where \ $m = [\log_2\,p_{n + 1}]$\,, \ there is more unit,
that comparison of corresponding series possessed necessary
efficiency. In other case advantage of series of the second
configuration not so obvious, if is generally present.

Really example without searched determined parameter
$q = 2m+1$ takes place
$$
sr_n^{(II)}(2m+1) = 2^m + p_{n+1}\,, \quad sr_n^{(I)} (2m+1) \ = \
\max_{1 \le i \le m} \left\{p_{n+i} + \frac{p_{n+q-i} -
p_{n+i}}{2}\right\}\,,
$$
and \ $sr_n^{(I)} (2m+1)$ \ can approximately estimate by value
\ $p_{n+m+1}$\,. \ Then superiority of the second configuration
is connected to a obvious inequality \ $2^m > p_{n+m+1} -
p_{n+1}$ \ growing in accordance with increase \ $m$ \ and not
always valid for initial \ $m$\,. \ As far as parameter \ $m$ \
is connected with \ $n$ \ logarithmically, it results in
reasonably large values \ $n + q$\,, \ at which firm advantage
of series \ $SR_n^{(II)}(q)$ \ comes.

The achieved advantage will be saved for reasonably large \
$q$\,, \ but thus the kind of the second configuration \
$Kf_n^{(II)}$ \ should be transformed by inclusion of units,
according to values \ $2^k\,p_{n + s}$\,. \ It compels to limit
observed set of parameters \ $q$\,. \hfill $\Box$
\vspace*{3mm}

Fillings method and proven theorems 13 -- 17 permit to formulate
exclusively important statement concerning maximum series in the
system \ $SP_1$\,, \ satisfying to mentioned condition of
absolute independence from imaging principle. Nevertheless it
does not mean non-necessity or mistake of it, and opposite it
independently confirms and it pays attention to universality
at appeal to any classes of systems.
\vspace*{3mm}

{\bf Theorem 18.} Values of maximum series \ $MSR_n$ \ and \
$MSR_n(q)\,, \ q > 0$ \ are bending for the greatest series of
configurations \ $Kf_n^{(I)}$ \ and \ $Kf_n^{(II)}$\,. \ At
initial \ $n$ \ they are expressed through the first series \
$SR_n(q)$ \ of fillings, then they are connected with series of
the first configuration, and in result at \ $q = 2m$ \ they are
already only derivative of greatest series of the second
configuration:
$$
MSR_n (q) \subset SP_1: \qquad \left\{\, Kf_n^{(II)} \ \supset \
SR_n^{(II)} (2m) \ = \ MSR_n (2m)\, \right\}\,, \eqno{(6.19)}
$$
where quantity of generated grids \ $n \ge 15$ \ and \ $m =
[\log_2\,p_{n+1}]$\,.
\vspace*{3mm}

{\sl Proof.} Determining role of greatest series of the second
configuration \ $Kf_n^{(II)}$ \ becomes absolute only for large
\ $n$\,. \ At initial and small \ $n$ \ the series of the second
axial configurations accept auxiliary and supporting
participation in formation of maximum series. In system \ $SP_1$
\ we shall observed formation of initial maximum series \ $MSR_n
= MSR_n (0)$ \ as the most important for reception of many
further conclusions. Moreover their values by simple relation
are connected to series \ $MSR_n (q)$ \ for reasonably wide
spectrum of parameters \ $q$\,.

Generation of maximum series of filling \ $Z_n \subset SP_1$ \
passes through three stages. At first (the zero stage)
consecutive zerofilling of next units in fixed interval results
to occurrence of the first series \ $SR_n^1 (0)$\,, \ which and
become maximum. However this stage is quickly finished. Already
for \ $n = 4$ \ the first infringement is observed and it
appears chronic at growth of parameter \ $n$\,.

As was specified above and it is consistently confirmed by
examples for various systems, the step-by-step filling next
(right) units by zero of following grids is algorithm, realizing
one of variants of the greatest series of filling in the given
interval. A little that, such process really determines maximum
series \ $MSR_n (0)$ \ for many systems (without multiple
zeroes) and for many fillings in any systems.

Algorithm of sequential filling of units by grids of increasing
modules in system \ $SP_1$ \ results to maximum series \ $MSR_n
\equiv MSR_n (0)$ \ for \ $1 \le n \le 3$\,. \ And really for
these \ $n$ \ equality \ $p_{n-1} = 0.5\,(p_{n+1} - 1)$ \ takes
place. However then advantage of maximum series over
the first series thus regulated filling begins to grow.

It is necessary to specify, that there will be such order of
grids product of concrete filling, which occurrence of maximum
series as the first series of considered interval provides, if
to use described algorithm of formation.

So greatest series \ $SR^1$ \ of axial configuration of filling
\ $Z_n$ \ of length \ $0.5\,(p_{n+1} - 1)$ \ already is enough
extended, to claim for the special attention. If to remind, that
it is rather close from it (through series \ $S_n^0$ \ of length
unit) symmetric places such series, the arisen series \
$SR_n(2)$ \ with two units \ $(q = 2)$ \ becomes object,
claiming for extremes of characteristics. This moment
determines transition to following first stage of generation of
maximum series, connected with series of the first axial
configuration \ $Kf_n^{(I)}$ \ in period \ $PZ_n \subset SP_1$\,.

According to theorem 15 greatest series \ $SR_n^{(I)}(2)$ \ of
the first configuration has length \ $sr_n^{(I)}(2) =
p_{n+1}$\,. \ It immediately gives the lower estimation of
maximum series value \ $MSR_n$\,, \ as far as \ $msr_{n+s} (q -
s) \ge sr_n (q)$\,, \ $s \ge 0$\,. \ In particular, we receive \
$msr_{n+2} \ \ge \ sr_n^{(I)}(2)$\,, \ whence follows that \
$msr_n \ \ge \ p_{n-1}$\,.

Concrete check has demonstrated, that such lower estimation of
maximum series \ $MSR_n$ \ is upper for parameters \ $0 \le n
\le 7$\,, \ that is \ $msr_n \ = \ p_{n-1}$\,. \ At the same
time maximum series satisfy \ $msr_8 \ > \ p_7$\,, \ $msr_{11} \ > \
p_{10}$\,, \ though for some other \ $n$ \ equality \ $msr_n \ = \
p_{n-1}$ \ is restored.

Theorems 15 and 17 permit to reveal reasons of enough satisfactory
approximation of greatest series of the first axial
configuration to absolute values of corresponding maximum
series. As there is demonstrated above it is explained by
interval of series \ $SR_n^{(I)}(2)$\,, \ where each new grid
appears generated, differently it zerofills not less than two
units. If exactly, it zerofills as time two units, and also
does not bring one new multiple zero. Certainly, at increase of
series interval occurrence of multiple zeroes inevitably, but
they are formed with the help of the previous grids.

At the same time theorem 15 admits interpretation of separate
infringements of equality \ $msr_n \ = \ p_{n-1}$ \ for observed
values \ $n$\,. \ The greatest series value of the first
configuration with \ $q$ \ units expressed by following formula
$$
sr_n (q) \ = \ \max_{1 \le i \le [q/2]} \
\left\{p_{n+i} \ + \ \frac{p_{n+q-i} - p_{n+i}}{2}\right\}\,, \qquad
q \ \ge \ 4\,,
$$
can appear by initial decentralized series sum, that is maximum
value is reached in the formula for parameter \ $i$\,, \ not
equal to \ $ [q/2]$\,. \ For example, for \ $q = 4$ \ such
parameter \ $i = 1$\,. \ It means, that if in set compulsory
of non-generated grids for given formula realization there will
be such variant, at which one of grids will appear generated (it
zerofills two units), that there is filling \ $Z_{n+q-1}$\,, \
the series of which surpasses centralized series.

In particular, \ $sr_6 (2) = 19 = p_7 = sr_8' (0)$\,. \ But \
$sr_4 (5) = 20$\,, \ that is explained by advantage \ $0.5\,(p_8
- p_7) = 2$ \ over \ $0.5\,(p_7 - p_6) = 1$\,. \ As far as
generated grid is found for filling in variant with five units,
it has resulted to  value \ $msr_8 (0) = 20$\,. \ Case is quite
analogous: \ $sr_9 (2) = 31 = p_{10} = sr_{11}' (0)$\,. \ Here
is \ $sr_7 (5) = 33$\,, \ that is explained by advantage \
$0.5\,(p_{11} - p_{10}) = 3$ \ over \ $0.5\,(p_{10} - p_9) =
1$\,. \ As generated grid was found and here, \ $msr_{11} = 33$\,.

Naturally and hereafter for large parameter \ $n$ \ similar
effects can be observed. They little decrease maximum series \
$MSR_n \subset Kf_n^{(I)}$\,, \ that is value \ $msr_n$
\ relatively centralized variant \ $sr_n' (0) = p_{n-1}$ \
obviously following from relation for the first configuration \
$sr_{n-2}(2) = p_{n-1}$\,. \ Searches of such cases of variant
estimations would acquire greater sense, if the second
configuration \ $Kf_n^{(II)}$ \ has not interfered in generation
of maximum series \ $MSR_n$ \ and \ $MSR_n (q)$ \ for \ $q = O(n)$\,.

For first \ $Kf_n^{(I)} \subset SP_1$ \ configuration of \
$n$-filling (6.8) with two permanent central units each new grid
will be generated (it zerofills equally two units) in interval of
greatest series \ $SR_n^{(I)}(2)$\,. \ According to theorems 7
and 9 lengths of its greatest series \ $SR_n^{(I)}(q)$ \ for \
$q \ge 2$ \ units is equal \ $sr_n^{(I)}(q) = p_{n + q/2} +
\Delta$\,, \ where \ $\Delta = o (p_n)$\,. \ Then for
the greatest series of this configuration we receive
$$
SR_n^{(I)} (q) \subset Kf_n^{(I)}: \qquad \lim_{n,q \to \infty} \
\sup_q \,\frac{sr_n^{(I)} (q)}{p_{n + q}} \ = \ \lim_{n \to \infty}
\frac{sr_n^{(I)} (2)}{p_{n + 2}} \ = \ 1\,,
$$
and by that the ability to be generated for such kind of all
grids results to estimation of maximum series with help of the
first axial configuration. Received estimation is presented by
value \ $msr_n \sim p_n$\,. \ It is possible once again to note,
that such estimation was not surpassed for observed examples of
parameter \ $n$\,.

Besides from generation algorithm of the first configuration
follows, that zeroes frequency will be upper just in interval \
$SR_n^{(I)}(2)$\,. \ Its appreciable exceeding inevitably
results to occurrence of multiple zeroes for grids of greatest
modules. But it causes fall of essential factor of grids.
Purely, this phenomenon and quite precisely reflects the last
given relation.

However maximum observable series of the first configuration is
not occasion for extrapolational conclusions. Attentive consideration of
the second configuration \ $Kf_n^{(II)} \subset SP_1$ \ of expression
(6.10) for parameter \ $m = [\log_2\,p_{n + 1}]$ \ demonstrates
insufficiency of series of the first configuration \ $Kf_n^{(I)}$\,. \
Although among series of kind \ $p_{n+i+1} - p_{n+i}$ \ will meet and
units, corresponding to numbers \ $2^{m+s}\,, \ 2^sp_{n+j}$\,.

According to the theorem \ 16 \ and formula \ (6.10) the lengths \
$(sr_n)$ of greatest series \ \ $SR_n^{(II)}(2m)\,, \ SR_n^{(II)}(2m+1)$
\ and \ $SR_n^{(II)}(2m+2)$ \ are equal accordingly
$$
sr_n^{(II)} (2m) =
\max\{\,2^{m+1},\,2^{m-1} + p_{n+1}\,\}\,,
$$
$$
sr_n^{(II)} (2m+1) = 2^m + p_{n+1}\,, \quad sr_n^{(II)}
(2m+2) = \max\,\{2\cdot p_{n+1}, \ 2^m + p_{n+2}\}\,.
$$

The central core as proven here statements, as main numerical
characteristic of the second axial configuration \ $Kf_n^{(II)}
\subset SP_1$ \ is following idea:
\vspace*{3mm}

{\sl Series of the second axial configurations \ $SR_n^{(II)}(q)$ \ of
system \ $SP_1$ \ with units \ $q = 2m, \ 2m+1, \ 2m+2$ \ act by maximum
series \ $MSR_n (q)$ \ prime filling \ $Z_n \subset SP_1$ \ with
corresponding quantity of units for all \ $n \ge n_0$\,, \ since
some \ $n_0$\,.}
\vspace*{3mm}

This thesis requires careful consideration and confirmation
taking into account, that the initial values \ $n$ \ of fillings demonstrate
advantage of the first series of sequential filling (zero stage)
and series with two units for the first axial configuration \
$Kf_n^{(I)} \subset SP_1$\,, \ realizing by the first stage.

The second stage of transformation of maximum series coincides
with the coming superiority of series of the second axial
configuration \ $SR_n^{(II)}(q)$\,. \ It is realized for \ $n$ \
practically the depriving researchers of concrete check
opportunities, as resources of computer at searching algorithms
are rather limited. Nevertheless theorems 5 and 6 grant lower
estimations of maximum series, which already permit many.

Lengths of greatest series \ $SR_n^{(II)}(q) \subset SP_1$ \ mentioned by
last relations contain \ $2m \le q \le 2m+2$ \ units.
According to the rule of determination of parameter \ $m$ \ for
small \ $n$ \ these values \ $q$ \ can even it surpass.

We shall find those values \ $n$ \ for which in configuration \
$Kf_n^{(II)}$ \ there will be series of length \ $2^{m-1}$ \ first.
According to proven states, the core of configuration with \
$2 m$ \ units does not already change and will increase by
following series of length \ $2^m$\,. \ In such case for each parameter
\ $m$ \ its value \ $n = n_m$ \ will be found by scheme \
$p_{n+1} = \min\limits_i \{p_i > 2^m\}$\,. \ We shall give initial
values \ $n$\,:
\vspace*{1mm}

\begin {tabular}{||c|r r r r r r r r r r r||}
\hline \hline
m & 2 & 3 & 4 & 5 & 6 & 7 & 8 & 9 & 10 & 11 & 12 \\
\hline
$n = n_m$ & 1 & 3 & 5 & 10 & 17 & 30 & 53 & 96 & 171 & 308 & 559 \\
$p_n$ & 3 & 7 & 13 & 31 & 61 & 127 & 251 & 509 & 1021 & 2039 & 4093 \\
$p_{n+1}$ & 5 & 11 & 17 & 37 & 67 & 131 & 257 & 521 & 1031 & 2053 & 4099 \\
\hline \hline
\end{tabular}
\vspace*{5mm}

From the table it is visible, that only for \ $m \ge 6$ \ value
\ $n_m$ \ confidently surpasses considered parameters \ $q$\,. \
But also there are not enough it, as far as that to take
advantage of the second axial configuration for construction of
the greatest series, it is necessary part of grids (naturally
greatest modules) to send for liquidation all without exception
\ $2m$ \ central units of configuration \ $Kf_n^{(II)}$\,.

From representation of the second configuration (6.10) it is
possible to conclude, that for \ $n = n_m + s$\,, \ where \ $s =
0,\,1,\,...\,,O(n_m)$\,, \ relatively greatest series \
$SR_n^{(II)}(q)$ \ contains \ $q = 2m$ \ units. For \ $n =
n_{m+1} - s$\,, \ where \ $s = 1,\,2,\,...\,,O(n_m)$\,, \
opposite, relatively greatest series \ $SR_n^{(II)}(q)$ \
contains \ $q = 2m+2$ \ units. In intermediate variants there
can quite appear most acceptable \ $q = 2m+1$\,.

The fixed places (points), in which units of considered central core of
configuration \ $SR_n^{(II)}(q) \subset SP_1$ \ stand, permit to
estimate opportunity of filling,
that is zerofilling of these units by grids of greatest modules,
that series without units \ $sr_n (0)$ \ to generate by corresponding
redistribution.

It is easy to find odd distances between units of central core
of configuration \ $SR_n^{(II)}(q) \subset SP_1$ \ presented only
by numbers of kind \ $2^k - 1$ \
or \ $2^k + 1$\,, \ where \ $k \le m$\,. \ It means, that in
considered interval (especially for large value \ $n$\,) primes as
modules of grids can not find more than one, satisfying these
conditions. Really, interval is limited by degree of two and
for one \ $k$ \ values \ $2^k - 1$ \ and \ $2^k + 1$ \ can not
be primes simultaneously.

Thus practically all \ $2m$ \ grids \ $S(p_i) \subset Z_n$ \ of large
modules \ $p_i$ \ should be
directed to generation of filling \ $Z_n$ \ with series of the
second axial configuration \ $SR_{n-2m}^{(II)}(2m)$ \ of length
\ $sr_{n-2m}(2m)$ \ already not containing units. According to
given relations this value is equal
$$
sr_{n-2m}^{(II)} (2m') \ = \ \max\{\,2^{m'+1},\,2^{m'-1} +
p_{n-2m+1}\,\}\,, \quad m' = [\log_2p_{n-2m+1}]\,.
$$

For considered parameters \ $n$ \ value \ $m'$ \ is possible
only per unit less \ $m$\,, \ therefore at estimation of
greatest series is allowable to consider \ $n$ \ such, for which
\ $m = m'$\,. \ In view of all these conditions, expressions
and conclusions in filling \ $Z_n$ \ there will be series
without units of length \ $2\,p_{n-2m+1} + \Delta$\,, \ where \
$\Delta = O(p_{ n-2m+1})$\,.

However that this series of the second configuration \
$sr_{n-2m}^{(II)} (2m') \subset Kf_n^{(II)}$ \ exceeds serial
structure of the first configuration \ $Kf_n^{(I)} \subset SP_1$\,,
\ is obviously
necessary execution of inequality \ $2\,p_{n-2m+1} > p_{n-1}$\,. \
According to theorem 6 such inequality will be executed at \ $n \ge 26$\,.
\ Advantage of the second configuration series for \ $n \to
\infty$ \ increases, approached to coefficient two.

Thus maximum series \ $SR_{n-q}^{(II)}(q)$ \ with quantity of
units \ $2m \le q \le 2m+2$ \ automatically means too maximum
series without units, constructed by described way in interval of
central core of configuration for \ $n \ge 26$\,.

For proof of maximum series \ $SR_{n}^{(II)}(q)$ \ for \ $2m \le
q \le 2m+2$ \ it is necessary again to address to concept of
generated grids, in this case in interval of central core of
configuration \ $Kf_n^{(II)} \subset SP_1$\,. \ Each new grid \
$S(p_n)$ \ in interval of configuration \ $Kf_n^{(II)}$ \ of length
\ $2\,p_{n+1}$ \ is responsible for one zero of multiplicity two
(additional multiplicity) and also for zerofilling of two units.
The comparison as though for the benefit of series of central core of
the first configuration \ $Kf_n^{(I)} \subset SP_1$\,, \
where at same two eliminated units does not occur multiple zero.

This conclusion has hurried character. According to it, in
general becomes inexplicable occurrence of advantage of series
of the second configuration \ $Kf_n^{(II)} \subset SP_1$\,. \ The
reason that in extended central
interval of length \ $4\,p_{n+1}$ \ the same grid \ $S(p_n)$ \
zerofills too units, corresponding to values \ $2\,p_n$\,.
\ Thus in such extended interval property of grid \ $S(p_n)$ \
to be generated grid reflected by four eliminated (zerofilled) units
and one multiple zero.

From elementary product of grid \ $S(p_n)$ \ with grid \ $S (p_1
= 3)$ \ follows, that or such filling \ $Z_n \subset SP_1$\,, \
limiting on efficiency
and opportunity of grids to be generated grid in extended
interval, exceeds it is impossible, or for even more extended
interval it will take place, but in such case it will be
executed by the same structure of the second axial configuration
\ $Kf_n^{(II)} \subset SP_1$\,.

There can arise question, how the eliminated second pair of
units of extended interval of length \ $4\,p_{n + 1}$ \
influences to occurrence of series of length \ $2\,p_{n + 1}$\,.
\ Answer is extremely simple. In series of length \ $2\,p_{n + 1}$ \
eliminated second pairs of units of grids group of smaller
modules have come. It predetermines creation of the second
configuration with central series of such extent.

It is necessary to take into account, that the minimum grid \ $S(3)$ \
dictates and determines impossibility of excess of essential
coefficient for product with grid \ $S(p_n)$ \ in interval of
length \ $4\,p_{n+1}$ \ by formula (4 eliminated units -- 1
multiple zero). Exactly such quantity of units is
zerofilled by the second axial configuration. From here becomes
explicable impossibility of excess twice of length of maximum
series \ $MSR_n(0)$ \ of greatest grid module: \ $msr_n \not>
p_n$\,.

And really, we shall consider occurrence in half-configuration \
$Kf_n^{(II)}$ \ of zeroes of multiplicity one and higher,
introduced in filling \ $Z_n$ \ by grid \ $S(p_n)$\,. \
Sign \ $\emptyset$ \ is zero of increased multiplicity:
\vspace*{5mm}\\
\begin{tabular}{c@{~~}c@{~~}c@{~~}c@{~~}c@{~~}c@{~~}c@{~~}c@{~~}c@{~~}c@{~~}c@{~~}c@{~~}c@{~~}c}
0 & $p_n$ & $2p_n$ & $3p_n$ & $4p_n$ & $5p_n$ & $6p_n$ & $7p_n$
& $8p_n$ & $9p_n$ & $10p_n$ & $11p_n$ & $12p_n$ & ... \\
$\emptyset$ & o & o & $\emptyset$ & o & $\emptyset$ & $\emptyset$ &
$\emptyset$ & o & $\emptyset$ & $\emptyset$ & $\emptyset$ &
$\emptyset$ & ... \\
.. & ... & $\frac{4}{1}$ & $\frac{4}{3}$ & $\frac{6}{3}$ &
$\frac{6}{5}$ & $\frac{6}{7}$ & $\frac{6}{9}$
& $\frac{8}{9}$ & $\frac{8}{11}$ & $\frac{8}{13}$ &
$\frac{8}{15}$ & $\frac{8}{17}$ & ... \\
\end{tabular}
\vspace*{5mm}

Here symbol \ "o" \ means zero of one multiplicity of grid \
$S(p_n)$ \ for filling, which is responsible for increase of
zero series, as far as zeroes of increased multiplicity can not
already play such role. Naturally, zero of multiplicity one of
grid \ $S(p_n)$ \ stands in place \ $kp_n$ \ only in that case,
when it takes place \ $(k,\,\prod_1^{n-1}\,p_i) = 1$\,.

Correlation of quantity of generated zero of multiplicity one to zero
increased multiplicity in the same interval appreciably
decreases with growth of interval \ $I_n$\,. \ If to take into account
half-configuration of given chain, since interval of length \
$4\,p_{n+1}$ \ we shall receive the correlation \ 4/1\,, \ which as
will remain greatest in third line of the table. Even appeal to
initial filling \ $n = 2$ \ with grid \ $S(p_2) = S(5)$\,, \
when only zeroes of kind \ $15 k$ \ will be multiple, does not
correct situation, as the subsequent correlation is \ $\frac{8}{3} < 4$\,.

Thus correlation \ $4/1$ \ of quantity of unitary zero to multiple
quantity for interval of length \ $I_n = 4\,p_{n+1}$ \ is impossible
to exceed, and it is impossible and in extended interval,
including structure of the second configuration \
$Kf_n^{(II)} \subset SP_1$\,. \ But if second axial configuration \
$Kf_n^{(II)}$ \ realizes exactly such limiting generation,
which restricts length of series by value \ $sr_{n}^{(II)}(q)$ \
for \ $2m \le q \le 2m+2$\,, \ it and will be maximum in period
\ $PZ_n$ \ of filling \ $Z_n \subset SP_1$\,.

Assuming that common problem is fixed to create product of grids
in interval under condition of greatest zerofilling of units
set. In such case it is impossible to escape product of next
grid \ $S(p_n)$ \ with first and basic grid \ $S(3)$\,. \ First
stage of creation of such product without multiple zeroes is
characterized by reception of the first axial configuration \
$Kf_n^{(I)}$\,, \ in which each new grid appears generated
that is it zerofills two units in interval of greatest series
\ $SR_n^{(I)}(2)$\,.

But we shall notice, the length of such series \ $sr_n^{(I)}(2)
= p_{n+1}$ \ appears not limiting for large \ $n$\,. \ Therefore
has to consider other product, other configuration, already with
multiple zero in interval of series.
The product of grids \ $S(3)\,, \ S(p_n)$ \ in interval of
length \ $4\,p_{n+1}$\,, \ when central zero of multiplicity two
appears common, permits in given variant to zerofill at once
four units. Easily to see, such construction will be executed
for each new grid.

Above is proven that product of grids \ $S(3)\,, \ S(p_n)$ \ in
any interval of greater length not capable to ensure large
concentration of zeroes of multiplicity one. It is necessary to
note, for all that effect of intermediate grids is not taken
into account. It does not touch on all four new essential zero,
but results in additional multiple zeroes in interval, exceeding
of length \ $4\,p_{n+1}$\,.

Thus described structure of sequential product of grids provides
maximum concentration of unitary zeroes in interval with length
\ $4\,p_{n+1}$\,. \ But this structure and is the second axial
configuration \ $Kf_n^{(II)}$\,. \ If to take into account, that
its greatest series with central units \ $SR_n^{(II)}(2m+2)$ \
has length about \ $2\,p_{n+1}$\,, \ that is far from \
$4\,p_{n+1}$\,, \ it and means that the series with \ $(2m+2)$ \
units of considered second axial configuration are maximum in
all period of filling.

The outstripped interval of extreme property \ $4\,p_{n+1}$ \
relatively maximum series, not exceeding \ $2\,p_{n+1}$\,, \
consists that zeroes of grid \ $S(p_n)$ \ form the basis for the
subsequent maximum series, though they do not include in series,
remained unitary. Thus the realization of maximum series
permanently overtakes interval of extreme properties of grids
to be generated. It creates preconditions of successful
preservation of established relation of intervals ($\approx 2$).

The second axial configuration realizes ineradicable condition
of preservation of property of grids to be generated in interval
of length \ $2\,p_{n+1}$ \ and increased property -- in interval
of length \ $4\,p_{n+1}$\,. \ As far as at increase \ $n$ \ this
extremity will also find reflection in generated the greatest
series \ $SR_{n+s}^{(II)}(q)$\,, \ that such property can not be
already abolished. Interval of length \ $4\,p_{n+1}$ \ is least
with such high essential coefficient (remaining all appreciably
more length of greatest series) and the received series of the
second axial configuration \ $SR_{n}^{(II)}(q)$ \ for \ $2m \le
q \le 2m+2$ \ is maximum.

If to take into account, that the third axial configuration is
not present, the central thesis about maximum series \
$SR_{n}^{(II)}(q)$ \ is proven. For the proof of theorem and
expression (6.19) now it is enough to remind, that grids of the
greatest modules, sent for elimination of units of central core
of configuration, at the best once (that is one grid) can
zerofill two units (in the other cases -- only one).

Thus, reception of maximum series \ $MSR_n(q)$ \ for \ $q < 2m$
\ or even for \ $q = 0$ \ for large \ $n \ge 26$ \ is directly
connected with series \ $SR_{n}^{(II)}(q)$ \ for \ $2m \le q \le
2m+2$ \ and maximum objects in period of filling are generated
with the help of typical series of the second configuration of
previous \ $n' < n$ \ by the same way, which they were submitted in
model. It is redistribution of zeroes of grids in filling \
$Z_n$\,.

Let series \ $SR_n^{(II)} (q)$ \ is maximum \ $MSR_n (q)$ \ series and
$$
MSR_n (q)\,: \qquad sr_n^{(II)} (q) \ = \ msr_n (q)\,, \qquad q = 2m
\quad and \quad q = 2m+2\,,
$$
that is maximum series \ $msr_{n+2m} (0)$ \ is created in its
basis. Then at transition from \ $n$ \ to \ $n+2$ \ series \ \
$SR_{n+2}^{(II)} (2m)$ \ arises without fail for creation
maximum series \ $msr_{n+2m+2} (0)$\,. \ This is connected with
grids \ $S(p_{n+1})$ \ and \ $S(p_{n+2})$\,. \ They are used
ineffectively at generation of maximum series \ $msr_{n+2m}
(0)$\,. \ Moreover let is formed maximum series \ $MSR_{n+2m} (0)$ \
on basis of the greatest series \ $SR^{(II)}_n (2m)$ \ with help \
$2m$ \ grids \ $S(p_{n+i})$ \ for \ $1 \leq i \leq 2m$\,.

Then series \ $SR^{(II)}_{n+2m} (q)$\,, \ where \ $q = 2m$ \ or
\ $q = 2m+2$ \ is source of generation of next maximum series \
$MSR_{n+2m+q} (0)$\,. \ It follows from inequality
$$
sr_{n+2m+q} (0) \ < \ sr^{(II)}_{n+2m} (q)\,, \qquad SR_{n+2m+q}
(0) \ \supset \ MSR_{n+2m} (0)\,,
$$
so far as grids \ $S(p_{n+j})\,, \  1 \leq j \leq 2m+q$ \ can
zerofill only \ $q + 2m + \epsilon$\,, \ where \ $\epsilon = O(m)$ \
units. Decentralized series on basis of the maximum series \ $MSR_{n+2m}
(0)$ \ can direct to addition of a little more one zerofilling
of units by each grid. For all that series \ $SR_{n+2m+q} (0)$ \
increases by only average series of filling in period.

Axial configurations \ $Kf_{n}^{(I)}$ \ and especially \
$Kf_{n}^{(II)}$ \ are natural, unique and active factories of
maximum series. Natural packing of grids in the second axial
configuration maximum effectively and clearly sequentially
realizes activity of every new grid at generation of the maximum
series \ $MSR_n (2m)$\,. \ Any infringement of packing leads to
loss of maximum. \hfill $\Box$
\vspace*{3mm}

The proven theorem permits quite essentially and considerably more
precisely to judge distribution of maximum series with units and
without them in given filling.
\vspace*{3mm}

{\bf Theorem 19.} Maximum series \ $MSR_n$ \ and \ $MSR_n(q)$ \
for all \ $n,\,q$ \ in prime system \ $SP_1$ \ have the upper estimation
$$
MSR_n (q): \qquad msr_n (q) \ < \ C_n\,p_{n+q}\,; \quad C_n <
1\,: \ \ n + q < 25\,;
$$
$$
C_n \ > \ 1\,: \qquad n \ > \ 15 \quad and \quad n + q \ > \ 28\,;
\quad 2 \ < \ C_n \ \stackrel{n \to \infty}{\Longrightarrow} \
2\,. \eqno{(6.20)}
$$
\vspace*{0mm}

{\sl Proof.} For initial \ $n \le 3$ \ value of maximum series \
$MSR_n \equiv MSR_n (0)$ \ coincides with the first
series \ $SR_n^1$ \ of sequential filling:
$$
MSR_n \subset SP_1: \qquad sr^1 (0) \ = \ p_{n-1} \ = \
\frac{1}{2}\,(p_{n+1} - 1)\,,
\quad 1 \le n \le 3\,.
$$

But already for \ $n = 4$ \ this equality is infringed and
advantage goes to series of the first axial configuration \
$Kf_n^{(I)}$\,. \ According to theorems 17 and 18 as well as with
given calculations, for \ $n < 25$ \ equality is valid
$$
MSR_n (2) \subset SP_1: \qquad msr_n (2) \ = \ sr_n^{(I)} (2)\,,
\qquad 3 \le n < 25\,,
$$
whence follows \ $msr_n (0) \ < \ p_n$\,. \ If to take into
account, that parameter of units of series \ $q > 0$ \ does not
change estimations, only sharply them easing at increase, the
first part of the theorem is proven.

However the further increase of main parameter of filling \ $n$
\ results in loss of advantage of series of the first
configuration \ $Kf_n^{(I)}$\,, \ which goes to second \
$Kf_n^{(II)}$\,. \ From the same theorems 9 and 10 it is possible
to receive equality, especially not pay attention to details to
extreme exact expression, that is having chosen one of variants
of maximum series \ $MSR_n (2m+2)$\,:
$$
msr_n (2m + 2) \ = \ sr_n^{(II)}(2m+2) \ = \ 2\cdot p_{n+1}\,;
\quad n > 15\,, \ m = [\log_2 p_{n+1}]\,.
$$

For \ $n$ \ close \ 12 -- 15 \ the values of corresponding
maximum series of both configurations are reasonably close one
another, that is connected with proximity of values \ $n$ \ and
\ $2 m$\,. \ Thus and coefficient \ $C_n$ \ of theorem will be
close to 1. But from the last equality clearly, that at the
increase \ $n$ \ value \ $2m$ \ as \ $2m = O(n)$ \ all less
influences for estimation of maximum series \ $MSR_n \subset SP_1$\,:
$$
MSR_n (0) : \qquad msr_n \ < \ 2\cdot p_{n-2m}\,;
\quad n + 2m > 28\,, \ m = [\log_2 p_{n+2m'}]\,,
$$
as well as \ $msr_n \ > \ p_n$\,, \ and where value \ $m'$ \ not
more, than per unit differs from \ $m$\,. \ To the point from
connection of systems \ $SP_1$ \ and \ $SP_0$ \ and received
inequality immediately follows
$$
MSR_n (0)\ \subset \ SP_0: \quad \qquad msr_n \ < \ 4\cdot p_{n-2m-1}
$$
under the same easily attainable conditions. Taking into account, that
parameter \ $q > 0$ \ inequality do not infringed, from ultimate
expression the proof of the second part of the theorem and
relation (6.20) follows. \hfill $\Box$
\vspace*{3mm}

{\bf Theorem 20.} The validity of theorem 19 for system \ $SP_1$ \
means validity of the main theorem for any non-singular system \
$SS$\,, \ and with the same upper coefficient two, as well as in
unimproved expression (4.7).
\vspace*{3mm}

{\sl Proof.} As far as in system \ $SP_1$ \ there is expression
$$
\frac{n}{1 - \gamma_n} \ = \ n\,\prod_{i=1}^n \left(1 +
\frac{1}{p_i - 1}\right)\ = \ C_n\,p_n\,, \qquad C_n \ < \ 1\,,
\eqno{(6.21)}
$$
from theorems 13--18 and especially from theorem 19 follow
validity of the main theorem in the form (4.7) for this system.
According to lemma 2 if the main theorem is valid in system \
$SP_1'$\,, \ it is fair and in any other system, and its
estimations will remain by majorizing. As far as under theorem 4
systems \ $SP_1'$ \ and \ $SP_1$ \ differ whole by one grid, the
characteristics of their maximum series are rather close. Taking
into account, that \ $C_n$ \ in (6.21) does not reach unit, the
upper bound 2 in estimation (4.7) will remain and for system \
$SP_1'$\,. \ Besides possibility of exception for system \
$SP_1'$ \ is refuted by inequality \ $msr_{n+s}[SP_1] > msr_n
[SP_1']$\,, \ where \ $s = O(n)$\,.

It is necessary to note, that in monograph [1] proofs of the main theorem
are placed also, giving for system \ $SP_1$ \ estimations more
weak, than value of expression (6.21): \ $C^{(1)}_n\,p_n \ln \ln n\,;
\ \ C^{(2)}_n\,p_n \ln \ln n \ln n\,; \ \ C^{(3)}_n\,p_n\,, \ \
C^{(3)}_n \ < \ 3$\,. \hfill $\Box$
\vspace*{5mm}

\begin{center}{\Large \bf 7. Two-sided fillings and main theorem}
\end{center}
\vspace*{3mm}

The statements of 5 chapter grant some other approach to the proof of
the finishing form of the main theorem produced by theorem 20. It is
thus necessary to note a decisive role of complex (summarized)
characteristics \ $n$-fillings in estimation of its major numerical
parameters.
\vspace*{3mm}

{\bf Theorem 21.} Upper estimations of maximum series \ $MSR_n(q)$ \
depend on density of zero of filling \ $\gamma_n$ \ and parameters
\ $(n,\,q$ \ for all classes of systems.
$$
MSR_n (q,\,SS): \qquad msr_n (q) \ \le \ FF\{n,\,q,\,\gamma_n,\,SS,\,Z_n\}\,,
\quad \forall\,\{Z_n \subset SS\}\,. \eqno{(7.1)}
$$
\vspace*{1mm}

{\sl Proof.} Certainly, only non-singular fillings and systems are meant.
The statement and formula estimations would not have equally the price,
if not unity of law \ $FF$ \ outside of dependence from a class of system
\ $SS$\,. \
Reasonably the weak difference of concrete functions \ $FF_1(SS_1)\,, \
FF_2(SS_2)$ \ for the most various systems is meant too. The discrepancy of
functions \ $FF_1 \neq FF_2$ \ is in complete dependence from
availability whether or not representative sets of multiple zeroes in a
appropriate system. The systems from first to third types are various.

Such function \ $FF$ \ (7.1) does not depend on a particular set of grids, if
given \ $n$-filling does not leave the same class of systems. And if other
set of grids has same or though close density of zeroes \ $\gamma_n$\,, \
estimations of maximum series will coincide or are close. Such unity
permits to consider a number of problems without concreteness of this or
that system, as far as for them it should expect uniform conclusions.

The function of the main theorem \ $MT$ \ from parameters \ $(n,q)$ \ and
density (average frequency) of zeroes for the period of filling in the
formula of the main theorem acts unified majorant of values of maximum
series \ $MSR_n(q)$\,. \ And irrespective of availability of multiple
zeroes as in period, as in series. But if multiple zeroes have managed
to avoid or any image pass to variant without multiple zeroes in period
(even only in interval), the global problem estimate of maximum series
can be considered permitted. \hfill$\Box$
\vspace*{3mm}

Will not hinder to specify occurrence of the universal formula (5.9,\,7.1)
of main theorem \ $MT$\,, \ playing by a determining role in the
fillings method.
\vspace*{3mm}

{\bf Theorem 22.} The coefficient \ $\tau_n$ \ of the main theorem \ $MT$ \
does not exceed unit \ $\tau_n \le \tau = 1$ \ for degree-systems, without
multiple zeroes and with primary growth of modules. It determines
one-sideness of their fillings \ $Z_n^{(1)}$\,.
$$
\left\{SS_d,\,SS',\,SS_{[2]}\right\}: \quad msr_n (q) \ \le \
\tau\,\frac{n + q}{1 - \gamma_n} + 1\,; \quad \gamma_n \ = \
\frac{H_n}{PZ_n}\,, \quad (\tau = 1) \ \Rightarrow \ Z_n^{(1)}\,.
\eqno{(7.2)}
$$
\vspace*{1mm}

{\sl Proof.} First occurrence of the formula of the main theorem on
the basis of exact expression of value of a maximum series \ $MSR_n(q)$ \
for a class of degree-systems \ $SS_d$ \ can be complemented a little
by other reasons. We shall consistently consider some turning-points for
the clearing of exposition.
\vspace*{1mm}

1. We shall evaluate important value, if frequency (density) of zeroes \
$\gamma_n$ \ in period \ $PZ_n$\,, \ period and other characteristics
of \ $n$-fillings are
$$
PZ_n = H_n + E_n\,, \ \ \gamma_n = \frac{H_n}{PZ_n}\,: \quad
r_n \ = \ \frac{1}{1 - \gamma_n}\,; \quad r_n \ = \ \frac{PZ_n}{E_n} \ =
\ \frac{H_n}{E_n} \ + \ 1\,. \eqno{(7.3)}
$$

Value \ $r_n$ \ from expressions (7.3) it is possible to interpret by
quantity of zeroes per unit of a period, summarized with unit.
Thus \ $r_n$ \ according to definition is interval, stipulated by
zeroes of average grid.

2. We shall consider examples. The value \ $r_n$ \ equally to unit at \
$n$ \ infinite grids, and value of maximum series coincides with the
upper estimation in the main theorem \ $msr_n \ = \ \frac{n}{1 - \gamma_n}
+ 1 \ = \ n+1$\,. \ Thus estimation of a series with coefficient \
$\tau = 1$ \ is unimproved. At the same time the value of maximum series
\ $msr_n$ \ will remain same, that is \ $n+1$ \ and at \ $a_1 \ge n+2$\,.
\ Though in such case density of zeroes, obviously, different from zero
\ $(\gamma_n' \neq 0)$\,.

The interval of a series with zeroes of one average grid is equal two at \
$\gamma_n \ = \ \frac{1}{2}$\,, \ that is at \ $r_n \ = \ \frac{1}{
1 - \gamma_n} \ = \ 2$\,. \ It means, that the estimation of maximum series
at such density \ $\gamma_n$ \ is equal \ $2n + 1$\,, \ and
\ $msr_1(q-1) \ = \ 2q$\,, \ that is estimation more length of series
whole per unit (the affinity of values is doubtless) at identical
quantity forming grids or grid with units.

Interval of a series with zeroes of one average grid fractional at density
of zeroes \ $\gamma_n \ = \ \frac{1}{3}$\,, \ that is at \ $r_n \ = \
\frac {1}{1 - \gamma_n} \ = \ \frac{3}{2}$\,. \ It does not interfere
formation of the upper estimation of maximum series \ $\frac{3n}{2} + 1$
\ or \ $\frac{3(1 + q)}{2} + 1$\,, \ and it more values \ $msr_1(q)$ \ for
\ $1$-filling -- grid \ $S(3)$\,.

3. The considered examples demonstrate and confirm a role of value
\ $r_n$\,, \ zero reflecting zeroes of an averaging grid of \ $n$-filling.
Then the value \ $r_n + 1$ \ will appear by the upper estimation of
maximum series \ $MSR_1$\,, \ found with the help of such grid. The
estimation is well grounded as far as from expression (7.3) follows:
the quantity of zeroes per unit of filling takes into account the
frequent contribution of all grids of \ $n$-filling. From here the
numerical value \ $r_n + 1$ \ is absolute majorant of series, formed
by one grid. This estimation is achievable \ ($r_n + 1 = 2$) \ for
infinite grid. In other cases the estimation is not achievable.

4. Attraction of the second unit or value \ $n = 2$ \ results in summarized
estimation \ $2r_n + 1$ \ of the heaviest series by two grids. This
estimation of maximum series with the help \ $r_n$ \ zeroes precisely such
a conditional grid takes into account summarized density of all grids of
filling. The estimation of maximum series \ $MSR_n$ \ by all \ $n$ \ grids
arises completely similarly. Taking into account, that the frequent
contribution of each initial grid \ $S(a_i)$ \ is included in values \
$r_n$\,, \ received estimation \ $n\,r_n + 1$ \ is obliged to surpass
actual value \ $msr_n$ \ of the most maximum series \ $MSR_n$\,. \
It occurs in a reality.

5. The address to variant \ $q > 0$\,, \ certainly, does not reduce
efficiency of estimations of maximum series of the main theorem.
Any growth \ $q$ \ provides linear increase of estimation, which can become
only superfluous. The classes of systems, for which indicated reasons
(for coefficient \ $\tau = 1$ \ of the main theorem), are described by the
theorem 11.

6. If to consider conditional grids from \ $r_n$ \ intervals was held, there
is the question about parameters provided that all \ $n$ \ such grids
will form conditional filling with the same density of zeroes.
There is expression immediately following from sum of zeroes of all
\ $n$ \ grids, as far as these grids are identical:
$$
\left\{\frac{n\,r_n}{PZ_n'} \ = \ \gamma_n\right\} \ \Rightarrow \
\frac{n}{(1 - \gamma_n)\,PZ_n'} \ = \ \gamma_n\,; \quad
PZ_n' \ = \ \frac{n}{\gamma_n\,(1 - \gamma_n)}\,; \quad a_i' = PZ_n'\,,
$$
that is modules of such grids \ $a_i' = PZ_n' > n\,r_n$ \ coincide with
the conditional period.

7. Creation of conditional \ $n$-filling with conditional grids of \
$r_n$ \ zero for a conditional period not only results to majorizing
upper estimations of maximum series, but it also explains sources of
origins of the main theorem.
\vspace*{1mm}

Alongside with the theorem 11 statement and expression (7.2) prove,
that by initial premises of a coefficient \ $\tau = 1$ \ in the
formulation of the main theorem and one-sidedness of filling \
$Z_n^{(1)}$ \ a system correlation between modules acts. It is realized in
classes of systems of grids \ $SS_d,\,SS',\,SS_{[2]}$\,. \ The
interdependence of the described characteristics is reasonably obvious.
\hfill$\Box$
\vspace*{3mm}

Rather in detail the investigated systems of grids of the previous theorems
can not affect central interest, which cause systems of the third type with
other components -- grids of modules of completely other kind.
\vspace*{3mm}

We shall remind about achieved. Estimate possibility of the upper values of
maximum series common (summarized) density of zeroes is the main purpose
of the fillings method. The theory of numbers in sieving process tries
to operate with modules of separate grids \ $S(a_i)$\,. \ Besides the
fillings method is oriented to a period of products of grids, and
other methods consider at the best interval of length as a square of the
heaviest module of filling. Would seem, unremovability of multiple
zeroes in somehow appreciable interval is a insuperable obstacle to
further reasons and conclusions. Especially, for major systems (for
example, \ $SP_1$), \ the summarize set of all multiplicities in many
times surpasses a initial period.
\vspace*{3mm}

On a way of the decision of the problem of essential set of multiple
zeroes, so characteristic for many important systems and fillings, it is
necessary to specify the rule of transformation once more.
\vspace*{3mm}

{\bf Definition 30.} If \ $PZ_n$ \ is period of filling \
$Z_n$\,, \ $H_n$ \ is quantity of zeroes and \ $E_n = PZ_n - H_n$ \ is
units for this period at density of zeroes \ $\gamma_n$\,, \ transition
to \ $n$-fillings with majorizing frequency of zeroes \ $\gamma_n^*$ \
of multiplicity unit for that a set of grids \ $\{S(a_i)\}$ \ is executed
under the scheme
$$
PZ_n = H_n + E_n\,, \ \ \gamma_n = \frac{H_n}{PZ_n}\,: \quad
H^*_n \ = \ PZ_n \sum_{i=1}^n\,\frac{1}{a_i}\,, \ \ \gamma_n^* \
= \ \frac{H_n^*}{H_n^* + E_n}\,. \eqno{(7.4)}
$$
\vspace*{1mm}

Completely obviously, that \ $\gamma_n \ \le \ \gamma_n^*$\,, \ and the
equality can be observed only and only in case, when initial \ $n$-filling
has not multiple zeroes.
\vspace*{3mm}

{\bf Theorem 23.} The density of zeroes \ $\gamma_n^*$ \ in the
formulation of the main theorem provides absolute majorizing estimation \
$msr_n (q)$ \ of maximum series with a coefficient \ $\tau = 1$ \ in the
class of any non-singular systems.
$$
MSR_n (q) \ \subset \ \forall\,SS\,: \quad \left\{\gamma \ = \ \gamma_n^* \
= \ \frac{H_n^*}{H_n^* + E_n}\right\} \ \Rightarrow \
\left\{msr_n (q) \ \le \ \frac{n + q}{1 - \gamma_n^*} + 1\right\}.
\eqno{(7.5)}
$$
\vspace*{1mm}

{\sl Proof.} The statement reflects idea of a opportunity of local
redistribution of zeroes high multiplicity in lowered multiplicity,
down to unitary. The expansion of a allocated interval thus occurs,
and the quantity of units remains constant. The decrease of number of
units in a interval would mean not redistribution of zeroes, but
product with unknown new grid. It contradicts the principle of filling.
Really the limited opportunities of similar redistribution (for the interval
whole period they are reduce to zero), can not be obstacle to
idea of reception of imaginary filling without multiple zeroes in period.

The disposal from zeroes of multiplicity higher unit is, purely, decision
of the problem of upper estimation of maximum series. From here a reason
about replacement of actual density of zero for a period, where zero
of multiplicity \ $k > 1$ \ is one zero. The former frequency \ $\gamma_n$
\ varies other \ $\gamma_n^*$\,, \ and at its formation multiple zeroes
are redistributed by a set of zeroes of multiplicity unit. Received in a
result obviously the higher conditional density appears by majorizing
frequency of zero in the formulas of estimation of a maximum series.
Thus the theorem 3 about the upper estimation of a series \ $MSR_n(q)$ \
is applicable to such conditional filling (already without multiple zeroes).

However the density (frequency of zeroes) found thus appears obviously
redundant for practically important systems in applications (especially
for systems of third type). It is
explained unremovability of multiple zeroes in intervals of maximum series,
and even at reasonably small \ $n$\,. \ For example, for a system of
primes \ $SP_1$ \ achievement rather small value \ $msr_n \sim 15$
\ during growth \ $n$ \ not in forces to avoid multiple zero in interval of
such length. Product of modules of the first two grids equally 15, and zero
of multiplicity two is unremoved.

From here supervision follows. The increase of main parameter \ $n$ \ of
filling and consequently value \ $msr_n(q)$ \ (even at \ $q = 0$) \ in
systems of the third type provides a fast increase of a set of those
multiple zeroes, which can not participate in redistribution. Though it
is conditional mental operation.
\vspace*{3mm}

Let any non-singular system \ $SS$ \ and fillings \ $Z_n$ \ are given.
The construction of majorizing density of zeroes at preservation of other
parameters of filling \ $(n,q)$ \ consists of expansion number of zeroes
by them multiplicities for increasing period. The quantity of units
remains constant at such transformation. Then the received density
reflects density of filling without multiple zeroes
or even degree-filling.
\vspace*{3mm}

The deep sense of the generalized formula of the main theorem consists of
confirmation of a decisive role of density of zeroes \ $\gamma_n$ \ for
the upper estimation of maximum series \ $MSR_n(q)$\,. \ Such estimation
does not depend already from a class of a system and particular grids
\ $\{S(a_i)\}$\,, \ which enter in filling. The idea of the formula of
the upper estimation is incorporated in reasonably transparent reasons.

The variants considered above of systems \ $SS_d\,, \ SS'$ \ and their
fillings have not zeroes of multiplicity higher unit. Then from definition
28 and expressions (7.4,\,7.5) coincidence \ $\gamma_n = \gamma_n^*$ \
and already proven case of the main theorem follows.

If initial grids \ $\{S(a_i)\}$ \ are infinite (infinite module
\ $a_i = \infty$)\,, \ the upper estimation of series coincides their
length at \ $\tau = 1$\,. \ The variant of density of zeroes, different
from a zero \ $\gamma_n > 0$\,, \ is also considered by the theorem 22.
The theorem talks, that on each unit of filling drops \ $r_n \ = \
\frac{1}{1-\gamma_n} \ > \ 1$ \ -- interval with zeroes. So conditional
filling of \ $n$ \ equal grids creates, and value \ $r_n$ \ is the main
characteristic of a compact arrangement of zeroes, including all frequent
components of former \ $n$ \ grids. From here the simple summation of
values \ $r_n$ \ provides majorizing character of numerical estimations
of maximum series.

Or else, the created conditional grid in difference from standard contains
not one zero in a own period, equal to module, but \ $r_n$ \ at unified
averaged period. Then the heaviest series is formed by simple association
of zeroes each from \ $n$ \ equal grids. Such process of formation conducts
to creation of majorizing estimation of maximum series, and at all not to
creation of series of zeroes \ $MSR_n$ \ of initial filling. Certainly,
in common case the number \ $r_n$ \ is not integer, but also
it not a obstacle to conditional construction of a redundant modernized
series, consequently and upper estimation.

Now simplely to notice, that the described mechanism of the account of
multiple zeroes and density of zeroes of multiplicity not above unit
is situated also in a scheme of redistribution of zeroes at formation
of the upper estimation of maximum series \ $MSR_n(q)$\,. \ The count
of quantity of zeroes in the period of \ $n$-filling can be considered
by arithmetic consecutive summation of zeroes of all grids. It reduces
process to the scheme (7.4), to a sum of multiplicities of zeroes, and to
density \ $\gamma_n^*$\,. \ And the discrepancy with \ $\gamma_n$ \ is
not the certificate of falsehood.

The received estimation is absolutely upper for a class of any non-singular
systems and fillings. Certainly, such conclusion does not remove the
subsequent conclusions concerning redundancy of the last found estimation
(7.5) for systems with abundance of multiplicities. This consequence of
growth of a set of unremovable multiple zeroes and parallel increase of a
set of grids of large modules, participating in creation of density
\ $\gamma_n^*$\,. \ But they introduce in comparison weak contribution
to formation of a concrete series \ $MSR_n$\,. \hfill$\Box$
\vspace*{3mm}

The obvious incompleteness of effect of multiple zeroes, even promptly
growing, compels critical revision of idea of frequency \ $\gamma_n^*$ \
as absolute majorant of the theorem 23 and expression (7.5). Too frank
redundancy of density \ $\gamma_n^*$ \ forces to resort to restrictions.
As a example it is possible to result a extreme system \ $SS$ \ with
abundance of multiple zeroes \ $\{(p_i-1) \cdot S(p_i)\}$ \ at \
$i \ \ge \ 0$\,.
\vspace*{3mm}

{\bf Theorem 24.} The obvious redundancy of a set of multiple zeroes \
$H_n^*$\,, \ that is density of zeroes \ $\gamma_n^*$ \ at creation
of majorizing estimation \ $msr_n (q)$ \ maximum series, returns the basic
characteristic to frequency of zeroes \ $\gamma_n$\,.
$$
\left\{H_n^* \gg H_n^{RR} > H_n; \ \gamma_n^* \gg \gamma_n^{RR} >
\gamma_n \Rightarrow MSR_n (q)\right\}: \quad
\left[ msr_n (q) \le \tau_n \frac{n + q}{1 - \gamma_n} + 1 \right].
\eqno{(7.6)}
$$
\vspace*{1mm}

{\sl Proof.} Even the fluent sight on advancing growth of a set of passive
(at formation of maximum series) multiple zeroes in relation to value \
$H_n^*$ \ specifies insecure of density of zeroes \ $\gamma_n^*$ \
as the candidate of basic value. The fact is that for any system valid
density \ $\alpha_n$\,, \ received with the help of \ $msr_n$ \ from the
formula of the main theorem, all further deviates from \ $\gamma_n^*$ \
with growth \ $n$\,.

In expression (7.6) effective multiple zeroes \ $H_n^{RR}$ \ and density \
$\gamma_n^{RR}$ \ can appreciably concede to number \ $H_n^*$ \ and
frequency \ $\gamma_n^*$\,. \ It restores support of density of zeroes
\ $\gamma_n$ \ in the formula of estimation of maximum series with some
factor \ $\tau_n$\,, \ it is possible, different from unit.
\hfill$\Box$
\vspace*{3mm}

It is necessary to take into account an alternate kind of algorithm of
filling for confirmation of a determining role of the characteristic
\ $n$-filling -- frequency \ $\gamma_n$\,.
\vspace*{3mm}

{\bf Theorem 25.} Supportness of density of zeroes \ $\gamma_n$ \ results to
to the majorizing coefficient two \ $(\tau = 2)$ \ of the main theorem
at frequency of zeroes \ $\gamma_n$\,. \ It reflects two-sided algorithm
of \ $n$-filling at formation of the upper estimation of maximum series \
$MSR_n(q)$ \ in a class of all non-singular systems \ $\{SS\}$\,.
$$
\left\{\gamma_n^* \gg \gamma_n^{RR} \Rightarrow \gamma_n \Rightarrow
MSR_n \left[q,\,Z_n^{(2)} \subset SS_{(III)} \right]
\right\}: \quad \left[ msr_n (q) \le 2 \frac{n + q}{1 - \gamma_n} + 1
\right]. \eqno{(7.7)}
$$
\vspace*{1mm}

{\sl Proof.} This major result is in detail discussed and submitted
above. Here the coefficient two is specified as a natural border,
predetermined by two-sided algorithm of filling. But the statement
cannot be perceived by the declaration of intentions. The two-sideness
of filling is destroyed by equivalent replacement of one grid
\ $S(a_{i+1})$ \ by two grids \ $S(2\,a_{i+1})$\,, \ each of which and
both are together the components of one-sided scheme. But the
one-sided algorithm of filling permits to construct a series of zeroes
twice smaller distance. From here and there is the coefficient two.

Argument about replacement of one grid \ $S(a_{i + 1})$ \ by two grids
\ $S(2\,a_{i + 1})$ \ should perceive by the only explanation of a
qualitative difference between one-sided and two-sided fillings algorithms.
Really, the continuation of such scheme results to replacement of a next
grid \ $S(a_{i + 2})$ \ by four grids \ $S(4\,a_{i + 2})$ \ with
preservation of former zeroes frequency \ $\gamma_n$\,. \ Thus one-sided
algorithm of filling is saved, but parameter of quantity of grids of
filling \ $n'$ \ appreciably grows. Further growth \ $n'$ \ results in
essential distortion of estimation \ $msr_{n'} \ < \ \frac{n'}{1 - \gamma_n}
+ 1$ \ of maximum series \ $MSR_{n'}$\,. \ Its value becomes more \
$2\,\frac{n}{1 - \gamma_n} + 1$ \ even at preservation of one-sided
fillings algorithm \ $Z_n^{(1)}$\,.

Let the inequality \ $a_{i + 1} < 2\,a_i$ \ takes place for all grids
of \ $n$-filling of described system \ $SS$ \ for \
$i = 1,\,2,\,...,\, n-1$\,. \ Then consecutive growth of
parameter \ $n$\,: \ $n = 2,\,3,\,... $ \ according to consequence 1 for
each stage permits find alternate one-sided and two-sided algorithms of
fillings. At the same time renewed for each step the replacement of the
grid \ $S(a_{i+1})$ \ by two grids \ $S(2\,a_{i+1})$ \ results to
estimation \ $msr_n\,: \ (i \ \Rightarrow \ 2i \ \Rightarrow \ 2n-1)$ \
of maximum series \ $MSR_n$
$$
\{S(a_{i})\}\,: \qquad \left\{S(a_{i+1}) \Rightarrow 2 \cdot S(2\,a_{i+1})
\right\}_1^{n-1} \ \Rightarrow \ \left\{\frac{2n-1}{1 - \gamma_n} + 1 \ < \
2\,\frac{n}{1 - \gamma_n} + 1\right\}\,.
$$

The coefficient \ $\tau = 2$ \ appears by limit (at \ $i = 1,\,2,...,\, n-1$)
\ of sequential reference to one-sided algorithm in a result of
replacement of the senior grid of current \ $(i+1)$-filling by two grids.

According to the previous theorem, excessive redundancy of zeroes density
\ $\gamma_n^*$ \ concerning real effective density \ $\gamma_n^{RR}$ \
conducts to supported frequency \ $\gamma_n$\,. \ For systems of the third
type \ $SS_{(III)}$ \ it means two-sided algorithm of filling \ $Z_n^{(2)}$
\ at formation of maximum series \ $MSR_n(q)$ \ in the expression (7.7).
In turn, such fact predetermines establishment of the majorizing
coefficient \ $(\tau = 2)$ \ in the formulation of the main theorem.
\hfill$\Box$
\vspace*{3mm}

Following statement specifies deep connections between imaging
principle, two-sided character of fillings algorithm and formulation
of the main theorem within the framework of fillings method for systems
of all types, among which it is necessary to allocate systems of the
third type \ $SS_{(III)}$\,.
\vspace*{3mm}

{\bf Theorem 26.} The imaging principle meets with unremoval redundancy
of multiple zeroes in systems the second and especially the third type
\ $SS_{(III)}$ \ and is transformed to two-sided algorithm of \
$n$-fillings, that results to majorizing coefficient \ $\tau = 2$
\ of the main theorem \ $MT$\,.
$$
\left\{MSR_n \left[q,\,\left(W_n^*\right)\right] \ \Leftrightarrow
\ MSR_n \left[q,\,Z_n^{(2)} \subset SS_{(III)} \right]
\right\}: \quad \left[ msr_n (q) \le 2 \frac{n + q}{1 - \gamma_n} + 1
\right]. \eqno{(7.8)}
$$
\vspace*{1mm}

{\sl Proof.} All results of the previous chapters and fillings method as
a whole specify validity of a put forward statement. It is thus
necessary to take into account majorizing character of related systems
of primes \ $SP_0\,, \ SP_1'\,, \ SP_1$ \ on mixing criterion. It
guarantees maximum of a constant \ $\tau$ \ in the formulation of the main
theorem. Such determined by two different ways value of the coefficient
is equal \ $\tau = 2$\,.

In the expression (7.8) process of creation of maximum series \
$MSR_n(q)$ \ with the help of imagings \ $(W_n^*)$ \ and then
account of their redundancy, results together with the introduction
of algorithm of two-sided filling, to occurrence of majorizing
coefficient \ $\tau = 2$ \ in the formula of the main theorem.
\hfill$\Box$
\vspace*{5mm}

\begin{center}{\Large \bf 8. The central theses of research and method}
\end{center}
\vspace*{3mm}

We shall move a intermediate result. We shall allocate and shall remind
the most important and turning points of the offered fillings method,
directly leading as to the formulation of the main theorem, as to its proof.
The determining character of consequences from it for number theory is
obvious.
\vspace*{3mm}

Let's restore turning-points of the fillings method immediately
carrying on to the statement of the main theorem and to the
proof of it.
\vspace*{3mm}

1. The major numerical characteristics of \ $n$-fillings and fillings method
directly depend from zeroes distribution of given grids and from
correlation (on base of divisibility) of grid's modules.

2. The source object is the strip region of binary elements. It
consists from \ $n$ \ of grids and has volume \ $(n;\,\infty)$\,.

3. The strip region is periodic, and the period of this strip is
equal \ $PZ_n$\,.

4. The main numerical characteristics of period of elements
strip region (quantity of units \ $E_n$ \ and zeroes \
$H_n^*$\,) are constants.

5. The imaging of the strip region on line inevitable erases constancy
of zeroes quantity (from \ $H_n^*$ \ up to \ $H_n$) \ on mapped
period and length of period, accordingly from \ $PZ_n^*$ \ up to
\ $PZ_n$\,.

6. The numerical relations of filling, including extremal properties,
accessible on some interval, depend from zeroes and units of imagined period.
The received estimations are corrected by length of an interval.

7. The imaging \ $ \stackrel{f(k)}{\Longrightarrow}$ \ of
multiple zeroes of the elements strip region are transformed in
the density of zeroes \ $(\gamma_n - \gamma_n^{**} - \gamma_n^*)$\,.

8. The precise values in expressions for the maximum series of zeroes
are found for fillings without multiple zeroes.

9. The upper-bound estimates of the maximum series depend from
density of zeroes in the interval of updated length.

10. The majorizing density of zeroes \ $\varrho_n$\,, \ not
exceeding \ $\gamma_n^*$\,, \ leads to an evaluation \
$msr_n(q) \ < \ (n + q)/(1 - \varrho_n) + 1$\,.

11. The proof of the main theorem is carried in classes of
systems and fillings - from degree-system up to systems of the third
type.

12. The precise values of all numerical characteristics of
fillings, including maximum series, are found for all class of
degree-systems.

13. The main theorem is proved for arbitrary systems of the
first and second type on the basis of outcomes for systems
without multiple zeroes.

14. The transition to systems of the third type is carried out
because of redundancy of the array of multiple zeroes.

15. The possibilities of the fillings method allow to construct
number of the independent proofs of the main theorem.

16. On the basis of the fillings method the proof of the main
theorem in three forms is created, including strong decisive
form without the support of the imaging principle with
majorizing multiplicative constant two for all classes of systems.

17. The direct application of the main theorem leads in
fillings and systems of prime numbers \ $SP_0$ \ and \ $SW_0$ \ to
outcomes noticeably majorizing reached earlier.

18. The constant two \ $(\tau = 2)$ \ of the main theorem \ $MT$ \
is confirmed by a principle of two-sided filling for all systems of the
third type, and also does not lean obviously on the imaging principle.

19. At the same time unremoval connection between the imaging
principle and two-sided scheme of filling the main theorem reveals.
\vspace*{3mm}

In the basis both fillings method and imaging principle
are incorporated rule of constancy of main numerical
characteristics of period and rule of adequacy of carry of
possibilities of an interval to the same period of filling.

The quantity of units remains by the absolutely invariable value
of mapped period. But imaging of the strip region can not save
simultaneously length of period and quantity of zeroes in systems
with multiple zeroes.
\vspace*{3mm}

Whole way lead to the proof of the main theorem, reasonably
convincingly testifies to independence of the fillings method
and originality of problem statement. Nevertheless it is
necessary precisely to separate offered from known sieving
process.
\vspace*{3mm}

Already comparison quite not neutral names of methods specifies
cardinal difference of initial points at formation of research
domain. If sieving process picks out primes rejecting other
numbers as rubbish, fillings method concentrates attention as
time on these rejected numbers, which are reflected by zeroes of
grids -- major elements of subsequent constructions.

Just zeroes, deleted elements, them multiplicity and frequency
in period of fillings appear by determining objects of further
structures and numerical analysis. This initial divergence in
determination of main object of fixed attention is ineradicable
reason of basic difference of methods.

There are all basis to consider, that fatal unimprovement of
many received classical estimations and it is enough easily
guessed proximity to unimproved many other, is explained
by absolute non-removal of basic value \ $x$ \ as bound of
sieving process algorithm. And really, any sieve by fixing of
\ $x$ \ concludes self in rigid cage.

Two compared methods are not reduced one to other, though at the
stage of initial model sieve is some fragment of filling. There are
especially far from fillings method of trigonometrical sums and
group of methods, connected with research of famous Riemannian
zeta-function \ $\zeta(s)$\,.

All numerous differences of sieving process and fillings method,
considered in monograph [1], result in the following conclusion:
\vspace*{3mm}

{\sl The fillings method in essence and even in private variants
do not reduce to the sieving process because of discrepancy of their
premises and initial objects.}
\vspace*{3mm}

\begin{center} {\Large \bf 9. Main achievements of the fillings
method} \end{center}
\vspace*{3mm}

All without exception described and proved statements lead to
the main theorem. The constants enter into their formulas.
They are different, but they are limited by several units. All received
results by fillings method follow qualitatively and sometimes even
quantitatively from first (weak) form of the main theorem with complete
imaging of zeroes. However given results of number theory strengthen
by unimpoved form of the main theorem with coefficient two.
\vspace*{4mm}

Fillings method is oriented to study of connection and dependence
of natural numbers as modules of grids in systems, however
its achievements, presented here, primes distributions concern.
As the first step we are addressed to the statement, antiquity
which permits to name its by mathematical symbol.
\vspace*{3mm}

{\sl Confirmation of statement.} The series of primes \ $p_i$ \ has
not completion.
\vspace*{3mm}

{\sl Third proof of the famous theorem} (after Euclid and
Euler), given in [1] is received not so much because from its
independence, necessity or special significance for further, how
many with purposes of demonstration of force and efficiency of the
fillings method. The address to primes system \ $SP_0$ \ demonstrates
increased set of units in interval \ $[p_n,\,p_p^2]$ \ of regulated
filling. They indicate primes above \ $p_n$\,.
\vspace*{3mm}

Concrete objects of the study usually lie on fixed intervals
regulated or semi-regulated fillings. Therefore all further
proofs base on statement, validity of which is clear from
definitions and whole fillings method:
\vspace*{3mm}

{\sl Values \ $msr_n (q)$ \ of maximum series \ $MSR_n (q)$ \ of unregulated
fillings majorize series values of other kinds of fillings algorithms.}
\vspace*{3mm}

{\bf A.} The maximum series in prime systems \ $SP_1$ \ and \ $SP_0$ \ are
$$
\begin{array}{l}
msr_n (SP_1) \ < \ 2\,p_{n - s}\,, \\
msr_n (SP_0) \ < \ 4\,p_{n - s - 1}\,, \\
\end{array} \ \qquad where \quad 1 \ < \ s \ = \ O(n)\,. \eqno{(9.1)}
$$
\vspace*{1mm}

{\sl Proof} follows from the main theorem and theorem 19, it is whence
possible to receive reasonably exact estimation for parameter \ $s$\,.
\hfill $\Box$
\vspace*{3mm}

{\bf B.} The distance between prime numbers satisfies to an
estimation:
$$
p_{n},\,p_{n-1} \ \in \ \{P\}: \quad \qquad p_{n} \ - \ p_{n-1} \
< \ C_n\,\sqrt{p_n}\,, \qquad C_n < 4\,. \eqno{(9.2)}
$$
\vspace*{3mm}

{\sl Proof.} Let some integer \ $N \in {\bf N}\,, \ N \ge 9$ \
is given. Proceeding from this value, we shall find quite
uniquely index of prime \ $p_m$ \ and corresponding interval \
$I_m = I_m(N)$ \ under the offered simple scheme:
$$
m \ = \ 1 \ + \ \max\limits_{p_i \leq \sqrt{N}} \ i\,, \qquad I_m \
= \ [\,p_{m-1} + 1\,, \  p_{m}^2 - 1\,]\,, \eqno{(9.3)}
$$
whence \ $m \ge 2$\,. \ We shall consider sieve of Eratosthenes
and corresponding regulated \ $m$-filling, generated by grids \
$S(p_0),\,S(p_1),\,...\,,\,S(p_{m-1})$\,. \ According to expression
(9.3) for given dependent intervals
$$
[\,p_{m-1} + 1\,, \ p_{m}^2 - 1\,] \ = \ I_m \ \ \supseteq \ \
I_N' \ = \ [\,p_{m-1} + 1\,, \ N - 1\,] \eqno{(9.4)}
$$
all units and only they correspond to primes. According to
received relations (9.3,\,9.4), upper bounds of intervals \ $I_m$ \ and
\ $I_N'$ \ are connected by inequality \ $p_{m}^2 \ \ge \ N$\,. \
It compels to address to the principal statement \ $MT$\,.

According to the main theorem greatest possible distance between
units (that is between primes) in interval \ $I_m$ \ always
can not exceed maximum series \ $msr_m$ \ of \ $m$-fillings
zeroes for system \ $SP_0$\,. \ Then we shall receive from
expressions (9.1) as estimates of maximum series \ $MSR_n (q)$ \
for all primes \ $p_n < p_m^2 - 1$\,:
$$
p_{n} \ - \ p_{n-1} \ < \ C_m\,p_m\,, \qquad m \ = \ 1 +
\max\limits_{p_i \leq \sqrt{N}} \ i\,, \ C_m < 4\,. \eqno{(9.5)}
$$

Now we shall distinguish responsible moment. At reception of
decisive conclusion there is no necessity to consider whole
interval \ $I_m \ = \ [\,p_{m-1} + 1\,, \  p_{m}^2 - 1\,]$\,, \
though just for it the inference is valid concerning maximum
series as possible distances between units. But regulation of \
$m$-filling (sieve) and way of construction of interval \ $I_m$
\ from (9.3) demonstrate, that for \ $N < p_{m-1}^2$ \ found
parameter \ $m$ \ decreases. Thus the constructive interval of
\ $N$ \ estimation, on which is searched of greatest distance
between primes, has bounds \ $p_{m-1}^2 + 1 \ \le \ N \ \le \
p_{m}^2 - 1$\,.

We shall take into account received conclusion, let \ $s = 0$ \
in expression (9.1) little coarsened estimation, and as far as with
sufficient precision
$$
p_m \ = \ \sqrt{N} \ + \ O(m)\,, \qquad N \ = \ p_n \ + \
o(N)\,, \eqno{(9.6)}
$$
we receive inequality (9.2).

The opportunity of decrease of coefficient \ $C_n < 4$ \ of
inequality (9.2) is problematic. Within the framework of filling
method and support on estimation of maximum series in filling
period (9.1) it is impossible, as far as according to theorem 14
limit (for \ $n \to \infty$) \ of coefficient \ $C_n$ \ is just 4.
Reserve of estimation decrease can consist in concrete
distribution of series in interval \ $J_m \ = \ [2,\,p_m^2]$\,.

From inequality (9.2) also follows, that between squares of
numbers \ $(M+2)^2$ \ and \ $M^2$ \ there is minimum one prime.
Clearly, this statement is valid and for squares of next primes.
Thus sequence of primes isn't limited. \hfill$\Box$
\vspace*{3mm}

{\bf C.} There is such integer \ $n_0 \gg 25$, that for
$n > n_0$ \ the inequality \ $p_{n+1} - p_n > \sqrt{p_n}$ \
is valid unlimited quantity in spite of Legendre's hypothesis.
\vspace*{3mm}

{\sl Proof.} According to theorem 14 for reasonably large \ $m \
(m \gg 26)$ \ inequality \ $msr_m(SP_0) \ \ge \ 4 \cdot p_{m -
2s - 2}$ \ where \ $s = O(m)$ \ is executed. Thus in period of \
$m$-filling  extensive file of series \ $SR_m$ \ meet without
fail, length of which lies in bounds \ $p_m < sr_m < C\,p_m$\,,
\ where \ $C = 4 - \varepsilon_m$\,. \ Already one such case of a
long series, fallen in initial interval \ $I_m$ \ (9.4,\,9.6)
of filling, refutes Legendre's hypothesis.

However the reasons, which have resulted to inequality (9.2) are
invalid at establishment of upper bound of maximum series \ $MSR_n$ \ in
interval \ $I_m$\,. \ So long series is not obliged to meet in interval \
$I_m$\,. \ At the same time series of practically such length
really lies in initial segment of extended interval \ $J_m \ = \
[2,\,N]$\,, \ but already regulated filling, instead of prime numbers.

This fact of regulated filling in \ $SP_0$ \ sharply limits
opportunity of achievement of series such length in interval \
$I_m$ \ with upper bound \ $p_m^2 - 1$ \ even for enough large \
$m$\,. \ For confirmation we shall notice, that according to
theorem 11 we have \ $msr_m \approx 2\cdot p_m$ \ for \ $m \sim
25$\,, \ but series of length \ $sr_n > p_m$ \ place on distance
about share of period from beginnings of fixed filling. And only
at inclusion of such series in following intervals \ $I_n$ \ it
is possible to expect qualitative change of estimations.

We shall consider major relations, describing condition of
Legendre's hypothesis. As basic Legendre's bound acts value \
$p_m \sim \sqrt{p_n}$\,:
$$
\frac{p_m - 1}{p_m}\,, \qquad \frac{p_{n+1} - p_n}{\sqrt{p_n}}\,,
\qquad \frac{p_m^2 - p_m}{p_m}\,, \qquad \frac{msr_m = C_m p_m
}{p_m}\,. \eqno{(9.7)}
$$

The first expression (9.7) of length of the first series of current \
$m$-filling to bound demonstrates firm proximity to unit. The
second expression gives checked values. Except the first infringements,
explicable by fillings method, for all \ $n > 8$ \ this value
less unit. The third expression demonstrates growth of length of
filling interval and search of infringement of hypothesis.
Special attention is deserved fourth value (9.7).

According to important theorem 19 expression of length of maximum series \
$msr_m$ \ to Legendre's bound already for small \ $m$ \ retires
from unit and up to \ $m \sim 26$ \ saves stable proximity to 2.
It is explained by role of the first axial configuration \
$Kf_m^{(I)}$ \ at generation of maximum series \ $MSR_n$\,. \ But
only at \ $m > 27$ \ qualitative transition to basic role of the
second axial configuration \ $Kf_m^{(II)}$ \ occurs, then last
expression (9.7) aspires to constant 4.

As far as checked values of the second expression (9.7) demonstrates
striking proximity to unit for some parameters \ $n$\,, \ for guaranteed
change of situation qualitative leap is necessary. Such leap
determines analysis of last correlation (9.7). It consists in
transition to the second axial configuration \ $Kf_m^{(II)}$ \
at generation of maximum series \ $MSR_n \ \to \ MSR_n (q)$\,.

In such case it is allowable to estimate order of those values,
since which it is possible to expect infringement of Legendre's
bound. As far as \ $p_{25}^2 > 10000$\,, \ it is necessary find
domain, near to length of period and \ $R$ \ as product of
primes
$$
R \ \sim \ \prod_{i=0}^{nn} \,p_i\,: \qquad
nn \ \sim \ \max \{j: \ p_j < 10000\} \ \sim \ 1085\,.
$$

Checks for such intervals \ $R$ \ of regulated filling for primes,
limited by value \ $p_{nn} \sim R$\,, \ are presented impossible
for observed time. At the same time and they are only beginning,
as coefficient \ $C_m$ \ in expression (9.7) begins essentially to surpass
insufficient value 2 only for parameters \ $m > 28$\,. \ Nevertheless,
fillings method immediately specifies the thesis, refutes
Legendre's hypothesis.

From here directly follows, that such series \ $SR_n$\,, \ superior
by length \ $sr_n$ \ initial series of filling, in result is met in
observed interval, that determines availability of pair of primes connected
by relation \ $p_{N+1} - p_N > \sqrt{p_N}$ \ at natural condition \
$p_{N+1},\,p_N \in I_{nn}$\,. \ Taking into account further increase
of coefficient \ $C_m \to 4$ \ in expressions (9.7), the
quantity of such cases is not limited. \hfill $\Box$
\vspace*{3mm}

{\bf D.} The sequence of twins \ $B_s$ \ has not completion on the
numerical axis.
\vspace*{3mm}

{\bf E.} The distance between twins satisfies to an evaluation
with the other constant:
$$
B_{s+1} \ - \ B_s \ < \ C_s\,\sqrt{B_s} \ \ln B_s\,, \quad (C_s < C\,, \
C > 1)\,. \eqno{(9.8)}
$$
\vspace*{3mm}

{\bf F.} The quantity of twins \ $\beta i\,(N)$ \ of smaller
value \ $N$ \ lays in boundaries \ \
$$
\frac{C'\,N}{\ln^2 N} \ < \ \beta i\,(N) \ < \
\frac{C''\,N}{\ln^2 N}\,, \qquad \ (C' < C'')\,. \eqno{(9.9)}
$$
\vspace*{3mm}

{ \bf G.} The infinite series, introduced by reverse values
of twins \ $B_s$\,, \ is convergent sequence.
$$
B_s \ \in \ \{B_s\} \ \subset \ \{P\}\,: \qquad \quad \sum_{s=1}^\infty \
\frac{1}{B_s} \ = \ C_{B} \ < \ \infty\,, \eqno{(9.10)}
$$
and value of constant \ $C_{B}$ \ accessible to numerical estimation.
\vspace*{3mm}

{\sl Proof} of connected statements \ {\bf D -- G} \ we shall
conduct in common. We shall consider \ $(2n+1)$-filling in
system \ $SW_0$ \ and frequency of zeroes \ $\gamma_{2n+1}$ \
for it:
$$
Z_{2n+1}: \qquad \gamma_{2n+1} \ = \ \frac{H_{2n+1}}{PZ_{2n+1}} \ = \ 1 \ - \
\frac{1}{2}\prod\limits_{i=1}^n \left( 1 - \frac{2}{p_i}
\right)\,. \eqno{(9.11)}
$$

Final estimation of maximum series for double prime system \
$SW_0$ \ according to the main theorem and by analogy to theorem
20 and formula (6.21) is
$$
MSR_{2n+1} (q): \qquad msr_{2n+1}(q) \ < \ C\,(2n + q +
1)\,\ln^2 n\,, \eqno{(9.12)}
$$
where value \ $C$ \ is limited and is reasonably well estimated
lower and upper with help of frequency from expression (9.11).
\vspace*{3mm}

Now fillings method permits proceed to practical estimation of
the numerical characteristics of double systems. At first we
shall consider kinsfolk of rank one.

Sieve for picking out of twins is regulated \ $(2n+1)$-filling in
system \ $SW_0$\,, \ received by product of two regulated \
$(n+1)$-fillings (sieves of Eratosthenes) in system \ $SP_0$\,,
\ shifted relatively one another by two positions (then grid \
$S(2)$ \ for them common). All units in interval \ $I_n = [p_n +
1\,, \ p_{n+1}^2 - 1]$ \ and only they correspond senior from
twins \ $B_s = p_k$\,.
\vspace*{3mm}

The statement \ {\bf D} \ immediately follows from here and from
expression (9.12) as far as \ $p_n \sim c\,n\,\ln n$ \ and even
for small \ $n$ \ in interval \ $I_n$ \ is always executed
$$
B_s \ \in \ I_n\,: \qquad p_{n+1}^2 - p_{n+1} - 2 \ > \ C\,(2n + q +
1)\,\ln^2 n \eqno{(9.13)}
$$
for some \ $q \ge 0$\,. \ At the same time any \ $q$ \ even
\ $q = 0$ \ means availability \ $B_s \in I_n$\,. \ Then
infinity of twins is consequence of relation (9.13) and infinity
primes. \hfill $\Box$
\vspace*{3mm}

The statement \ {\bf F} \ also follows from expressions (9.12,\,9.13).
Estimations of maximum series values \ $msr_{2n+1} (q,\,SW_0)$ \
permit to find bounds of \ $q = \beta i (N)$ \ for interval \
$I_n$ \ after equalization \ $p_{n+1}^2 = N$\,. \ Precision of
approximation is corrected by constants \ $C',\,C''$ \ from
correlation (9.9). \hfill $\Box$
\vspace*{3mm}

Statement \ {\bf G} \ and estimation \ $C_B$ \ from expression (9.10)
follow from expression (9.9) of statement {\bf F} as far as twins \
$B_s = c_s\,S\,\ln^2 S$\,, \ where factor \ $c_s$ \ is limited and
coefficient \ $c_s \to c$ \ at parameter \ $s \to \infty$\,. \hfill $\Box$
\vspace*{3mm}

Proof of statement \ {\bf E} \ basically repeats {\bf B}. For
estimation of greatest distance between next twins \ $B_{s+1}\,,
\ B_s$ \ it is necessary to apply to interval \ $I_n$\,. \ If
put \ $B_s \sim p_n^2$\,, \ the inequality (9.8) follows from
(9.12) for \ $q = 0$ \ as \ $\ln B_s \sim 2\,\ln n$\,. \hfill $\Box$
\vspace*{3mm}

{\bf H.} The distance between next Smith's numbers satisfies
to the following estimation \ $S_m - S_{m-1} \ < \ C\,\sqrt {S_m}\,\ln
\sqrt {S_m}$\,. \ Smith's number is prime \ $S_m = p_k$ \ and number
\ $0.5(p_k - 1) = p_l$ \ is prime too.
$$
S_m \ \in \ \{S_m\} \ \subset \ \{P\}\,: \qquad \quad S_m \ - \
S_{m-1} \ < \ C\,\sqrt{S_m} \ \ln\,S_m\,. \eqno{(9.14)}
$$
\vspace*{1mm}

{ \bf I.} The sequence of Smith's numbers \ $S_m $ \ is not
limited on an axis. The infinite series, introduced by
reverse Smith's numbers \ $S_m$\,, \ is convergent sequence.
\vspace*{3mm}

{\sl Proof} of related statements \ {\bf H} \ and \ {\bf I}
\ preferably to unite. The Smith's numbers important for many
appendices are offered and justified in [7].

That to pick out numbers of Smith it is necessary to generate
system \ $S2P$ \ as association of grid \ $S(2)$ \ and grids \
$S(2p_i)$ \ of double modules of system \ $SP_0$\,. \ In system
\ $S2P$ \ regulated \ $(n+1)$-filling repeats regulated \
$n$-filling of prime system \ $SP_0$ \ with doubled lengths of
all zeroes series. Product of two regulated fillings: \ $Z_n \in
SP_0$ \ and \ $Z_{n+1} \in S2P$\,, \ shifted by two positions,
determines \ $2n$-filling (initial grid \ $S(2)$ \ thus common)
of system \ $S3W$\,, \ for which units in interval \ $I_n$ \
correspond to Smith's numbers and only to them.

In result frequency of zeroes \ $\gamma_{2 n}$ \ of filling \
$Z_{2n} \subset S3W$ \ is equal
$$
S3W: \qquad \gamma_{2n} \ = \ \frac{H_{2n}}{PZ_{2n}} \ = \ 1 \ - \
\frac{1}{4}\prod\limits_{i=1}^{n-1} \left( 1 - \frac{2}{p_i}
\right)\,. \eqno{(9.15)}
$$

The comparison of expressions (9.11) and (9.15) permits to conclude,
that all inferences concerning Smith's numbers including expression
(9.14), can be transferred from statements \ {\bf D} -- {\bf G} \
and expressions (9.8--9.10) about twins distribution. \hfill $\Box$
\vspace*{3mm}

{\bf J.} (Goldbach's Conjecture). Even though one pair of prime
numbers \ $p_k\,, \ p_i$ \ (sometimes they can coincide)
satisfy to the equality \ $2J \ = \ p_k + p_i$ \ for each integer
\ $J \geq 2$\,, \ that is when \ $(J \in {\bf N})$\,.
\vspace*{3mm}

{\bf K.} The lower bound of Goldbach's different representations
\ $G(2J)$ \ for numbers \ $2J$ \ and \ $(J \in {\bf N})$ \ does not exist
at increase of argument \ $J$\,:
$$
\lim\limits_{J \to \infty} \ \ \inf\limits_{J \in [J_k,\,J_{k+1}]}
\ \ G(2J) \ = \ \infty\,; \qquad \forall\, (J_k < J_{k+1})\,. \eqno{(9.16)}
$$
\vspace*{3mm}

{\bf L.} The quantity of Goldbach's different representations
\ $G(2J)$ \ for the even value \ $2J \gg 12$ \ lays
in enough narrow limits for some constants:
$$
\frac{C'\,J}{\ln^2 J} \ < \ G\,(2J) \ < \ \frac{C''\,J\,\ln \ln
J}{\ln^2 J}\,, \qquad \ (C' < C'')\,. \eqno{(9.17)}
$$
\vspace*{3mm}

{\sl Proof} of statements \ {\bf J -- L} \ concerning
representation of any even number \ $2J$ \ as sum of two primes
\ $p_k + p_i$ \ is also better to unite.

We shall consider sieve of Eratosthenes as regulated variant of
\ $(n+1)$-filling with one difference: the first zeroes each
from \ $n$ \ odd grids corresponding to number (and module) \
$p_i\,, \ 1 \leq i \leq n$ \ are replaced by units. Then for
received structure, we shall designate which as \
$(n+1)^*$-filling, all units in interval \ $I_n \ = \ [2\,, \
p_{n+1}^2 - 1]$\,, \ and only they, are primes. For the offered
interval \ $I_n$ \ it is usual classical sieve of Eratosthenes.
This object is initial for next construction.

We shall consider units of product of two regulated fillings
(sieves of Eratosthenes): \ $(n+1)^*$-filling of system \ $SP_0$
\ and such \ $(n+1)^*$-filling of the same system, constructed
on the first filling, up to from value \ $2 J$ \ in opposite
direction. Units of such product and only units correspond to
Goldbach's representations for even numbers (concrete for number
\ $2 J \ \leq \ p_{n+1}^2 + 1$\,) \ in interval \ $[2\,, \ 2J -
2]$\,.

The quantity of various decompositions of Goldbach for number \
$2J$ \ is designated \ $G(2J)$\,. \ Now it is necessary find
unique suitable parameter \ $n$ \ of double filling:
$$
n \equiv n(2J) \ = \ \max\limits_{i \ \ge \ 1} \
\left\{i\,: \ \ p_i \ \leq \ \sqrt{2J - 2}\right\}, \quad
p_n^2 < 2J-1 \leq p_{n+1}^2. \eqno{(9.18)}
$$
\vspace*{3mm}

Statement \ {\bf J} \ (that is Goldbach's Conjecture in the
initial known form of uniqueness) follows from expressions
(7.12) and (7.18) as far as maximum series value appreciably
less \ $J$\,: \ $msr_{2n+1} < J$ \ or even \ $msr_{2n+1} \ll
J$\,. \ It means, that in fixed interval \ $ [2,\,J] $ \ units will be
always found out, which correspond to Goldbach's decompositions of
number \ $2J$\,. \ It proves the given statement. \hfill $\Box$
\vspace*{3mm}

Statement \ {\bf K} \ (Goldbach's Conjecture in strengthened
formulation) follows from advancing growth of parameter \ $J \sim p_n^2$ \
in relation to value of maximum series \ $msr_{2n+1} \sim
C_w\,p_n\,\ln n$\,. \ Thus quantity of decompositions \ $G(2J)$ \ grows.
\hfill $\Box$
\vspace*{3mm}

Statement \ {\bf L} \ in essential degree repeats the statement
\ {\bf F} \ for twins, but divisibility \ $J$ \ by primes
decreases frequency of zeroes and increases the upper bound of
quantity of Goldbach's decompositions \ $G(2J)$\,. \hfill $\Box$
\vspace*{3mm}

{\bf M.} The kinsfolk of all ranks \ $r\,, \ (1 \leq r < \infty)$ \
without exception there are on the numerical axis \ (that is integers \
${\bf N})$\,.
\vspace*{3mm}

{\bf N.} The distance between next kinsfolk of rank \ $r\,: \
(BR_s^r\,, \ BR_{s-1}^r$ \ are their senior prime numbers),
satisfies to an evaluation ($C^r$ \ is some constant):
$$
BR_s^r \ - \ BR_{s-1}^r \ < \ C^r\,\sqrt{BR_s^r}\,\ln
\sqrt{BR_s^r}\,, \qquad r \ge 1\,. \eqno{(9.19)}
$$
\vspace*{3mm}

{\bf O.} The sequence of kinsfolk of rank \ $r\,: \ \{BR_s^r\}_s$
\ is infinite (not completion) for all \ $r \ \geq \ 1$\,. \ For
example, twins.
\vspace*{3mm}

{\bf P.} The quantity of kinsfolk of rank \ $r$\,: \ $\{\beta\rho
i (N) \ = \ s\,, \ BR_s^r \ \leq \ N\} $ \ lays in boundaries
($C_1^r < C_2^r$ \ -- \ some constants):
$$
\frac{C_1^r\,N}{\ln^2 N} \ < \ \beta\rho i (N) \ < \
\frac{C_2^r\,N}{\ln^2 N}\,, \qquad r \ge 1\,. \eqno{(9.20)}
$$
\vspace*{3mm}

{\bf Q.} The series of numbers, introduced by reverse kinsfolk of rank
\ $r$\,, \ is convergent sequence. It includes and variant of twins.
$$
BR_s^r \ \in \ \{BR_s^r\} \ \subset \ \{P\}\,: \qquad \sum_{s=1}^\infty \
\frac{1}{BR_s^r} \ = \ C_{r} \ < \ \infty\,, \eqno{(9.21)}
$$
where value (constant) \ $C_{r}$ \ at concrete and not larger \ $r \ge 1$ \
quite accessible to numerical estimation. There is example of twins.
\vspace*{3mm}

{\sl Proof} of statements \ {\bf M -- Q} \ concerning distribution
kinsfolk of rank \ $r \ge 1$ \ (definition 28) in numerical axis
we shall present consistently.
\vspace*{3mm}

Statement \ {\bf M} \ follows from the analysis of generation of
configurations (definition 29) of type \ $Kf(i,\,r-i)$\,,
\ $r \ge 3\,, \ 1 \le i \le r - 1$ \ in systems \ $SP_0$ \ and \
$SP_1$ \ in period of fillings. Thus the existence \ $BR_s^r$ --
kinsfolk of rank \ $r$\,, \ that is configurations \ $Kf(r)$ \
for initial \ $r$ \ is known.

Already the first grid \ $S(2)$ \ results to generation of
important and typical configurations \ $Kf(1,\,2)$\,, \
$Kf(2,\,1)$\,. \ But not all configurations \ $Kf (i,\,r-i)$ \
are possible. For example, configurations \ $Kf(1,\,3v+1), \ v
\ge 0$ \ can not meet. However if kinsfolk of rank \ $r$ \ have
already arisen in period, according to divisibility there will
be configurations \ $Kf(i,\,r-i+1)$ \ though for some \ $i$\,.
\ After this configuration \ $Kf(r+1)$ \ will be generated by
next grid.

If to take into account, that some required configurations \
$Kf(r)$ \ are generated from configurations of kind \
$Kf(i,\,r-i)$\,, \ provided that quantity of configurations was
\ $K_n$ \ in period \ $PZ_n$\,, \ in period \ $PZ_{n+1}$ \ them
will become
$$
(p_{n+1} \ - \ 2)\,K_n \ \leq \ K_{n+1} \ \leq \ (p_{n+1} \ - \ 2 \ + \
\epsilon )\,K_n\,, \quad 0 \ < \ \epsilon \ \leq \ 1\,, \eqno{(9.22)}
$$
whence occurrence of values \ $BR_s^r$ (kinsfolk of rank \ $r$) \ in
interval of primes follows. There is variant of fillings and system
\ $SW_0$\,. \hfill $\Box$
\vspace*{3mm}

Statement \ {\bf N} \ repeats the statement \ {\bf E} \ for
twins. At the same time it is necessary to take into account,
that not all units of double filling (the second is shifted on \
$2r$), correspond kinsfolk of rank \ $r$\,. \ However capacity
of set of exceptions relatively small. \hfill $\Box$
\vspace*{3mm}

Statement \ {\bf O} \ repeats \ {\bf D} \ and follows from
expression (9.22). \hfill $\Box$
\vspace*{3mm}

Statement \ {\bf P} \ reminds the statement \ {\bf F} \ for twins,
but it should take into account the remark in {\bf N}. \hfill $\Box$
\vspace*{3mm}

Statement \ {\bf Q} \ similarly \ {\bf G} \ and follows from
{\bf P}. \hfill $\Box$
\vspace*{3mm}

{\bf R.} If given object \ $BK_s$ \ is the typical \ $s$-th configuration
of \ $m$ \ primes of length \ $M \ = \ 2\,\sum\limits_{i = 1}^{m-1}
\ r_i$\,, \ that the configurations \ $BK_s$ \ of primes lay in
boundaries, \ $(C' < C'')$ \ are some constants:
$$
C'\,s\,\ln^m s \ < \ BK_s \ < \ C''\,s\,\ln^m s\,, \qquad
s \ \gg \ M\,. \eqno{(9.23)}
$$
\vspace*{3mm}

{\sl Proof.} Picking out of configuration of \ $m-1$ \ ranks,
that is \ $m$ \ primes, requires product of \ $m$ \ fillings \
$Z_{n+1}$ \ of system \ $SP_0$\,, \ consistently shifted by \
$2r_i$ \ positions. So complex system including at \ $m$ \ grids
\ $S(p_i)$ \ is generated. For each \ $n$ \ under the formulas,
similar (9.11) it is possible find corresponding frequencies and
estimations, whence relation (9.23) follows immediately. \hfill $\Box$
\vspace*{3mm}

{\bf S.} The distance between next configurations of \ $m$ \ primes
($\,BK_s$\,, \ $BK_{s-1}$ \ are their senior prime numbers),
satisfies to an estimation, \ $C_b$ \ is some constant\,:
$$
BK_s \ - \ BK_{s-1} \ < \ C_b\,\sqrt{BK_s}\,\ln^{m-1}
\sqrt{BK_s}\,, \quad s \ \gg \ m\,. \eqno{(9.24)}
$$
\vspace*{3mm}

{\sl Proof.} Distance between configurations is estimated as well as
between twins \ ({\bf E}), if to take into account transition to
system with \ $m$ \ grids \ $S(p_i)$ \ and bounds (9.23). From here
inequality (9.24) follows, partial case of which relations
(9.2), (9.8) and (9.19) act. This result is in boundaries of
fillings method. \hfill $\Box$
\vspace*{3mm}

Purely, all indicated conclusions follow already from the first form of the
main theorem, however they are specified by third, strongest form. The
level of results achieved earlier with the help of sieving process and
some other methods can be in works of the number theory, for
example, in monograph [3].
\vspace*{5mm}

\begin{center}{\Large \bf 10. Graphic illustrations} \end{center}

\begin{center}
{\small (Diagrams don't insert in text because large volume. They can be
produce always)}
\end{center}
\vspace*{3mm}

All enclosed pictures and diagrams have especially illustrative purpose.
They do not claim for other role, and consequently it should not try on
their basis to make far conclusions. At the same time and such demonstrations
can put an idea into some change of a researched direction. In the similar
plan exhibition patterns of the fillings method can bring doubtless
advantage.

For this reason first of all the necessary explanation becomes obligatory
for creation regulated and unregulated fillings with occurrence of zeroes
series of a various length, including maximum. As initial filling we
shall graphically present complete period of \ $3$-filling in system \
$SP_1$\,, \ that is product of three grids with series of zeroes from
unit up to five.

Strict cycle of such filling \ $Z_3 \subset SP_1$ \ and independence
of arrangement of series \ $SR_3 (0)$ \ from the order of appearance
of grids in the period \ $PZ_3$ \ at product causes
coincidence regulated and unregulated fillings. It permits to present
the period on a circle, consistently having allocated for each series
sector, proportional its length. In such case of allocated series of a
length four \ $(sr_3 = 4)$ \ and five (maximum, \ $msr_3 = 5$\,)
will be till two specimens in the period.
\vspace*{12mm}

Fig. 1. Ring of 3-filling series in system \ $SP_1$\,.
\vspace*{3mm}

Grids \ $S(3)\,, \ S(5)$ \ and \ $S(7)$ \ form sole filling with
the united cyclic set of series on period \ $PZ_3$ \ of length \ $105$\,:
$$
PZ_3 \subset SP_1: \quad \ 512123132123313213234212124323123133212313212151
$$
This period includes first axis configuration of series \
$Kf^{(I)}_3 (1,\,5,\,1,\,5,\,1)$ \ and second axis configuration
\ $Kf^{(II)}_3 (4,\,2,\,1,\,2,\,1,\,2,\,4)$ \ of zeroes series. The
axis (symmetric) configurations are important for evaluation of maximum
series \ $MSR_n = MSR_n (0)$ \ for system \ $SP_1$ \ and other
systems.
\vspace*{12mm}

Fig. 2. The regulated \ $4$-filling in system \ $SP_0$\,.
\vspace*{3mm}

The grid \ $S(0) = L_0$ \ is line of units. \ $S(p_i)$ \ is grid of the
module \ $p_i$\,. \ $V(2)\,, \ V(3)\,, \ V(4)$ \ are products of
\ $2,\,3,\,4$ \ grids accordingly. \ $V(4) = Z_4$ \ for fixed
part of filling period. The multiplicity of zeroes in \ $Z_4$ \
is marked by additional dots. One dot signifies multiplicity of
two zeroes. The line \ $Z_4^*$ \ is conditional enlargement of
filling by multiplicity of zeroes. This line is given for
explanation of imaging principle.
\vspace*{12mm}

Fig. 3. Series of zeroes and maximum series \ $MSR_n$ \ in
system \ $SP_0$\,.
\vspace*{3mm}

The grid \ $S(p_i)$ \ is grid of the module \ $p_i$\,. \ The period is \
$PZ_5 = 2310$\,. \ The regulated \ $n$-fillings \ $Z_2 - Z_5$ \
are given and series \ $SR_5(0)\,, \ sr_5 (0) = 12$ \ (not
maximum) is distinguished in line \ $Z_5$\,. \ The first maximum
series \ $MSR_5^0 = MSR_5 (0)$ \ at regulated \ $5$-filling \
$Z_5$ \ is situated on interval under numbers \ $113 - 127$ \
and it is equal \ $msr_5^0 = 14$\,. \ The grid \ $S_w (p_i)$ \
has displacement \ $w$ \ of the first unit. Other displacements lead to
maximum series with one unit \ $MSR_5^1 = MSR_5 (1)$\,. \ The
length of this maximum series is value \ $msr_5^1 = 24$\,.
\vspace*{12mm}

Fig. 4. Maximum series \ $MSR_{11}$ \ in system \ $SP_1$\,.
\vspace*{3mm}

The grids \ $S_w (p_i)$ \ of the module \ $p_i$ \ have
displacements \ $w = w_i$ \ at product. Period of this filling
is \ $PZ_{11} \sim 3.71 \cdot 10^{12}$\,. \ The maximum series of
zeroes \ $MSR = MSR_{11}^{\,0} = MSR_{11} (0)$ \ has length \
$msr_{11}^0 = 33$\,. \ This series has 5 zeroes of two
multiplicity and such zeroes are marked by additional dots.
\vspace*{12mm}

Fig. 5. Maximum series \ $MSR_{12}$ \ in system \ $SP_1$\,.
\vspace*{3mm}

This illustration repeats scheme of Fig. 4.
The grids \ $S_w (p_i)$ \ of the module \ $p_i$ \ have
other displacements \ $w = w_i$ \ at product. Period of this filling
is \ $PZ_{12} \sim 1.52 \cdot 10^{14}$\,. \ The maximum series of
zeroes \ $MSR = MSR_{12}^{\,0}$ \ has length \ $msr_{12}^0 = 37$\,.
\ This series has 8 zeroes of two multiplicity and such zeroes
are marked by dots.
\vspace*{12mm}

Fig. 6. Twins, kinsfolk of rank two and numbers of Smith.
\vspace*{3mm}

There are grids \ $S(2),\,S(3),\,S(5),\,S(7)$ \ in lines 1 -- 4.
This grids form regulated \ $4$-filling in system \ $SP_0$\,. \
Sieve of Eratosthenes without elements \ $2,\,3,\,5,\,7$\,, \ that
is as zeroes, is presented by line 5.

The product of two sieves of Eratosthenes \ $Z_4^{**}$ \ with
displacement \ $w = 2$ \ is presented by lines 6 and 7. Twins as
corresponding units are marked by dark color. For example
(5,\,3), \ (7,\,5), ... , (61,\,59). This product is variant
of double prime \ $7$-filling of system \ $SW_0$\,.

The product of two sieves of Eratosthenes \ $Z_4^{**}$ \ with
displacement \ $w
= 4$ \ is presented by lines 8 and 9. Kinsfolk of rank two (primes
\ $p_i + 4 = p_{i+1}$\,) as corresponding units are marked by
dark color. For example (11,\,7), \ (17,\,13), ... , (47,\,43).
However primes (7,\,3) are not kinsfolk of rank two. This
product is variant of double prime \ $7$-filling also.

Smith's number is prime \ $S_m = p_k$ \ at
prime \ $0.5(p_k - 1) = p_r$\,. \ The product of sieve of
Eratosthenes (\,$Z_4^{**}$\,, \ line 10) and filling \ $ZZ_4^*$ \
(line 11) pick out numbers of Smith by corresponding units.
These units are marked by dark color. The filling \ $ZZ_4^*$ \
has elements: \ $l_{2s+1}(ZZ_4^*) = l_s (Z_4^{**})\,, \
l_{2s}(ZZ_4^*) = 0$\,. \ For example, numbers of Smith are 5
(units on places 5 and 2), 7 (7 and 3), ... , 59 (59 and 29).
This product of two prime fillings is mixed filling.
\vspace*{12mm}

Fig. 7. Twins and kinsfolk of rank two in intervals.
\vspace*{3mm}

The coincidence of number of series in the period does not mean such
coincidence twins and kinsfolk of rank two in an interval. A
little that, there can arise suspicion, that abundance of objects
of one kind necessarily forces out other objects, in particular if
the interval is not greater. To this question and next diagram is devoted.

Boundaries of interval \ $I_n$ \ are given by formula
\ $\rfloor\{1.2(n+1) \ln (n+1)\}^2\lfloor$ \ at \ $(n,\,n+1)$\,, \
where sign \ $\rfloor \cdot \lfloor$ \ means nearest
integer. Quantity \ $B_n$ \ of twins (twin) and kinsfolk of rank
two (primes \ $p_i + 4 = p_{i+1}$\,) as (qvad) in interval \
$I_n$ \ are presented by two broken lines. In spite of it, these
functions have strong correlation tie.
\vspace*{12mm}

Fig. 8. Increase of total quantity of kinsfolk of three ranks.
\vspace*{3mm}

Intervals \ $I_n$ \ are the same for Fig. 7. Total quantity \
$B_n$ \ of kinsfolk of rank 1 (twins, sumtw), rank 2 (sumqv) and
rank 3 (primes \ $p_i + 6 = p_{i+1}$\,) as (sumsx) for interval
\ $\sum_1^n I_k$ \ are presented. Quantity of kinsfolk of rank 3 is
essentially more than kinsfolk of rank 2 or 1\,.
\ Diagrams of number of twins and kinsfolk of rank 2
(sumtw, sumqv), \ not superior borders \ $N$\,, \ practically coincide
on all interval of supervision.

Here is submitted (curve, instead of graph) advancing growth of
number of kinsfolk of rank 3, reflected by function \ (sumsx)\,. \ It
is impossible not to note high smoothness of this function, monotone growing
concerning function of twins quantity  and kinsfolk of rank two, not
superior value \ $\sum R_m$\,.

\vspace*{12mm}

Fig. 9. Number of kinsfolk of three ranks on axis.
\vspace*{3mm}

At the same time the interval representation of common number
of kinsfolk not absolutely precisely reflects character and features
about growth of these functions. Therefore the following drawing offers three curve
distributions of the same objects depending on a growing border \ $N$\,, \
instead of from the number of the heaviest interval, included in area
\ $[1,\,N]$\,.

Quantity \ $B_n$ \ of kinsfolk of rank 1 (twins), rank 2 (primes
\ $p_i + 4 = p_{i+1}$\,) and rank 3 (primes \ $p_i + 6 =
p_{i+1}$\,) in interval \ $(1,\,N)$ \ are presented. Difference
with Fig. 8 is concluded in axis of abscissa.
\vspace*{12mm}

Fig. 10. Kinsfolk of rank 3 in intervals.
\vspace*{3mm}

There is difference between kinsfolk of rank 3 \ ($p_{k + 1} - p_k = 6$) \
with any such primes (optionally next). For fillings in systems \
$SP_0$ \ and \ $SW_0$ \ they mean, accordingly, series of a length
six and units of product of two identical fillings of a unary system
\ $SP_0\,: \ Z_{2n} (SW_0) = Z_{n + 1} \& Z_{n-1}$\,, \ shifted to six
points.

This diagram is analogy of Fig. 7. Quantity \ $B_n$ \ of
kinsfolk of rank 3 (primes \ $p_i + 6 = p_{i+1}$\,) in interval
\ $I_n$ \ is presented by broken line as (single). Function (double) is
unification and averaging of function (single): \
$B'_{(4n-1)/2} (double) = \frac{1}{2}\{ B_{2n-1}(single) +
B_{2n}(single) \}$\,.
\vspace*{12mm}

Fig. 11. Relation between quantity of objects of rank 3 and
number of twins.
\vspace*{3mm}

Certainly, with the same success it was possible to replace twins by
kinsfolk of rank 3 two. In the diagram two functions of aspiration to
various asymptotes of the relation of quantity of kinsfolk of rank
three \ (sumsx) \ and number of units of product of fillings \
$ZZ_n^*$ \ with shift six \ (sums6) \ to the same number of twins in
interval \ $[3,\,N]$ \ are submitted at \ $N < 50000$\,. \ One asymptote is
\ $y = 2$\,, \ a situation with second a little more difficult.

Functions \ sums6 \ and \ sumsx \ reflect relations between
primes:
$$
sums6 (N) \ = \ \frac{S6(N)}{TW(N)}\,, \qquad sumsx (N) \ = \
\frac{K6(N)}{TW(N)}\,, \qquad p_i \ \le \ N\,,
$$
where function \ $S6(N)$ \ is quantity of pairs primes: \ $p_i -
6 = p_k$ \ and \ $k = i-1$ \ or \ $k = i-2$\,, \ $p_i \le N$\,; \
function \ $K6(N)$ \ is quantity of pairs primes (kinsfolk of
rank 3): \ $p_i - 6 = p_{i-1}$\,; \ function \ $TW(N)$ \ is
quantity of pairs primes (kinsfolk of rank 1, twins): \ $p_i - 2
= p_{i-1}$\,. \ Asymptotes are \ $\lim\limits_{N \to \infty} \
sums6 (N) \ = \ 2\,, \ \lim\limits_{N \to \infty} \
sumsx (N) = 1.96683...$\,. \ However estimation of function \
$sumsx (N)$ \ gives value \ $\sim 1.56$ \ for \ $N \sim 50000$
\ or \ $n \sim 50$\,, \ since \ $p_{50}^2 > 50000$\,. \ Diagram
confirms this estimation.
\vspace*{12mm}

\hspace*{-2mm}Fig. 12. Goldbach's Conjecture and number of
decompositions $G(2J)$.
\vspace*{3mm}

Goldbach's Conjecture about representation of even number \ $2J$ \
in a kind of a sum two primes \ $2 J = p_k + p_m$ \ also a
subject of fixed attention of the fillings method. For a illustration
how by a specific image filling in the double prime system \ $SW_0$ \
acts we shall again consider initial regulated filling in system \ $SP_0$\,.

First 5 lines of this diagram are grids \
$S(2),\,S(3),\,S(5),\,S(7)$ \ in lines 1 -- 4, and they
form regulated \ $4$-filling \ $Z_4$ \ in system \ $SP_0$ \
(sieve of Eratosthenes without elements \ $2,\,3,\,5,\,7$\,; \
line 5). \ This 5 lines of diagram repeat lines of Fig. 6.

The product of two sieves of Eratosthenes \ $Z_4^{**}$ \ is
presented in lines 6 and 7. The mirror reflection of sieve from
point 60 is given in line 7. Then units and only units
correspond to decompositions of Goldbach for product of sieves in
interval \ $[2\,,\ 2J - 2]$\,, \ where \ $2J = 60$\,. \ The
red units correspond to decompositions of Goldbach and dark
units correspond to commutants of these decompositions.

The lines 8,\,9 and 10,\,11 are analogous to lines 6,\,7 with \
$2J = 58$ \ and \ $2J = 56$\,. \ Quantities \ $G(2J)$ \ of Goldbach's
decompositions are various: \ $G(60) = 6\,, \ G(58) = 4\,, \
G(56) = 3$\,. \ Only one decomposition \ $(58 = 29 + 29)$ \ at $2J =
58$ \ has not commutant on this diagram.
\vspace*{12mm}

Fig. 13. Number of Goldbach's decompositions in intervals \
$J_n$\,.
\vspace*{3mm}

The diagram reflects behavior of three functions, connected with
quantity of Goldbach's representations for each even number
\ $2J$ \ in an interval \ $J_n \ni 2J$\,. \ Heaviest, the average
and least values from number of decompositions in interval set
three growing curves.

Values \ (maxG), \ (meanG), \ (minG) \ depending from parameter of
the number of an interval \ $n$ \ are submitted here. According
to theoretical development, smoothness of maximum value of Goldbach's
decompositions \ $maxG$ \ appreciably below similar functions of
average and minimum value. It is explained specific divisibility of
value \ $2 Jmax$ \ by a number of primes, unique in each interval.
So, for example, \ $60060 = 3\cdot 4\cdot 5\cdot 7\cdot 11\cdot 13$\,.

Intervals \ $J_n$ \ must have even boundaries, therefore they
several differ from intervals \ $I_n$ \ of Fig. 7: \ $A_n =
2\cdot\rfloor 0.5 \{ 1.2(n+1) \ln (n+1)\}^2 \lfloor$\,, \ $J_n =
[A_{n-1} + 2\,, \ A_n]$\,. \ There are maximum of Goldbach's
decompositions \ $maxG = maxG_n = \max\limits_{2J
\in J_n} \ G(2J)$\,, \ minimum of Goldbach's decompositions \
$minG = minG_n = \min\limits_{2J \in J_n} \ G(2J)$ \ and mean \
$meanG$\,:
$$
meanG_n \ = \ \left\rfloor \frac{2}{A_n - A_{n-1}}
\sum_{2J \in J_n} \ G(2J)\right\lfloor, \quad minG < meanG < maxG
$$
for \ $n > 1$\,. \ Values are for \ $n = 1$\,: \ $meanG_1 =
maxG_1 = 2\,, \ minG_1 = 1$\,. \ Functions \ $minG_n$ \ and \
$meanG_n$ \ are monotonically increasing sequences for all \ $1
\le n \le 50$\,, \ and smoothness of function \ $meanG_n$ \ is
essentially higher, than functions \ $minG_n$ \ and \ $maxG_n$\,.
\vspace*{12mm}

Fig. 14. Extremal and mean data of number of Goldbach's
decompositions in segments of axis \ $(2J \in J_n)$\,.
\vspace*{3mm}

This diagram repeats scheme of Fig. 13 \ with the exception of
abscissa axis. Arguments \ $2J$ \ of functions \ $minG (2J)$ \
and \ $maxG (2J)$ \ are found and belong to intervals \ $J_n \ni
2J$\,. \ Argument \ $2J(n)$ \ of function \ $meanG_n [2J(n)]$ \ is \
$2J(n) = 0.5\,[A_{n-1} + A_n]$\,. \ Linearity of functions
especially \ $meanG$ \ clearly visible and justifies choice of
intervals \ $J_n$\,.
\vspace*{12mm}

Fig. 15. Approximation of minimum quantity of Goldbach's
decompositions in intervals \ $J_n$ \ by normalizing function.
\vspace*{3mm}

The diagram reflects perfect approximation of all three considered functions.
Broken lines \ {\bf maxG\,, \ meanG} \ and \ {\bf minG} \ on this diagram
signify correlations of functions from Fig. 13 and 14 to one and the
same increasing function \ $2J/\ln^2 (2J)$\,. \ In particular
importance has approximation
$$
{\rm \bf minG} \ = \ \frac{minG_n(2J)}{2J}\, \ln^2 (2J)\,, \qquad
2J \ \in \ J_n\,.
$$
The others functions have secondary significance. These broken
lines have asymptotes: line \ $y = \ \sim 0.74$ \ for main
function \ {\bf minG}\,, \ $y = \sim 1.22$ \ for function {\bf
meanG} \ and \ $y = \ \sim 3$ \ for function {\bf maxG}\,. \ There are
empirical data, are confirmed by the theory.
\vspace*{3mm}

All given illustrations are applications of the fillings method
and do not contradict to these conclusions.

\vspace*{3mm}

\begin{center}{\Large \bf 11. Conclusion} \end{center}
\vspace*{3mm}

After whole stated it follows to accent attention to sources,
given such appreciable advantage of developed method before
sieving process and some other methods having application in
number theory. This superiority can not be explained in essence,
if initially to assume the fillings method affiliated.

The fillings method, described in [1] and briefly in [2],
variant of which is offered here, it is far before completion
still. It is faster allowable to say about initial stage of
research. At the same time and considered states, leaning on
main theorem, permit to answer rather wide range of problems of
primes distribution, is not solved by other methods, including
various modifications of sieving process, exhaustion which less
and less doubts causes.

Results of the fillings method already have found practical application at
creation of algorithms of randomness simulating [6,7]. A little that
just use of a sequence of Smith's numbers has allowed to generate
algorithms, in a limit ensuring random simulated sequences [8].
Research of theoretic-numerical congruent pseudorandom numbers generator
has resulted in a little unexpected conclusion about necessity of
modernization of ancient Euclid's algorithm [9].

It is possible interesting, how there were ideas about other
approach to problem. The fact is that sieving of numbers,
multiple to some, should be optionally take place to multiple
numbers. Only linear connection is quite necessary, but in such
variant concept of grid and its shift is generated.

There was representation about zeroes series, greatest series
and its determining role in forthcoming appendices. All this
appears in overlapping strips with put grids, that is in product
of grids. So purely practically was born concept of filling in
unregulated and regulated form. Naturally, from here there was
equally step before picking out classes of systems.

The fillings method has become to crystallize little late [2],
when the exact formulas for degree systems, then exact
estimations for systems without multiple zero and some other
were received. Then and there was idea about imaging of strip
region of multiple zeroes to extended period and imaging leads
to variant of filling without multiple zeroes. Qualitative might
of such idea was immediately confirmed by all accessible
examples.

Particularly all this has found reflection in formulation of the
main theorem of fillings method. Its validity and logic, clear
from preamble, all have passed through numerous stage of proof,
including even not leaning on main idea. Such increased
attention to the central theorem is explained by that just here
unique turning point of fillings method lies, which offers to
accept idea, and not just naked mathematical transformation.

It should note, that at reception of main conclusions,
concerning primes distribution and corresponding numerical
objects, only classical theorem about primes is used. And the
necessity arises only in it is possible to sharper estimation of
sum reverse primes for finite and fixed value \ $n$\,. \
Summation permits sharply to specify integrated
estimations of corresponding constants, which in turn find
reflection in all received expressions.

In either case fillings method does not use which was the
obliging assumptions or doubtful postulates at one of stages.
Its sharp mathematization is based on finiteness of all initial
objects and clear logic of finitely observed conclusions. In
this method essentially differs from other extreme formalized
methods, for example inevitably leaning on such concept as "set
of all primes".

For this reason rather the plenty of problems, \ enabling \ to \ hope \
on \ progress in solving, always it is not limited only to
classical problems of primes distribution. The fillings method
admits research of any integer objects, as elements of
constructed system. At the same time arbitrary modules vector of
filling grids can not be abstract from their coprime
characteristics, that it is reflected in next values of zeroes
frequency. From here follows, that aspects of primes
distribution will be always exclusively important and for the
fillings method.

Achieved results, enumeration of which is not limited by
classical problems of primes distribution, specifies high
efficiency of the fillings method. Thus it is impossible to
underestimate significance of the main theorem as base of all
research. The fillings method or method of imaging of
frequent functions in period to interval within the framework
of system of generated grids permits in another way to present
many problems of the number theory.
\vspace*{3mm}

The number of considered problems has received due base support, published
only in monograph [4]. On this basis insufficiency of mathematically-logic
constructions was proven at creation of defended statements [10].
In particular, it concerns so of not clear (imaginary) objects as
infinite or continuum set.
\vspace*{10mm}

{\bf Reference}
\vspace*{5mm}

{\normalsize
 [1]  Antipov M.V. The Fillings Method and Problems of Prime Numbers
Distribution. -- Novosibirsk, -- Printing-house of Siberian
Stateservice Academy, 2002, 503 pp.

 [2]  Antipov M.V. The fillings method and some problems of number
theory. -- Novosibirsk, 1984. -- 19 pp. -- (Preprint / AN USSR,
Sib. Br., Comp. Cent.; N 528).

 [3]  Prachar K. Primzahlverteilung. Wien, Springer - Verlag, 1957, 512 pp.

 [4]  Antipov M.V. The Restriction Principle. -- Novosibirsk, --
Printing-house of Siberian Branch of Russian Academy of
Sciences, 1998, 444 pp.

 [5]  Antipov M.V. The restriction principle and foundation of
mathematics. -- Novosibirsk, 1997. -- 112 pp. -- (Preprint / RAN,
Siberian Branch, Inst. of Comp. Math. and Math. Geoph.; N 1100).

 [6]  Antipov M.V. Reality and Pseudorandomness. -- Novosibirsk,
Manuscript of monograph, 1992, 420 pp.

 [7]  Antipov M.V. Sequences of Numbers for the Monte Carlo Methods
// Monte Carlo Methods and Appl., Vol. 2, N 3, pp. 219 -- 236 (1996)
VSP, Utrecht, Tokyo.

 [8]  Antipov M.V. Congruent Operator Simulation
of Continuous Distributions // Computational Mathematics and Mathematical
Physics, Vol. 42, {\bf N} 11, 2002, pp. 1572 -- 1580.

 [9]  Antipov M.V. Congruence operator of the pseudo-random
numbers generator and a modification of Euclidean decomposition //
Monte Carlo Methods and Appl., Vol. 1, {\bf N} 3, pp. 203 -- 219
(1995), Utrect, Tokyo.

 [10]  Antipov M.V. Mirages of Evidence. -- Novosibirsk, -- OOO "Omega
Print", 2006, 120 pp.}
\vspace*{5mm}

\begin{center}
{\Large \bf Contents}
\end{center}
\vspace*{3mm}

1. Introduction \hfill2

2. Some definitions \hfill4

3. Systems without multiple zeroes \hfill8

4. Imaging principle and main theorem \hfill15

5. Premises of an evidence of the main theorem \hfill21

6. Proof of the main theorem \hfill28

7. Two-sided fillings and main theorem \hfill47

8. The central theses of research and method \hfill54

9. Main achievements of the fillings method \hfill56

10. Graphic illustrations \hfill64

11. Conclusion \hfill70

Reference \hfill72

\vspace*{3mm}

Lavrentieva st. 6, Novosibirsk, 630090, Russia

E-mail: \ amv@osmf.sscc.ru

http://osmf.sscc.ru/$^{\sim}$amv
\vspace*{5mm}
\end{document}